\documentclass[preprint,12pt]{elsarticle}

\usepackage{amsfonts,amsmath,amsthm,amssymb}
\usepackage{fullpage,lineno,setspace, nicefrac, subcaption, caption}
\usepackage{mathrsfs}
\usepackage{tikz}
\usetikzlibrary{automata,arrows,positioning,calc}
\usepackage[mathscr]{eucal}
\usepackage{xcolor}
\usepackage{hyperref}

\DeclareGraphicsExtensions{.png,.jpeg,.jpg,.pdf}
\graphicspath{{./Figures/}}
\allowdisplaybreaks

\theoremstyle{plain}
\newtheorem{theorem}{Theorem}[section]
\newtheorem{definition}{Definition}[section]
\newtheorem{proposition}[theorem]{Proposition}
\newtheorem{lemma}[theorem]{Lemma}

\newtheorem{remark}[theorem]{Remark}

\newtheorem{notation}{Notation}

\theoremstyle{definition}
\newtheorem{case}{Case}[section]

\usepackage{xcolor}
\definecolor{darkgreen}{rgb}{0,0.5,0}
\usepackage[colorinlistoftodos]{todonotes}
\usepackage{soul}
\setulcolor{red} 
\setstcolor{red} 
\sethlcolor{yellow} 

\newcommand{\lb}{\left\{}
\newcommand{\rb}{\right\}}
\newcommand{\set}[1]{\left\{ #1 \right\}}

\newcommand{\Def}{\overset{\textbf{def}}{=}}
\newcommand{\RR}{\mathbb{R}}
\newcommand{\ZZ}{\mathbb{Z}}
\newcommand{\NN}{\mathbb{N}}

\newcommand{\eps}{\varepsilon}
\newcommand{\bOne}{\textcolor{black}{\mathbf{1}}}

\newcommand{\rhomax}{\rho_{\text{max}}}
\newcommand{\fpartial}[2]{\frac{\partial #1}{\partial #2}}
\newcommand{\diverg}[1]{\nabla \cdot #1}

\newcommand{\abs}[1]{\lvert #1 \rvert}
\newcommand{\norm}[1]{\left\| #1 \right\|}

\newcommand{\Riemann}{\mathcal R}

\newcommand{\sign}{\text{sign}}

\newcommand{\front}{\mathfrak D}
\newcommand{\temple}{\mathcal T}

\newcommand{\bigmid}{\bigl\vert}
\newcommand{\action}[2]{\langle #1, #2 \rangle}

\newcommand{\tup}[1]{\textup{(}#1\textup{)}}
\usepackage{amsopn}
\newcommand{\totvar}[2]{\underset{#2}{T.V.}\left(#1 \right)}
\newcommand{\posvar}[2]{\underset{#2}{P.V.}\left(#1 \right)}
\newcommand{\BV}{\textbf{BV}}
\DeclareMathOperator{\supp}{Supp}
\DeclareMathOperator{\ball}{Ball}
\DeclareMathOperator{\loc}{loc}

\newcommand{\sa}{\mathsf a}
\newcommand{\ft}{\mathfrak t}
\newcommand{\AC}{\mathcal A_{\mathfrak c}}

\newcommand{\zplus}{z^+}

\newcommand{\zminus}{z^-}
\newcommand{\hatBN}{\hat{\mathcal A}}
\newcommand{\checkBN}{\check{\mathcal A}}

\newcommand{\nnorm}[1]{{\left\vert\kern-0.25ex\left\vert\kern-0.25ex\left\vert #1 
    \right\vert\kern-0.25ex\right\vert\kern-0.25ex\right\vert}}

\journal{ArXiv}

\begin{document}

\begin{frontmatter}



\title{On the Existence of Solution of Conservation Law with Moving Bottleneck and Discontinuity in FLux}


\author{Hossein Nick Zinat Matin}

\affiliation{organization={University of California at Berkeley},
            addressline={}, 
            city={Berkeley},
            postcode={94720}, 
            state={CA},
            country={USA}}
            
\author{Maria Laura Delle Monache}

\affiliation{organization={University of California at Berkeley},
            addressline={}, 
            city={Berkeley},
            postcode={94720}, 
            state={CA},
            country={USA}}
\begin{abstract}
In this paper, a PDE-ODE model with discontinuity in the flux as well as a flux constraint is analyzed. A modified Riemann solution is proposed and the existence of a weak solution to the Cauchy problem is rigorously investigated using the wavefront tracking scheme. 
\end{abstract}



\begin{keyword}
Conservation Law, Traffic Flow, PDE-ODE Model, Discontinuous Flux, Moving Bottleneck Dynamic in Traffic Flow



AMS Classification: 35L65
\end{keyword}

\end{frontmatter}


\section{Introduction and Related Works}
\label{S:Dynamical_Model}
This paper is concerned with conservation law consisting of discontinuous flux and flux constraint. Such a dynamical model, which is presented by a coupled partial and ordinary differential equation (PDE-ODE), is important in various applications in engineering and physics. The main motivation for the present work stems from traffic flow dynamics in the presence of moving bottlenecks. From an application point of view, Autonomous vehicles (AVs), which are considered one of the most promising technologies for regulating the traffic flow condition \cite{wang2016multiple, milakis2017policy}, can be naturally considered as moving bottlenecks and hence their influence on traffic behavior can be studied through understanding the properties of the solution of a PDE-ODE model \cite{lattanzio2011moving, liard2020optimal}. Such consideration can be employed to control and improve the traffic conditions, e.g. fuel consumption or total travel time \cite{piacentini2020traffic, liard2020optimal, daini2022centralized, piacentini2019highway, wadud2016help}. In addition, various numerical methods have been proposed to approximate the solution of this model \cite{daganzo2005moving, delle2016numerical, daganzo2005numerical, chalons2018conservative, simoni2017fast}.
More precisely, the behavior of the traffic dynamic is strongly affected by the dynamic of the bottleneck. In this sense, moving bottlenecks can create congestion in the upstream traffic and change the pattern of the corresponding flow. The interaction of the bottleneck with traffic has been investigated both from engineering perspectives \cite{lebacque1998introducing, leclercq2004moving, giorgi2002traffic} and mathematical point of view \cite{lattanzio2011moving, delle2014scalar, delle2017stability, liard2019well, liard2021entropic}. 

From a theoretical point of view, in \cite{delle2014scalar} the presence of a moving bottleneck and the way it affects the traffic flow dynamics is modeled by a PDE-ODE model. In particular, the behavior of the bottleneck is captured by incorporating an ODE, defining the trajectory of the bottleneck, and a constraint that presents the capacity reduction of the road to the Lighthill-Whitham-Richards (LWR) Conservation Law \cite{lighthill1955kinematic, richards1956shock}. The convergence result is extended in \cite{liard2021entropic}. Goatin et al. \cite{goatin2022interacting} studies the interaction of several bottlenecks. The existence of the solution in the case of a time-varying desired speed of the bottleneck is shown in \cite{garavello2020multiscale}. The extension of LWR to the (second-order) Aw-Rascle (AR) model for the bottleneck problem is addressed in \cite{villa2017moving}.
The stability of the solution of the PDE-ODE model is discussed in \cite{delle2017stability} and is generalized rigorously as a well-posedness problem in \cite{liard2019well}.

The problem of discontinuity in the flux attracted many research papers due to its applications in more complex systems; e.g. variable maximal speed limits in traffic dynamics. In the presence of discontinuity, on the other hand, the system is resonant (coinciding eigenvalues) and hence non-strictly hyperbolic. The main consideration in all different methods is treating the discontinuity points such that existence and uniqueness can be deduced; for instance, the flux function is the same on both sides of the discontinuity points which is basically the Rankine-Hugoniot condition. However, this is not usually sufficient to prove the uniqueness and stronger conditions need to be imposed. We review some of the approaches that have been considered in the literature. This problem was first addressed in \cite{gimse1992solution}.
Diehl \cite{diehl1995scalar}, Burger et al. \cite{ burger2009engquist} and Garavello et al. \cite{garavello2007conservation} consider the flux function 
\begin{equation}\label{E:flux_special} 
f(x, \rho) = H(x) f_L(\rho) + (1 - H(x)) f_R(\rho))
\end{equation}
with a finite number of discontinuities, $f_L$ and $f_R$ are the values of the flux function on the right and left sides of the jump discontinuity points, and $H$ is the Heaviside function. Diehl \cite{diehl1995scalar} introduces a condition for the uniqueness of the solution using Riemann solutions. 
 
The convergence of a difference scheme to the entropy solution is proven in \cite{burger2009engquist}.
Temple \cite{temple1982global} considers a $2 \times 2$ system of conservation laws with discontinuity in the flux with application in polymer and oil. Introducing a Riemann solution, the convergence of the Glimm difference scheme to the solution of the Cauchy problem is shown by introducing a singular map. The author shows that the total variation and compactness can be applied to the solution in the new space and consequently the convergence for the main approximate solution then follows by the properties of the map. 

The case of $f(x, \rho) = \gamma(x) f(\rho)$ for strictly concave $f$ and a function $\gamma$ is investigated in the work of \cite{seguin2003analysis, towers2000convergence}. Tower \cite{towers2000convergence} uses the Godunov and Engquist-Osher (EO) flux and a singularity map similar to that of \cite{temple1982global} to show the convergence of the approximate solutions. In addition, the entropy condition for the case of $\gamma$ piecewise $C^1$ with a finite number of discontinuities is obtained. Adimurthi et al. \cite{adimurthi2005optimal} consider the flux function of the form \eqref{E:flux_special} and they show the $L^1$-contractive solution. They also show the a.e. convergence of the Godunov approximation to the solution of the conservation law for a particular set of discontinuities of the flux.

In this paper, we consider a conservation law with a discontinuity in the flux and flux constraint. To the best of our knowledge, this is the first exposition of such a problem. In a sense, we are generalizing the result of \cite{delle2014scalar} by considering the jump discontinuity in the flux. Such a generalization is motivated by some applications in traffic flow and in particular by characteristics of traffic dynamics in different regions of the road; e.g. different speed limits in each segment of the road. In this paper, we will follow the wavefront tracking scheme and use a similar approach to \cite{temple1982global} we show the existence of the solution to the Cauchy problem using the total variation of a homeomorphism. The relative complexity of such an approach, which will be discussed in detail later, stems from the convolution of classical and non-classical waves as well as jump discontinuities in the flux. Furthermore, the flux discontinuity generally results in a loss of uniform total variation and resonant systems and hence showing the convergence of the approximate solutions in this method requires more work.
\subsection{Dynamical Model}
From the mathematical point of view, we are concerned with analyzing the existence of a solution to the following dynamical model of traffic for $(t, x) \in \RR_+ \times \RR,$
\begin{equation}\label{E:system}
\begin{cases}
    \rho_t + \fpartial{}{x}[f(\gamma,\rho)]= 0   \\
   \rho(0, x) = \rho_\circ(x) \\
    f(\gamma(y(t)), \rho(t, y(t)) - \dot y(t) \rho(t, y(t)) \le F_\alpha(y(t), \dot y(t)) \Def \frac{\alpha \rhomax}{4 \gamma(y(t))} \left(\gamma(y(t)) - \dot y(t) \right)^2\\
    \dot y(t) = \omega(y(t), \rho(t, y(t)+))\\
    y(0) = y_\circ.
\end{cases}\end{equation}
where function $x \mapsto \gamma(x)$ is considered to be a piecewise constant function 
\begin{equation}\label{E:gamma_function}
    \gamma(x) = \begin{cases} \gamma_{r_\circ} &, x \in I_0 \\
    \gamma_{r_1} &, x \in I_1 \\
    \vdots\\
    \gamma_{r_M} &, x \in I_M \end{cases}
\end{equation} 
Here, $I_0 = (- \infty, \sa_1)$, $I_m = [\sa_m, \sa_{m+1})$ for $m = 1, \cdots, M$ with $\sa_{M+1} = \infty$ (right continuous function). Function $\gamma(x)$ can be interpreted as the \textit{speed limit} at location $x \in \RR$ which is constant on each interval $I_m$. 

The unknown function $(t, x) \in \RR_+\times \RR \mapsto \rho(t, x)\in [0, \rhomax]$ denotes the density function. The flux function $f:\set{\gamma_{r_\circ}, \cdots, \gamma_{r_M}}  \times [0, \rhomax] \to \RR_+$ is defined by
\begin{equation}\label{E:flux}
f(\gamma, \rho) \Def  \rho v(\gamma, \rho)
\end{equation}
where the mean traffic speed $v:\set{\gamma_{r_\circ}, \cdots, \gamma_{r_M}} \times [0,\rhomax] \to \RR_+$ is defined by
\begin{equation}
    \label{E:velocity}
    v(\gamma, \rho) \Def  \gamma \left(1 - \frac{\rho}{\rhomax} \right).
\end{equation}
It should be noted that in equation \eqref{E:velocity},  $v(\gamma, \rhomax) = 0$ and $v(\gamma, 0) = \gamma$ (the maximal speed).

To understand the contribution of the bottleneck in the dynamic of the traffic flow, let function $t \mapsto y(t)$ denote the trajectory of the moving bottleneck (slow-moving vehicle such as a controlled autonomous vehicle). Function $y\in \RR \mapsto V_b(y)\in \RR_+$ denotes the maximal speed of the moving bottleneck at the location $y \in \RR$. Following \eqref{E:gamma_function} we consider function $V_b(\cdot)$ to be a piecewise constant function on each interval $I_m ,\, m = 0 , \cdots, M$, and we denote it by
\begin{equation*} 
V_b^{(m)} \Def V_b(y), \quad  y \in I_m, \quad \forall m = 0, \cdots, M. \end{equation*}
Whenever the traffic condition allows, the bottleneck moves with its maximum velocity which satisfies
\begin{equation*}
    V_b^{(m)} < \gamma_{r_m}, \quad  m = 0 , \cdots, M.
\end{equation*}
On the other hand, when the surrounding density is more than a threshold, the velocity of the bottleneck would be adapted accordingly. More accurately, we define the bottleneck speed profile by
\begin{equation}
    \omega(y, \rho) \Def \min \set{V_b(y) , v(\gamma(y), \rho)}.
\end{equation}
where, $v(\gamma(y), \rho)$ is defined as in \eqref{E:velocity}. It can be noted that for $\rho \le \rho^*(y) \Def \rhomax \left(1 - \frac{V_b(y)}{\gamma(y)} \right)$ we have $V_b(y) \le v(\gamma(y), \rho)$ and the converse inequality holds for $\rho > \rho^*(y)$ and $\rho^*$ is a piecewise constant function of the form
\begin{equation*}
    \rho^*_m \Def \rho^*(y), \quad y \in I_m, \, m \in \set{0, \cdots, M}.
\end{equation*}
The bottleneck trajectory $t \mapsto y(t)$ follows from the dynamics of the form 
\begin{equation}\label{E:BN_trajectory}\begin{split}
    \dot y(t) &= \omega(y(t), \rho(t, y(t)+)), \quad y(0)= y_\circ\\
    &  =  \min \set{V_b(y(t)) ,\gamma(y(t)) (1 - \rho(t, y(t)+))}
    \end{split}
\end{equation}
The functions $V_b(\cdot)$ and $\gamma(\cdot)$ are piecewise constant. The solution to \eqref{E:BN_trajectory} is well-posed in Carath\'eodory sense (see \cite{delle2014scalar, colombo2003holder}); i.e. an absolutely continuous function $y$ which satisfies
\begin{equation*}
    y(t) = y_\circ + \int_0^t w(y(s), \rho(s, y(s)+) ds.
\end{equation*}
\subsection{The Dynamics of Moving Bottleneck}
\begin{figure}
    \centering
    \includegraphics[width=3in]{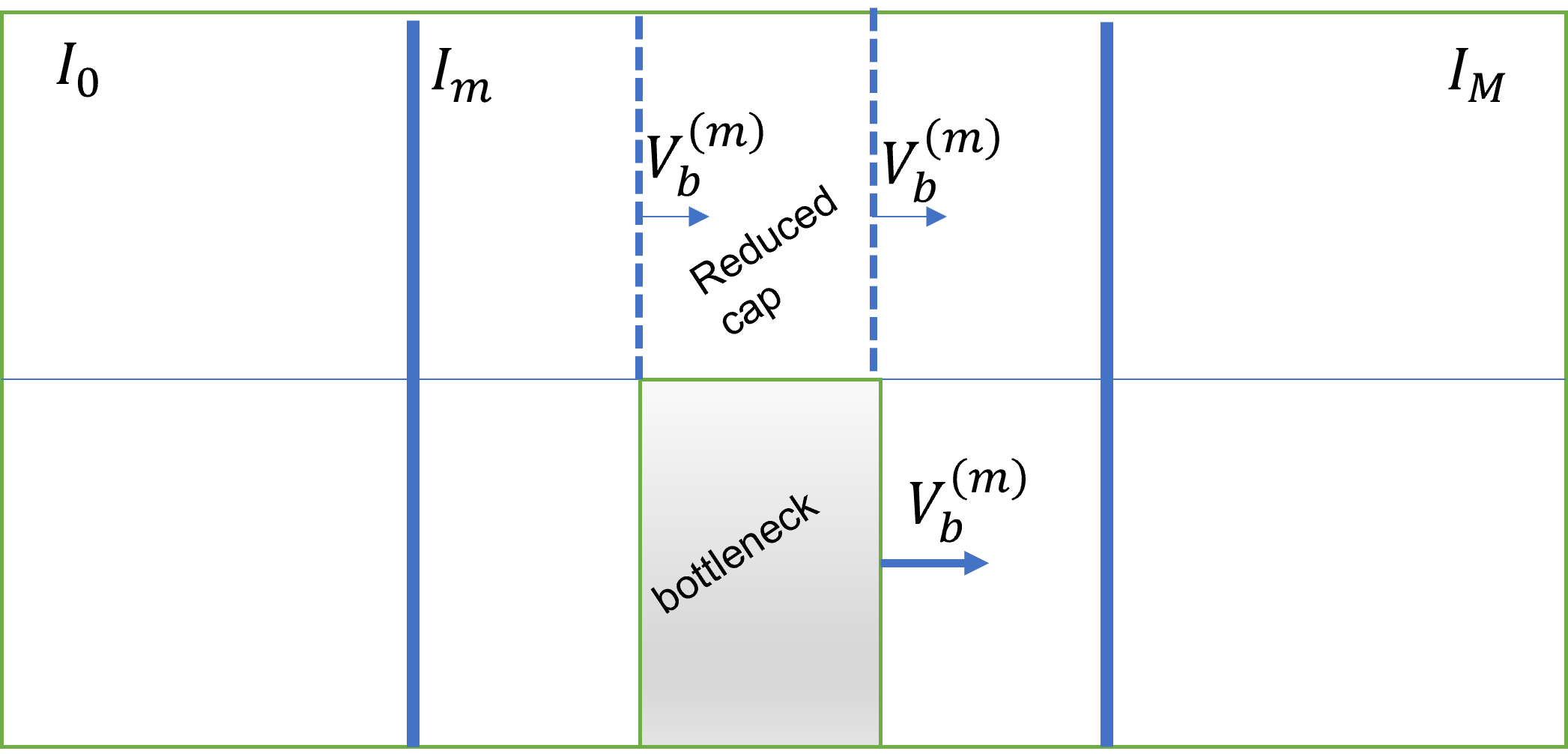}
    \caption{Capacity reduction at the location of the bottleneck.}
    \label{fig:bottleneck_cap_reduction}
\end{figure}

The goal of this section is to understand the dynamics of the bottleneck. In particular, at the location of the bottleneck the capacity of the road will be reduced (see Figure \ref{fig:bottleneck_cap_reduction}). 
To understand the dynamics in this situation, we change the coordinates to the location of the moving bottleneck by defining $X(t, x) \Def x - y(t)$ for any $t \in \RR_+$ and $x \in \RR$. In particular, $X(t,x)=0$ corresponds to the location of the bottleneck, $x = y(t)$ (the position of the bottleneck in the new coordinate will be the origin). 

The conservation of mass in new coordinates reads
\begin{equation}\label{E:conserv_new_coordinate}
    \fpartial{\rho}{t}(t, X) + \fpartial{}{X} \left[ f(\gamma(X), \rho(t, X)) - \dot y(t) \rho(t, X) \right]= 0.
\end{equation}
The capacity reduction of the road at the bottleneck location can be captured by scaling the maximum density; i.e. $\alpha \rhomax$, where $\alpha \in (0, 1)$ is the predefined reduction scale. 

\begin{figure}
    \centering
    \includegraphics[width=3in]{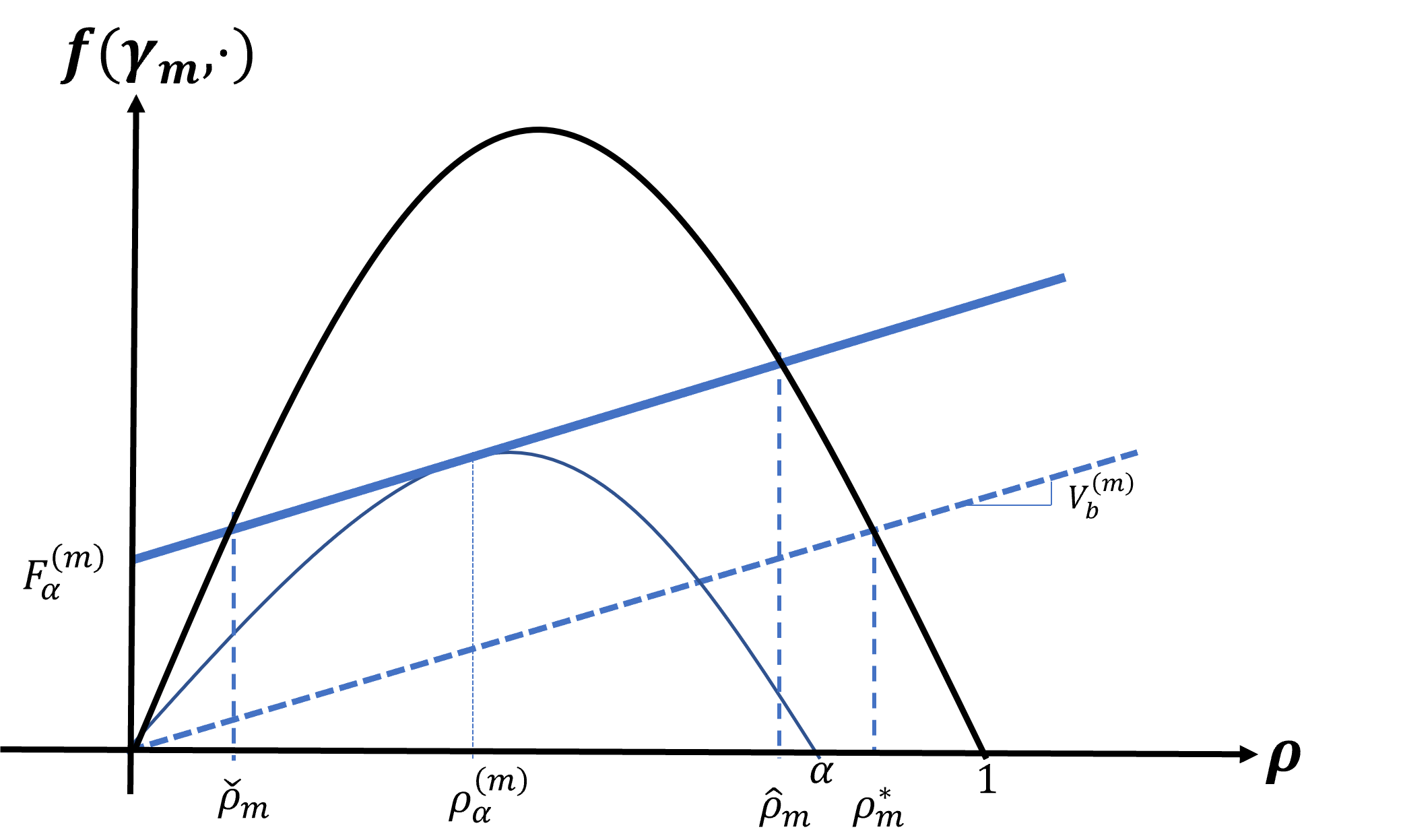}
    \includegraphics[width=3in]{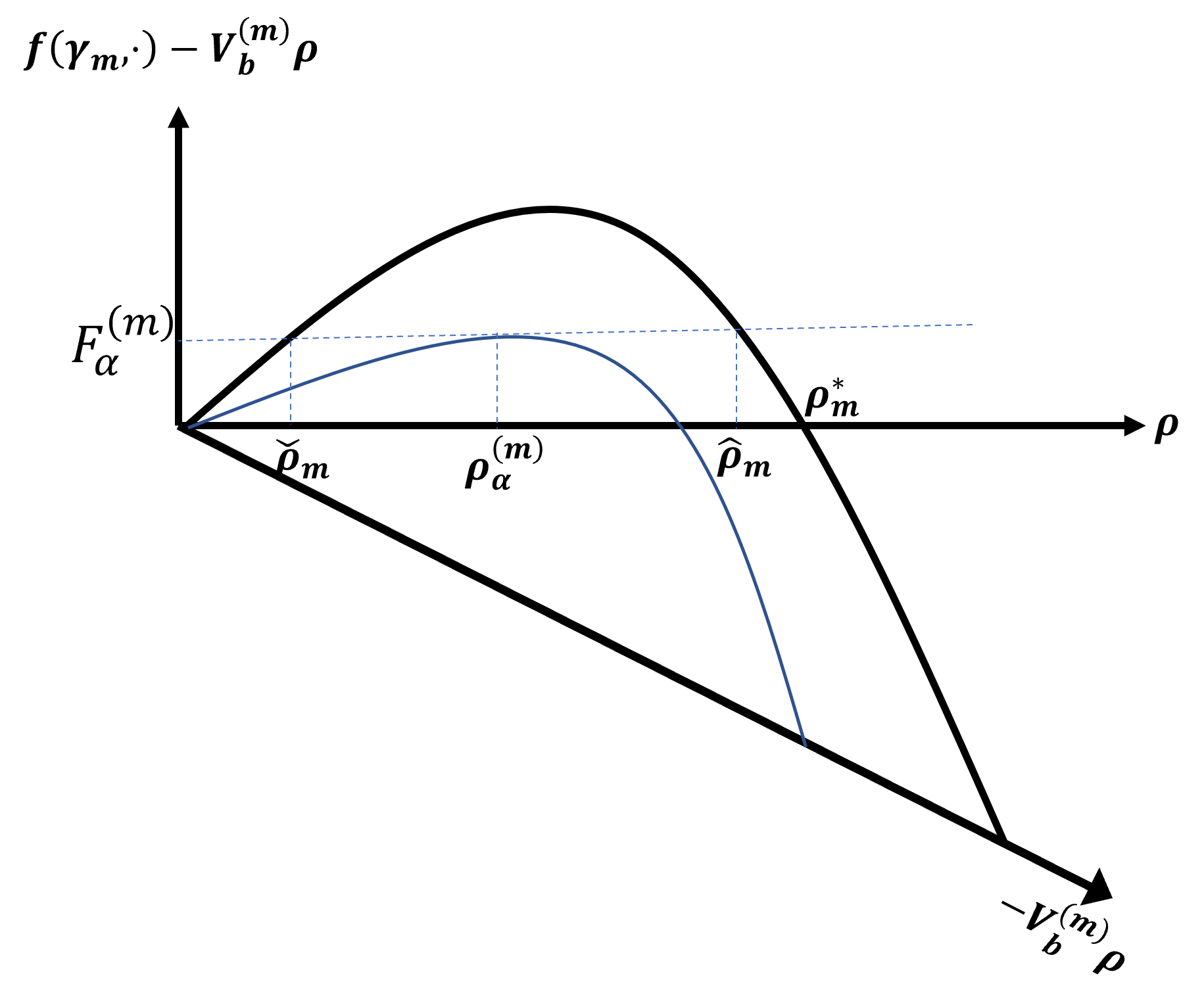}
    \caption{The left figure illustrates the bottleneck inequality. Any solution above this constraint violates the inequality. The right figure is obtained by a change of coordinates.}
    \label{fig:flux_newCoordinates}
\end{figure}

By definition of the flux function in \eqref{E:conserv_new_coordinate} and the reduction in the road capacity at the location of the bottleneck, we should have
\begin{equation}\label{E:flux_inequality}
    f(\gamma(y(t), \rho(t, y(t))) - \dot y(t) \rho(t, y(t)) \le \max_{\rho} \lb f_\alpha(\gamma(y(t), \rho)- \dot y(t) \rho \rb , \quad m \in \set{0, \cdots, M}
\end{equation}
where, function $(\gamma, \rho) \in \set{\gamma_{r_\circ},\cdots,\gamma_{r_M}} \times [0,\alpha \rhomax] \mapsto f_\alpha(\gamma, \rho) \in \RR_+$ is defined by
\begin{equation}
    \label{E:reduced_flux}
    f_\alpha(\gamma, \rho) \Def \gamma \rho \left(1 - \frac{\rho}{\alpha \rhomax} \right)
\end{equation}
and shows the reduction in the flux; consult the illustration of Figure \ref{fig:flux_newCoordinates}.
Simple calculations, show that the right-hand side of \eqref{E:flux_inequality} assumes its maximum at
\begin{equation*}
    \rho_{\alpha }^{(m)} \Def \frac{\alpha \rhomax}{2} \left(1 - \frac{\dot y(t)}{\gamma_{r_m}} \right), \quad \text{for any $m \in \set{0, \cdots, M}$}.
\end{equation*}
Therefore, Equation \eqref{E:flux_inequality} can be rewritten in the form of 
\begin{equation}
    \label{E:main_ODE_constraint}
    f(\gamma(y(t)), \rho(t, y(t))) - \dot y(t) \rho(t, y(t)) \bigmid_{y(t) \in I_m} \le \frac{\alpha \rhomax}{4 \gamma_{r_m}} \left(\gamma_{r_m} -\dot y(t) \right)^2,
\end{equation}
for any $m \in \set{0, \cdots, M}$. 
It should be noted that the inequality \eqref{E:main_ODE_constraint} is satisfied for $\dot y(t) = v(\gamma(y(t)), \rho(t, y(t))$ by \eqref{E:flux};
i.e. the left hand side of \eqref{E:main_ODE_constraint} vanishes.

For simplicity in the rest of the paper, we assume $\rhomax = 1$. We use the notation $F^{(m)}_\alpha(\dot y(t)) \Def F_\alpha(y(t),\dot y(t))$ when $y(t) \in I_m$. 
\section{Riemann Problem with Moving Bottleneck and Discontinuous Flux}
Let's start with a definition that will be used throughout the paper.
\begin{definition}\label{def:waves} Throughout these notes we adhere to the following definitions
\begin{itemize}
\item A bottleneck is active if it moves with its maximal velocity $\dot y(t) = V_b(y(t))$ and creates a queue in the upstream traffic \cite{munoz2002moving}. 

\item Fix $\gamma_{r_m}$. A non-classical shock only happen between $\hat \rho_m$ and $\check \rho_m$ where $\check \rho_m \le \hat \rho_m$ are the points of intersection of $V_b^{(m)} \rho(t, y(t)) + F^{(m)}_\alpha$ with the fundamental diagram of $\set{(\rho, f(\gamma_{r_m}, \rho): \rho \in [0, 1]}$. In addition, the non-classical shock only happens along the bottleneck trajectory. 
\end{itemize}
\end{definition}

First, we introduce a solution to the following Riemann problem:
\begin{equation}\label{E:Riemann}
    \partial_t \rho + \fpartial{}{x}[f(\gamma, \rho)] = 0, 
\end{equation}
with the initial data
\begin{equation}\label{E:Riemann_condition}
\rho(0, x) = \begin{cases}\rho_L &, x <0 \\ \rho_R &, x>0 \end{cases}, \quad 
    \gamma(x) = \begin{cases} \gamma_L &, x <0 \\
    \gamma_R &, x>0 \end{cases}.
\end{equation}
under the bottleneck constraint of
\begin{equation}\label{E:ODE_constraint}
    f(\gamma_R, \rho(t, y(t))) - \dot y(t) \rho(t, y(t))\le \frac{\alpha}{4 \gamma_R} \left(\gamma_R - \dot y(t) \right)^2,
\end{equation}
and the dynamic of 
\begin{equation}
    \dot y(t) = \min\set{V_b^{(R)}, v(\gamma_R, \rho_R)}, 
\end{equation}
for a.e. $t \in [0, \infty)$ and the initial value of $y_\circ = 0$ which implies that the location of the bottleneck $y(t) \in [0 ,\infty)$ for all $t \in \RR_+$. 

The solution of the Riemann problem \eqref{E:Riemann}-\eqref{E:Riemann_condition} is denoted by a vector-valued function 
\begin{equation*} \Riemann: (\rho_L, \rho_R; \gamma_L, \gamma_R) \in [0, 1]^2 \times \set{\gamma_\circ, \cdots, \gamma_{r_M}}^2 \to L_{\loc}^1(\RR; [0,1])\end{equation*}
and is defined according to the minimum jump entropy condition. The goal in this section is to construct a solution that in addition satisfies \eqref{E:ODE_constraint}. Let's briefly recall the Riemann solution in the case of discontinuity in the flux; see for example \cite{temple1982global, holden2015front} for more detail. 
Let's consider the initial data as in \eqref{E:Riemann_condition}. The main idea is to consider the solution $\rho(t, x)$ of the Riemann problem in the form of 
\begin{equation}
    \rho(t, x) = \begin{cases}
        v(t, x) &, x < 0 \\
        w(t, x) &, x >0
    \end{cases}
\end{equation}
where $v$ and $w$ are the Riemann solutions of with the initial data 
\begin{equation*}
    v_\circ(x) =\begin{cases}
        \rho_L &, x<0 \\
        \rho'_L&, x = 0
    \end{cases} \qquad,  w_\circ(x) = \begin{cases}
        \rho_R' &, x = 0 \\
        \rho_R &, x > 0
    \end{cases}
\end{equation*}
The task is to find $\rho'_L$ and $\rho_R'$ to define the Riemann solution $\rho(t,x)$. To do so, one main consideration would be that the speed of propagation at $x = 0^-$ should be negative and similarly at $x = 0^+$ positive. Therefore, we may define
\begin{equation}\label{E:envelope_1}
    h_L(\rho; \rho_L) = \begin{cases}
        \inf \set{h(\rho) \ge f(\gamma_L, \rho), h'(\rho) \le 0, h(\rho_L) = f(\gamma_L, \rho_L)} &, \rho \le \rho_L \\
        \sup \set{h(\rho) \le f(\gamma_L, \rho), h'(\rho) \le 0, h(\rho_L) = f(\gamma_L, \rho_L)} &, \rho \ge \rho_L
    \end{cases}
\end{equation}
\begin{equation}\label{E:envelope_2}
    h_R(\rho; \rho_R) = \begin{cases}
        \sup \set{h(\rho) \le f(\gamma_R, \rho), h'(\rho) \ge 0, h(\rho_R) = f(\gamma_R, \rho_R)} &, \rho \le \rho_R \\
        \inf \set{h(\rho) \ge f(\gamma_R, \rho), h'(\rho) \ge 0, h(\rho_R) = f(\gamma_R, \rho_R)} &, \rho \ge \rho_R
    \end{cases}
\end{equation}
Furthermore, 
\begin{equation*} H_L(\rho_L)\Def \set{\rho: h_L(\rho; \rho_L) = f(\gamma_L, \rho)} , \qquad H_R(\rho_R)\Def \set{\rho: h_R(\rho; \rho_R) = f(\gamma_R, \rho)}.\end{equation*} 
By construction, the maps $\rho \mapsto h_L(\rho;\rho_L)$ and $\rho \mapsto h_R(\rho;\rho_R)$ are decreasing and increasing respectively. Therefore, they can collide at most at one point, denoted by $\rho^\times$. Hence, $f(\gamma_L, \rho^\times) = f(\gamma_r, \rho^\times)$. Finally, $\rho_L'$ and $\rho_R'$ are chosen from minimum jump entropy condition of the form 
\begin{equation}\label{E:min_jump_entropy}
    \begin{split}
        \rho_L' &\Def \arg\min_{\rho} \set{\abs{\rho_L - \rho}: \rho \in H_L(\rho_L), \, h_L(\rho; \rho_L) = f^\times} \\
         \rho_R' &\Def \arg\min_{\rho} \set{\abs{\rho_R - \rho}: \rho \in H_R(\rho_R), \, h_R(\rho; \rho_R) = f^\times}
    \end{split}
\end{equation}

\begin{notation}\label{N:shock_speed}
In what follows, we define $\lambda(\rho_1, \rho_2)$, for any $\rho_1, \rho_2 \in [0,1]$, to be the speed of the shock wave $\rho_1$ and $\rho_2$. Such a shock front will be denoted by $\front[\rho_1, \rho_2]$ in these notes. In other words, by Rankine-Hugoniot condition 
\begin{equation}\label{E:shock_speed}
\lambda(\rho_1,\rho_2) \Def \frac{f(\gamma, \rho_1) - f(\gamma, \rho_2)}{\rho_1 - \rho_2}\end{equation}
for a fixed $\gamma$. 
\end{notation}
\begin{definition}[Riemann Solution] \label{def:Riemann_sol_new} A Riemann solver 
\begin{equation*} 
\Riemann^\alpha:[0,1]^2 \times \set{\gamma_\circ, \cdots, \gamma_{r_M}}^2 \mapsto L^1(\RR; [0,1])\end{equation*}
for \eqref{E:Riemann}-\eqref{E:ODE_constraint} is defined as follows
\begin{enumerate} 
\item If $f(\gamma_R, \Riemann(\rho_L, \rho_R; \gamma_L, \gamma_R)(V_b^{(R)})) > F_\alpha^{(R)} + V_b^{(R)} \Riemann(\rho_L, \rho_R; \gamma_L, \gamma_R)(V_b^{(R)})$, then
\begin{equation}
    \label{E:case_1_new_solution_top_area}
    \Riemann^\alpha (\rho_L, \rho_R; \gamma_L, \gamma_R)(\nicefrac{x}{t}) \Def \begin{cases}
        \Riemann(\rho_L, \hat \rho_R; \gamma_L, \gamma_R)(\nicefrac xt) &, \frac xt < V_b^{(R)}\\
        \Riemann(\check \rho_R, \rho_R; \gamma_R)(\nicefrac{x}{t}) &, \frac xt \ge V_b^{(R)}
    \end{cases} 
\end{equation}
 
 This case is associated with a non-classical shock.
\item If 
\begin{multline*}
    V_b^{(R)} \Riemann(\rho_L, \rho_R; \gamma_L, \gamma_R)(V_b^{(R)}) \le f(\gamma_R, \Riemann(\rho_L, \rho_R; \gamma_L, \gamma_R)(V_b^{(R)})) \\\
    \le F_\alpha^{(R)} + V_b^{(R)} \Riemann(\rho_L, \rho_R; \gamma_L, \gamma_R)(V_b^{(R)})
\end{multline*} then 
\begin{equation*}
    \Riemann^\alpha(\rho_L, \rho_R; \gamma_L, \gamma_R)(\nicefrac{x}{t})\Def \Riemann(\rho_L, \rho_R; \gamma_L, \gamma_R)(\nicefrac{x}{t}), \quad  y(t) = V_b^{(R)} t.
\end{equation*}

\item If $f(\gamma_R, \Riemann(\rho_L, \rho_R; \gamma_L, \gamma_R)(V_b^{(R)})) \le V_b^{(R)} \Riemann(\rho_L, \rho_R; \gamma_L, \gamma_R)(V_b^{(R)}) $, then 
\begin{equation*}
    \Riemann^\alpha(\rho_L, \rho_R; \gamma_L, \gamma_R)(\nicefrac{x}{t})\Def \Riemann(\rho_L, \rho_R; \gamma_L, \gamma_R)(\nicefrac{x}{t}), \quad  y(t) = v(\gamma_R, \rho_R) t.
\end{equation*}

\end{enumerate}

\end{definition}
Depending on the values of $\rho_L, \rho_R \lessgtr \tfrac 12$ and $f(\gamma_L, \rho_L) \lessgtr f(\gamma_R, \rho_R)$ different cases of Riemann solution can happen. Below, we will investigate of one of these cases and some other possibilities are postponed to the supplementary materials. 
\begin{case}\label{case:R_L_1}
Let $\rho_L < \frac 12$, $\rho_L < \rho_R$ and $f(\gamma_L, \rho) < f(\gamma_R, \rho)$. The Riemann solution of \eqref{E:Riemann}-\eqref{E:Riemann_condition} (without the bottleneck-constraints) can be explicitly written as
\begin{equation}
    \label{E:Riemann_case1}
    \Riemann(\rho_L, \rho_R; \gamma_L, \gamma_R)(\nicefrac{x}{t}) = \begin{cases}
    \rho_L &, x<0 \\
    \rho'_R &, \frac xt \in [0, \lambda(\rho'_R, \rho_R)) \\
    \rho_R &, \frac xt \ge \lambda(\rho'_R, \rho_R).
    \end{cases}
\end{equation}
and is illustrated in the Figure \ref{fig:case_1_no_ODE}. The solution in this case consists of a $\gamma$-jump between $\rho_L = \rho'_L$ and $\rho'_R$ followed by a shock wave between $\rho'_R$ and $\rho_R$. If this solution satisfies the bottleneck inequality, then $\Riemann^\alpha = \Riemann$. However, if the the solutions $\rho'_R$ and $\rho_R$ are located above the bottleneck inequality, then the Riemann solution needs to be redefined. Figure \ref{fig:case1_modified_Riemann} illustrates the modified Riemann solution of Case \ref{case:R_L_1}. 
\begin{figure}
    \centering
    \includegraphics[width=3in]{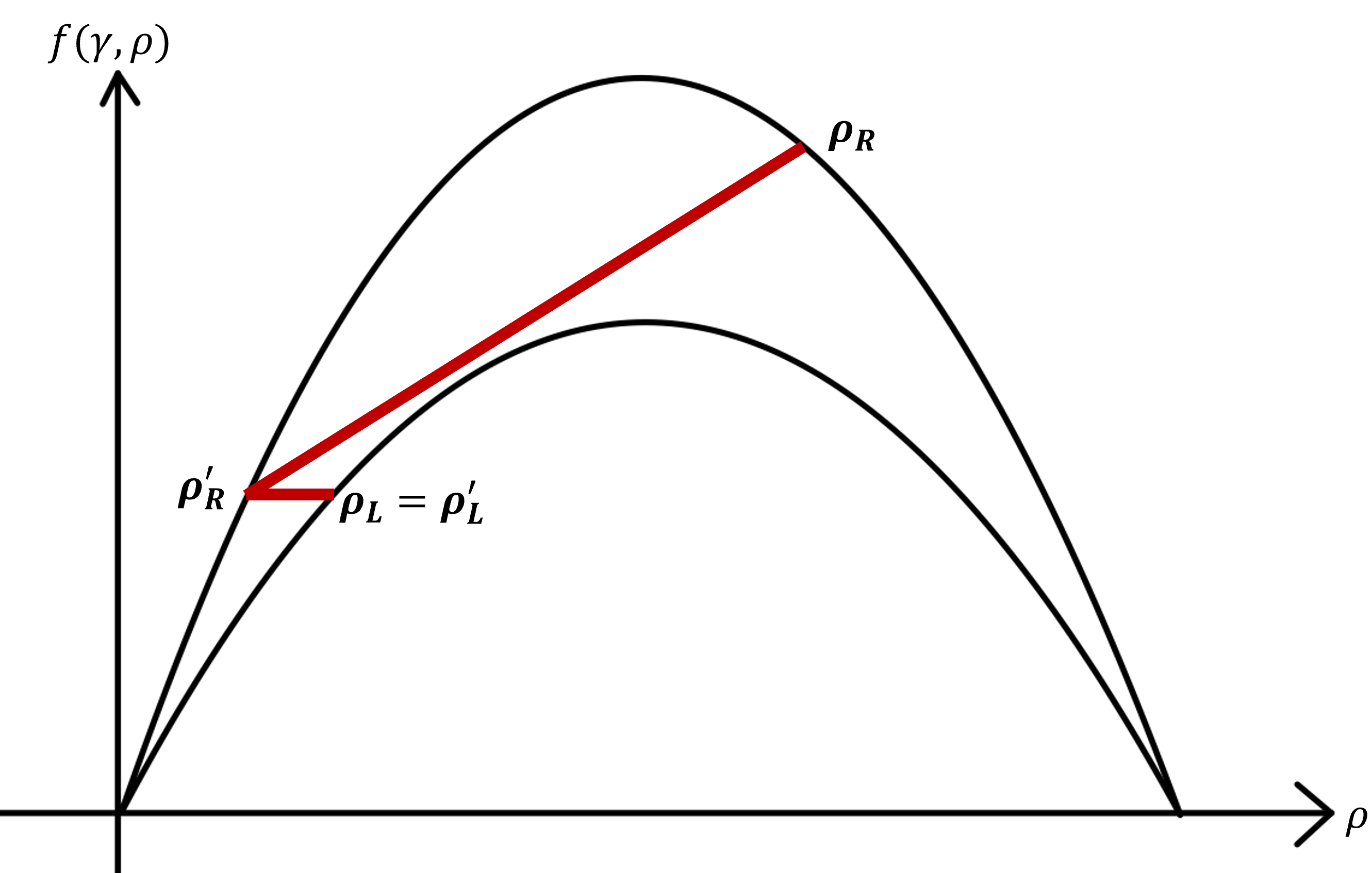}
    \caption{Riemann solution for $\rho_L < \frac 12$, $f(\gamma_L, \rho) < f(\gamma_R, \rho)$ and $\rho_R > \rho_L$.}
    \label{fig:case_1_no_ODE}
\end{figure}
\begin{figure}
    \centering
    \includegraphics[width=3in]{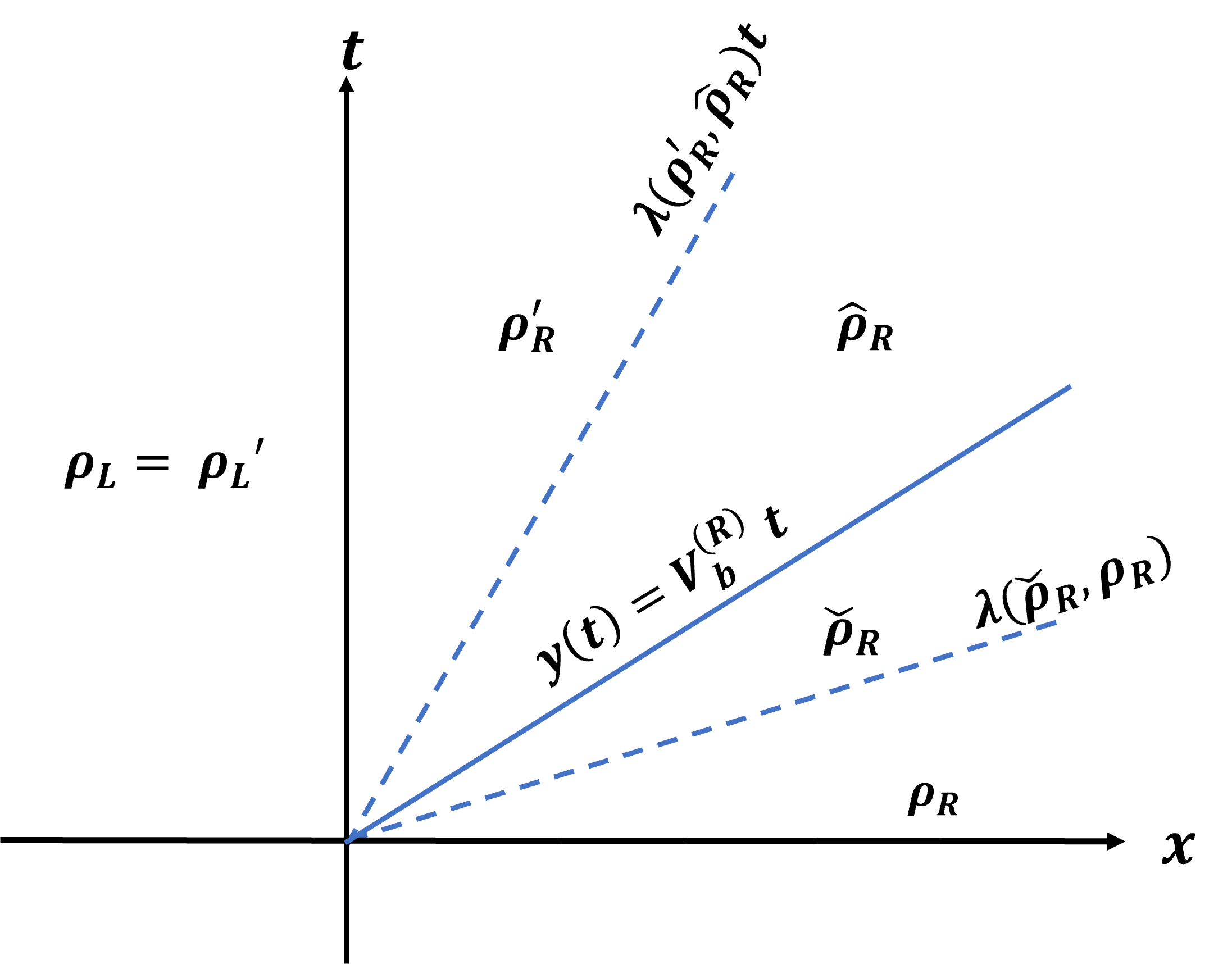}
    \caption{The Riemann solution $\Riemann^\alpha(\rho_L, \rho_R; \gamma_L, \gamma_R)$ in the case the Riemann solution.}
    \label{fig:case1_modified_Riemann}
\end{figure}
\begin{remark}\label{R:tempting_solution_case1}
In fact, it might be tempting to redefine the Riemann solution \eqref{E:case_1_new_solution_top_area} only on the location of the bottleneck (i.e. on $x\ge0$) without considering the solution on $x <0$. 
 However, this might violate the minimum jump entropy condition and the solution is not acceptable \tup{see figure \ref{fig:Case_1_solution_violation}}.  
\begin{figure}
    \centering
    \includegraphics[width=3in]{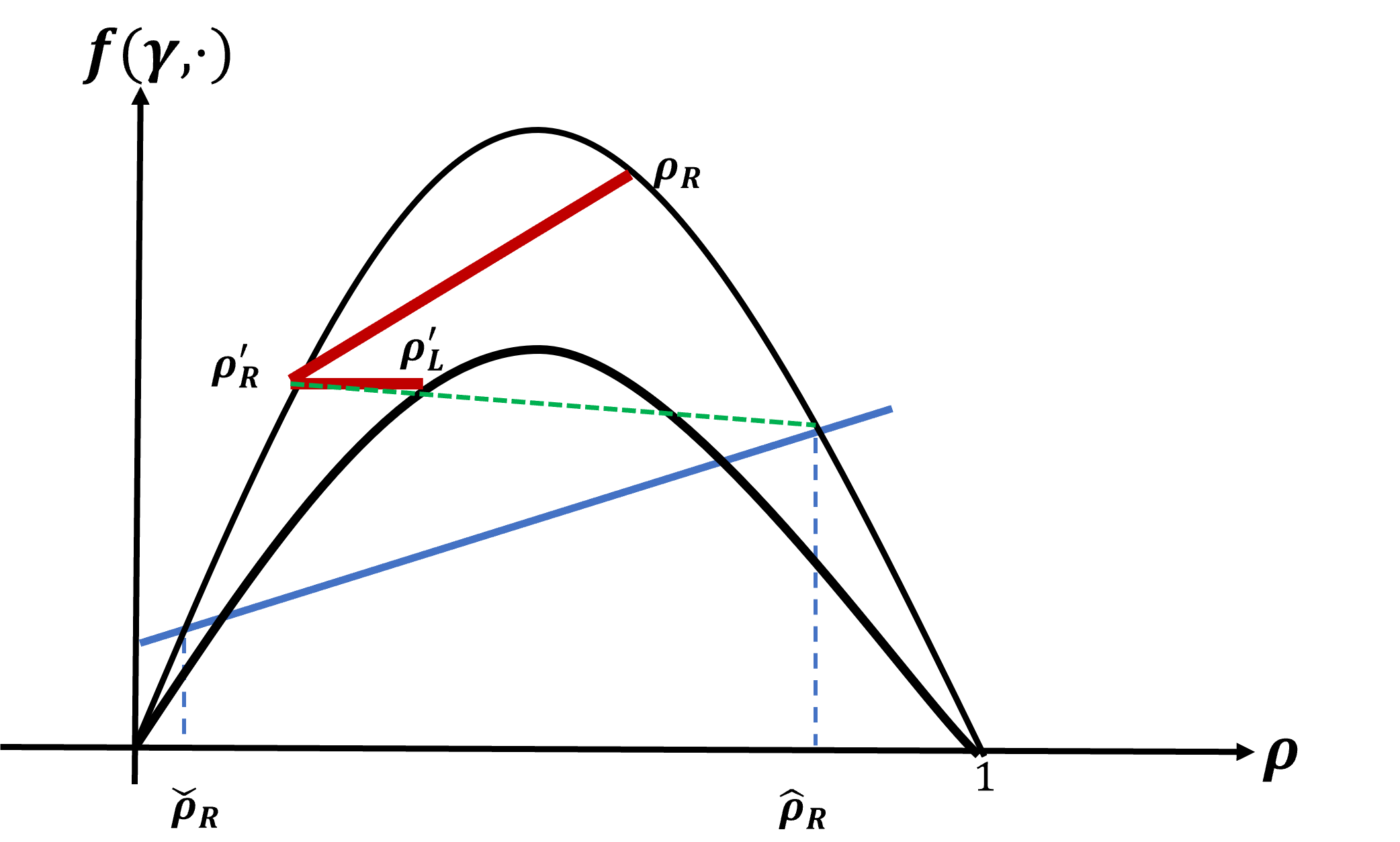}
    \caption{The possibility that $f(\gamma_R, \rho'_R) > f(\gamma_R, \hat \rho_R)$ (the green dashed line).}
    \label{fig:Case_1_solution_violation}
\end{figure}
\end{remark}
\noindent We refer the readers to \cite[Remark 1]{delle2014scalar} for some properties of the Riemann solution. 
\end{case}

The conservation law in the presence of the discontinuity in the flux can be written as a system of the form 
\begin{equation*}
    \begin{cases}
        \rho_t + \frac{\partial}{\partial x} [f(\gamma, \rho)] & = 0 \\
        \gamma_t &= 0
    \end{cases}
\end{equation*}
Then the Jacobian matrix for this system will have the eigenvalues of $\lambda_1= 0$ and $\lambda_2 = \tfrac{\partial f}{\partial \rho}$. In other words, in the presence of the discontinuity in the flux the strict hyperbolicity will be lost (coinciding eigenvalues and consequently resonant system) which implies that the solution will not be of bounded variation. Consequently, the compactness theorem cannot be employed to show the existence of the solution. Instead, it is customary to define a homeomorphism between the space of solutions determined by the front tracking scheme to another space in which the solution satisfies the requirements of the compactness theorem, \cite{gimse1992solution, temple1982global}. 


More accurately, for any fixed $\gamma \in \set{\gamma_{r_\circ}, \cdots, \gamma_{r_M}}$, we define a bijection $\rho \in [0, 1] \mapsto \psi(\gamma, \rho) \in \RR$ by 
\begin{equation}\label{E:homeo_map}\begin{split}
    \psi(\gamma, \rho) &\Def \sign\left (\tfrac 12  - \rho \right) \left(f (\gamma, \rho) - f(\gamma, \tfrac 12) \right)\\
    & = \tfrac 14 \gamma \sign\left(\rho - \tfrac 12 \right) (2\rho - 1)^2.
\end{split}\end{equation}
In particular, let 
\begin{equation}\label{E:state_spaces}
    \mathcal V \Def \set{(\gamma, \rho): \gamma \in \set{\gamma_\circ, \cdots, \gamma_{r_M}}, \rho\in  [0, 1]}, \quad \mathcal W \Def \set{(\gamma, z): \gamma \in\set{\gamma_\circ, \cdots, \gamma_{r_M}}, z \in [-\nicefrac \gamma 4, \nicefrac \gamma 4]}
\end{equation}
then, the map
\begin{equation}\label{E:homeo}
   \mathcal V \overset{\psi}{\mapsto} \mathcal W
\end{equation}
where for each fixed $\gamma \ne 0$, defines a homeomorphism.
In addition, \eqref{E:homeo_map} implies that for each $\gamma \ne 0$, the inverse function can be defined by 
\begin{equation}
    \label{E:psi_inverse}
    \psi^{-1}(\gamma, z) = \frac 12 \left(1+ \sign (z) \sqrt{\frac{4\abs{z}}{\gamma}} \right), \quad \gamma \in \set{\gamma_{r_\circ}, \cdots, \gamma_{r_M}}, \quad z \in [ \nicefrac{-\gamma}{4}, \nicefrac{\gamma}{4}].
\end{equation}
\subsection{Riemann Solution in $\mathcal W$}
We recall that the Riemann solution in $\mathcal V$ space which is introduced in Definition \ref{def:Riemann_sol_new}. Since we will be working in the $\mathcal W$ space, here we discuss the properties of the Riemann solution in this space. The following observations define the Riemann solution in $\mathcal W$-space (sometimes it is referred to as $(z, \gamma)$-state space in what follows).
\begin{enumerate}
    \item The $\rho$-waves (they are also called $z$-waves in $\mathcal W$-space are horizontal lines (parallel to the $z$-axis). 
\item The $\gamma$-waves have the slope of $\pm \tfrac 14$ in $\mathcal W$ space.
To see this, we need the following result.
\begin{lemma}\label{L:gamma_jump_state_closedness}
If $\rho'_L$ and $\rho'_R$ are connected by a $\gamma$-front and are determined by the minimum jump entropy condition, then 
\begin{equation*}
    \rho'_L , \rho'_R \le \frac 12, \quad \text{or $\rho'_L, \rho'_R \ge \frac 12$}.
\end{equation*}
\end{lemma}
\begin{proof}
By definition of the minimum jump entropy in \eqref{E:min_jump_entropy}, the proof is immediate. 
\end{proof}
For the following discussion, the readers are encouraged to consult the example illustrated in Figure \ref{fig:z_gamma_bottleneck_1}. As in the classical case, in $(z, \gamma)$-space the $\gamma$-fronts happen along a straight line with the slope of $\pm \frac 14$. To see this, we consider $\tilde \rho'_L$ and $\hat \rho'_R$ be determined by the minimum jump entropy condition as in \eqref{E:min_jump_entropy} and be connected by the $\gamma$ front. Then, by the first presentation of \eqref{E:homeo_map} we define
\begin{equation}\label{E:solution_case1}
\begin{split}
    \tilde z'_L & \Def \psi(\gamma_L, \tilde \rho'_L) = \sign(\tfrac 12 - \tilde \rho_L') \left(f(\gamma_L, \tilde \rho_L') - f(\gamma_L, \tfrac 12) \right)\\
    \hat z_R' &\Def \psi(\gamma_R, \hat \rho_R') = \sign(\tfrac 12 - \hat \rho_R') \left(f(\gamma_R, \hat \rho_R') - f(\gamma_R, \tfrac 12) \right).
\end{split}\end{equation}
By Lemma \ref{L:gamma_jump_state_closedness}, using the fact that $f(\gamma_L, \tilde \rho'_L) = f(\gamma_R, \hat \rho_R')$ and $f(\gamma, \tfrac 12) = \tfrac 14 \gamma$, we conclude that 
\begin{equation}\label{E:z_gamma_relation}
    \frac{\hat z_R' - \tilde z_L'}{\gamma_R - \gamma_L} = \pm \tfrac 14
\end{equation}
the line connecting $\tilde z_L'$ and $\hat z_R'$ has a slope of $\pm \tfrac 14$. This remains correct for any $\gamma$-jump. Such observation is crucial in calculating the total variation of the approximate solutions in this space.
\begin{figure}
    \centering
    \includegraphics[width=\textwidth]{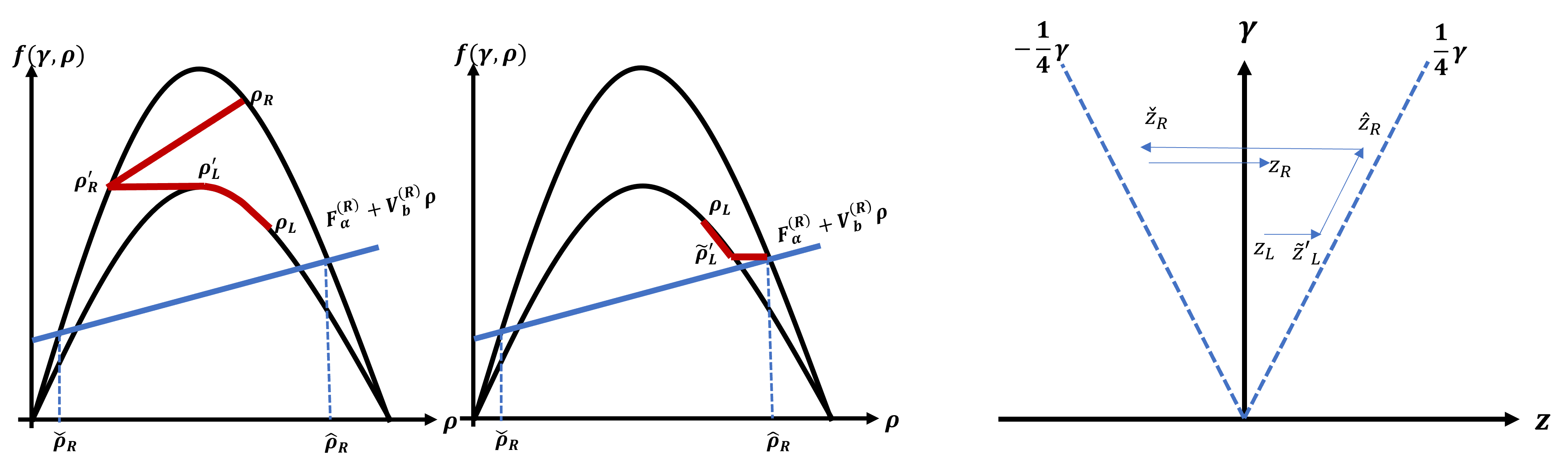}
    \caption{The left figure illustrates the solution of the Riemann problem without the moving bottleneck. The figure in the middle depicts the proposed Riemann solution in the presence of the moving bottleneck. The right figure is the solution of the proposed Riemann problem in $\mathcal W$-space.}
    \label{fig:z_gamma_bottleneck_1}
\end{figure}
\item Let's consider $\rho > \nicefrac{1}{2}$. Then, 
\begin{equation*}
    \psi(\gamma, \rho) = - \left(f(\gamma, \tfrac 12) - f(\gamma, \rho)\right) > -\tfrac 14 \gamma
\end{equation*}
as $f(\gamma, \rho) >0$. Similarly, argument for $\rho < \nicefrac{1}{2}$, we have that 
\begin{equation*}
    \psi(\gamma, \rho) = f(\gamma, \tfrac 12) - f(\gamma, \rho) < \tfrac 14 \gamma
\end{equation*}
Therefore, $z = \pm \tfrac 14 \gamma$ provides a bound for the relative variation of the solutions in the $(z, \gamma)$-space (see the right illustration of Figure \ref{fig:z_gamma_bottleneck_1}).


\item By the construction of the solution in $\mathcal W$-space, we note that if the direction along the horizontal line is to the left, the corresponding wave is either a rarefaction or a non-classical shock (which is along the bottleneck trajectory and going from $\hat z$ to $\check z$) and the right direction implies the shock wave.
\end{enumerate}
\section{Cauchy Problem}
In this section, we gradually build some required results which will be employed later to show the existence of the solution of the Cauchy problem \eqref{E:system}. 
Let's begin defining a solution.
\begin{definition}\label{def:existence}
A functional $(\rho, y) \in  C(\RR_+; L^1_{\loc}(\RR; [0,1])) \times W^{1,1}_{\loc}(\RR_+, \RR)$ is a solution to the Cauchy problem \eqref{E:system} with initial value $\rho_\circ \in  \BV(\RR; [0, 1])$, if 
\begin{enumerate}
    \item The density function $\rho$ satisfies
    \begin{equation*}
        \int_{\RR_+} \int_{\RR } (\rho \partial_t \varphi + f(\gamma, \rho) \partial_x \varphi) dx dt + \int_{\RR}\rho_\circ(x) \varphi(0, x) dx = 0, \quad \forall \varphi \in C_c^\infty(\RR_+ \times \RR);
    \end{equation*}
    \item On $(0, T) \times \RR \setminus \set{(t, y(t)): t \in \RR_+}$, for any $\varphi\in C_c^1([0, T] \times \RR; \RR^+)$ with $\varphi(t, y(t)) = 0$ and for any constant $c \in \RR$,
    \begin{multline}
    \label{E:entropy}
        \int\int\lb \abs{\rho - c} \varphi_t + \digamma(\gamma, \rho, c)\varphi_x \rb dx dt 
        + \sum_{i = 1}^{M} \int_0^T \abs{f(\gamma_i^+, c) - f(\gamma_i^-, c)} \varphi(t, \sa_i) dt \\
        + \int_\RR \abs{\rho_\circ(x)- c} \varphi(0, x)dx \ge 0 
    \end{multline}
    where, 
    \begin{equation*}
        \digamma(\gamma, \rho, c) \Def sgn(\rho -c) (f(\gamma, \rho) - f(\gamma, c) )
    \end{equation*}
    \item The bottleneck trajectory $y$ is a Carath\'eodory solution to the dynamic $\dot y$ in \eqref{E:system}, i.e. for a.e. $t \in \RR_+$, 
    \begin{equation}\label{E:BN_ODE}
        y(t) = y_\circ + \int_0^t \omega(y(s), \rho(s, y(s)+)) ds.
    \end{equation}
    In other words, $y \in \AC([0, T]; \RR)$ for any $T >0$, where $\AC$ is the class of absolutely continuous functions. 
    \item The flux constraint in \eqref{E:system} is satisfied in the sense that, for a.e. $t \in \RR_+$
    \begin{equation}\label{E:BN_cap}
        \lim_{x \to  y(t) \pm} f(\gamma(x), \rho(t, x)) - \rho(t, x)\dot y(t) \le F_\alpha(y(t), \dot y(t)). 
    \end{equation}
\end{enumerate}
\end{definition}
The main result of this paper is as follows:
\begin{theorem}[Existence of Cauchy Solution]\label{T:existence_Cauchy}
    Let $\rho_\circ \in\BV(\RR; [0, 1])$. Then the Cauchy problem \eqref{E:system} has a solution in the sense of Definition \ref{def:existence}.
\end{theorem}
In the rest of this paper, we construct the proof of this theorem. 
Broadly speaking, we employ the wavefront tracking scheme to show the existence of the solution in the sense of Definition \ref{def:existence}. We start with defining the approximate problems. More precisely, we approximate the initial condition $\rho_\circ$ by a sequence of piecewise constant functions $\set{\rho^{(n)}_\circ: n \in N_\circ}$ and the flux function $f(\gamma, \cdot)$ by a sequence of piecewise continuous function $\set{f^{(n)}(\gamma, \cdot): n \in N_\circ}$ for some $N_\circ \in \NN$ which will be discussed later in this section. Then, for each $n \in N_\circ$, we show that the solution of the corresponding approximate problem satisfies the conditions of the compactness theorem in $\mathcal W$ space, and hence as $n \to \infty$ the solution of the Cauchy problem exists in the limiting sense. 
\subsection{Construction of Grid Points}
For $n \in \NN$ we consider the discretization
\begin{equation}\label{E:grid}\begin{split}
   & \gamma^{(n)}_k \Def \tfrac{k}{2^n}, \quad k \in \NN\\
   & z^{(n)}_{k, j} \Def \tfrac 14\tfrac{j}{2^n}, \quad j \in \ZZ \cap [-k,k], \text{for any $k \in \NN$}. 
\end{split}\end{equation}
The corresponding grid points in $\mathcal V$ can be determined by 
\begin{equation*}
    \rho^{(n)}_{k, j} = \psi^{-1}(\gamma^{(n)}_k, z^{(n)}_{k, j}) = \tfrac 12 \left(1 + \sign(z^{(n)}_{k, j}) \sqrt{\tfrac{4 \abs{z^{(n)}_ {k, j}}}{\gamma^{(n)}_k}}\right)
\end{equation*}
The collection of these points for any fixed $k \in \NN$ (or equivalently, for any fixed $\gamma^{(n)}_k$) are denoted by a set
\begin{equation}\label{E:breakpoints}
    \mathcal U_k^{(n)} \Def \set{\rho^{(n)}_{k, j} = \psi^{-1}(\gamma^{(n)}_k, z^{(n)}_{k, j}): j \in \ZZ \cap [-k, k]}, 
\end{equation}
and for any set $\mathcal K \subset \NN$, 
\begin{equation}
    \mathcal U^{(n)}_\mathcal K \Def \bigcup_{k \in \mathcal K} \mathcal U_k^{(n)}.
\end{equation}
The image set $\psi(\mathcal U^{(n)}_k)$ of $\mathcal U^{(n)}_k$ contains the corresponding solution $z$.

\begin{figure}
    \centering
    \includegraphics[width=3.5in]{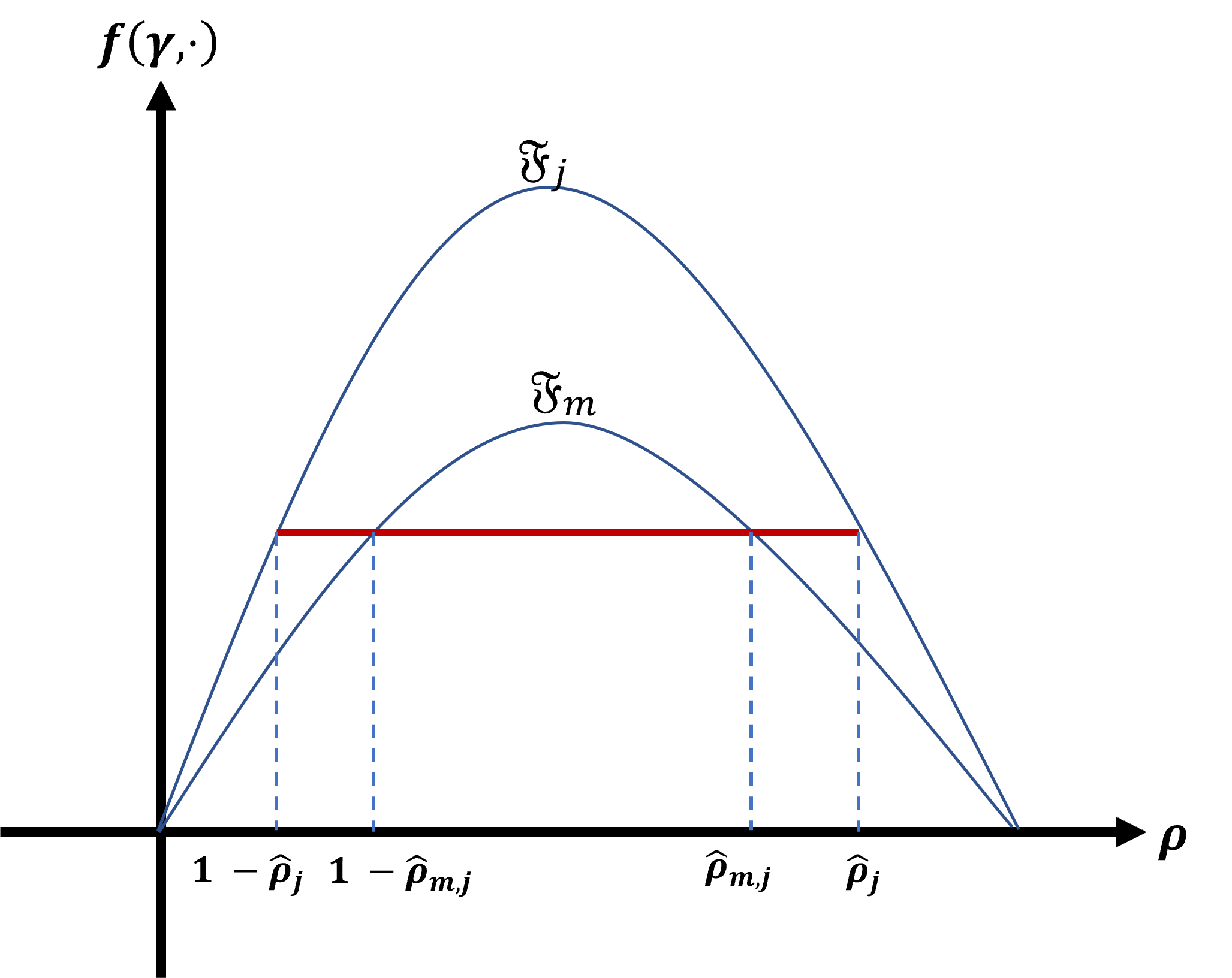}
    \caption{For each $m \in \set{0, \cdots, M}$ we add $\hat \rho_m$, $1 - \hat \rho_m$ and $\hat \rho_{m, j}$, $j \in \set{0,\cdots, M}$ on the fundamental diagram $\mathfrak F_m$}.
    \label{fig:collection_grids}
\end{figure}

Moreover, the points of intersection of the bottleneck constraints $F_\alpha^{(m)} + V_b^{(m)} \rho$ (i.e. the flux constraint when the bottleneck is in the $I_m$-region) with the fundamental diagram 
\begin{equation}\label{E:FD}
\mathfrak F_m \Def \set{(\rho, f(\gamma_{r_m}, \rho)): \rho \in [0, 1]}\end{equation}
are denoted by $\check \rho_m$ and $\hat \rho_m$ where $\check \rho_m < \hat \rho_m$, for any $m\in \set{0, \cdots, M}$. 
\begin{remark}\label{R:piecewise_gamma}
It is important to notice that, since $x \mapsto \gamma(x)$ is piecewise constant, $\check \rho_m$ and $\hat \rho_m$ are independent of the sequence index $n$. For the same reason, in \eqref{E:grid} we are mainly concerned with $\gamma^{(n)}_{k} = \gamma_{r_m}$ for $m \in \set{0, \cdots, M}$ and for some $k \in \NN$ and in particular, $k = k_m^{(n)} = \gamma_{r_m} 2^n$. 
\end{remark}
\noindent For any $m \in \set{0, \cdots, M}$, we define (readers are advised to consult Figure \ref{fig:collection_grids} for the following setup)
\begin{equation}\label{E:Ahat_set}\begin{split}
  \hat{\mathcal A}_m & \Def 
 \set{\hat \rho_m} \cup \set{1 - \hat \rho_{m}} \\
 \qquad &\cup \set{\hat \rho_{m, j}\in [0, 1]:  f(\gamma_{r_m}, \hat \rho_{m, j}) = f(\gamma_{r_j}, \hat \rho_j), j \in \set{0, \cdots, M} \backslash \set{m}}\\
 \qquad &\cup \set{1- \hat \rho_{m,j}\in [0, 1]:  f(\gamma_{r_m}, \hat \rho_{m, j}) = f(\gamma_{r_j}, \hat \rho_j), j \in \set{0, \cdots, M} \backslash \set{m}}
\end{split}\end{equation}
Let's elaborate on the concept of $\hat{\mathcal A}_m$ more carefully. 
The first set contains $\hat \rho_m$ which is directly determined by the intersection of the bottleneck constraint with the fundamental diagram $\mathfrak F_m$. The second set contains the symmetry of $\hat \rho_m$ on the fundamental diagram $\mathfrak F_m$ (We consider the symmetric points to keep the approximate functions symmetric). The element of the third set, i.e. $\hat \rho_{m, j}$, corresponds to projection of the point $(\hat \rho_j, f(\gamma_{r_j}, \hat \rho_j))$, $j \in \set{0, \cdots, M} \setminus \set{m}$, on the fundamental diagram $\mathfrak F_m$. We recall that preserving these stationary jumps from one fundamental diagram to the other is crucial to ensure that the solution of a Riemann problem determined by the minimum jump entropy condition remains a grid point. 
To collect all these points, we set 
\begin{equation}
    \hat {\mathcal A} \Def \bigcup_{m= 0}^{M} \hat {\mathcal A}_m.
\end{equation}
Similarly, we denote 
\begin{equation}\label{E:Acheck_set}\begin{split}
  \check{\mathcal A}_m & \Def 
 \set{\check \rho_m} \cup \set{1 - \check \rho_{m}} \\
 \qquad &\cup \set{\check \rho_{m,j}\in [0, 1]:  f(\gamma_{r_m}, \check \rho_{m, j}) = f(\gamma_j, \check \rho_j), j \in \set{0, \cdots, M} \backslash \set{m}}\\
 \qquad &\cup \set{1- \check \rho_{m, j}\in [0, 1]:  f(\gamma_{r_m}, \check \rho_{m, j}) = f(\gamma_j, \check \rho_j), j \in \set{0, \cdots, M} \backslash \set{m}}
\end{split}\end{equation}
and 
\begin{equation}
    \check {\mathcal A} \Def \bigcup_{m= 0}^{M} \check {\mathcal A}_m.
\end{equation}
Finally, we define $\mathcal A^*_m$ and $\mathcal A^*$ in a similar way. In addition, we will use the image sets $\psi(\hatBN_m)$ and $\psi(\checkBN_m)$ for the corresponding $\hat z_m$ and $\check z_m$, and comparably for $\hatBN$, $\checkBN$, $\mathcal A^*_m$ and $\mathcal A^*$. 


For the analytical purpose in this work (see proof of Theorem \ref{T:doty_n}), we need to have $\nicefrac{\underline \delta^{(n)}}{\bar \delta^{(n)}} = \mathcal O(1)$, where $\underline \delta^{(n)}$ and $\bar \delta^{(n)}$ are the minimum and maximum bounds on the grid points' distance (see Definition \ref{R:UL_delta} below). More precisely, we will need to ensure that $\underline \delta^{(n)} $ and $\bar \delta^{(n)}$ are of the same order. 
To do so, we consider the following procedure to update the collection of grid points. 
\begin{enumerate}
    \item Let 
    \begin{equation*}
        \delta_{\min} \Def \min_{m \in \set{0, \cdots, M}} \min \set{\abs{z- z'}:{z, z' \in \psi(\hatBN_m) \cup \psi(\checkBN_m) \cup \psi(\mathcal A^*_m) \cup \set{- \tfrac 14 \gamma_{r_m}, 0, \tfrac 14 \gamma_{r_m}} }}
    \end{equation*}
    In other words, $\delta_{\min}$ only considers the minimum distance between points $\hat z_m$, $\check z_m$ and $z_m^*$ which are fixed points and the extreme points of $z$ which by the construction of the Riemann solution should always remain as the grid points. 

    \item The value of $N_\circ$ can be chosen sufficiently large and uniquely such that $\delta_{\min} = \frac{\lambda}{2^{N_\circ + 2}}$ for $\exists \lambda \in [1, 2)$ (see \eqref{E:grid}) which implies that $\tfrac{1}{2^{N_\circ + 2}} < \delta_{\min}$. For the rest of the paper, we will be interested in $n \ge N_\circ$. 

    \item Update $\mathcal U_\NN^{(n)}$: For any $m, r \in \set{0, \cdots, M}$, and any $\bar z_{m,r} \in \psi(\hatBN_m) \cup \psi(\checkBN_m) \cup \psi(\mathcal A^*_m)$ (i.e. any of the points whose location is always fixed; see \eqref{E:Acheck_set}) if 
    \begin{equation*}
        \min \set{\abs{z^{(n)}_{m,j} - \bar z_{m, r}}: j \in [-m  ,  m] \cap \ZZ} < \frac{1}{2^{N_\circ + 2}}
    \end{equation*}
    where, $z^{(n)}_{m, j} \in \psi(\mathcal U^{(n)}_\NN)$, then remove $\arg\min_j \abs{z^{(n)}_{m,j} - \bar z_{m, r}}$ and all the associated points from all the grid points. Here, by the associated point of $z^{(n)}_{m, j}$ we refer to all $z^{(n)}_{m', j'} \in \psi(\mathcal U^{(n)}_\NN)$ such that $\frac{z^{(n)}_{m, j} -z^{(n)}_{m', j'}}{\gamma_{r_m}^{(n)} - \gamma_{r_m'}^{(n)}} = \pm \frac 14$ (the $\pm$ sign depends on whether the $z$ front is located on the positive or negative side) and their symmetric points with respect to the $\gamma$-axis (consult Figure \ref{fig:collection_grids}).
    \begin{definition}\label{R:UL_delta} Let's define the set of all grid points by
\begin{equation}
    \label{E:all_grids}
    \mathcal G^{(n)} \Def \mathcal U^{(n)}_\NN \cup \hat {\mathcal A} \cup \check{\mathcal A} \cup \mathcal A^*
\end{equation}
By the process of updating $\mathcal U^{(n)}_\NN$ (or comparably its image set), for any $n \ge N_\circ$ and $z_m, z'_m \in \psi(\mathcal G^{(n)})$ we have 
\begin{equation*} 
\underline \delta^{(n)} \le \abs{z_m - z_m'} < \bar \delta^{(n)}
\end{equation*}
where, $\underline \delta^{(n)} = \tfrac 12 \hat \delta^{(n)}$, $\bar \delta^{(n)} =  2 \hat \delta^{(n)}$, and $\hat \delta^{(n)} \Def \tfrac{1}{2^{n+ 2}}$. 
    \end{definition}
\end{enumerate}

\begin{notation} \label{N:proj_FD_m} 
For any $n \ge N_\circ$ and for a fixed $m \in \set{0, \cdots, M}$ \tup{and consequently, a fixed $\gamma_{r_m}$}, we denote by $\mathcal G^{(n)}_m$ the projection of $\mathcal G^{(n)}$ on the fundamental diagram $\mathfrak F_m$ and the elements are denoted by $ \rho^{(n)}_{m, j}$.
In addition, the grid points for each fixed $\gamma_{r_m}$ have increasing order, i.e. $\rho^{(n)}_{m,j} < \rho^{(n)}_{m, j+ 1}$
\end{notation}

Next, we define the piecewise linear approximation of the flux function $f$ by
\begin{equation}
    \label{E:flux_approx}
    f^{(n)}(\gamma_{r_m}, \rho) \Def f_{m, j} + \frac{f_{m, j+1} - f_{m, j}}{\rho^{(n)}_{m, j+1} - \rho^{(n)}_{m,j}}(\rho - \rho^{(n)}_{m, j}),
\end{equation}
for all $\rho \in [\rho^{(n)}_{m, j}, \rho^{(n)}_{m ,j+1}] \subset \mathcal G^{(n)}_m$ (defined as in Notation \ref{N:proj_FD_m}) and 
\begin{equation*}
    f_{m, j} \Def f(\gamma_{r_m}, \rho^{(n)}_{m, j})= f^{(n)}(\gamma_{r_m}, \rho^{(n)}_{m, j}).
\end{equation*}

\begin{notation}\label{N:PP_grid_points}
    Since the distance between the grid points is not the same, to address the corresponding distance of a grid point $z$ to the preceding and proceeding points, we define $\delta^{(n)}_-(z)$ and $\delta^{(n)}_+(z)$, respectively. 
\end{notation}
Using the grid point construction, the initial value function $\rho_\circ$ can be approximated by simple functions of the form  
\begin{equation}
    \label{E:iv_approx}
    \rho_\circ^{(n)}(x) \Def \sum_{m = 0}^{M} \sum_{\substack{j = - k^{(n)}_m \\ j \in \ZZ}}^{k^{(n)}_m -1} \rho_{m, j}^{(n)} \bOne_{E_{m, j}^{(n)}}(x), 
\end{equation}
where, $k_m^{(n)} = \gamma_{r_m} 2^n$ (see Remark \ref{R:piecewise_gamma}), $\rho^{(n)}_{m, j} \in \mathcal G^{(n)}_m$ and 
\begin{equation*}
    E_{m, j}^{(n)} \Def \rho_\circ^{-1} \left([\rho_{m, j}^{(n)}, \rho_{m, j + 1}^{(n)}) \right) \cap I_m.
\end{equation*}
which are disjoint sets. 
Using the bounded variation and measurability of $\rho_\circ$, for any fixed $n \in \NN$ and $m \in \set{0, \cdots, M}$ and by possibly rearranging on a set of Lebesgue measure zero, the approximation $\rho^{(n)}_\circ$ can be written in the form of 
\begin{equation}\label{E:piecewise_approx}
    \begin{split}
       & \rho^{(n)}_\circ(x) \bigmid_{x \in I_m} = \begin{cases} \rho^{(n)}_{m,j_1} &, x\in [x_\circ, x_1) \\
        \rho^{(n)}_{m, j_2} &, x \in [x_1, x_2) \\
        \vdots \\
        \rho^{(n)}_{m, j_{N_m}} &, x \in [x_{N_m -1}, \sa_{m + 1}) \end{cases}
        \end{split}\end{equation}
for $m \in \set{0, \cdots, M}$ where $\rho^{(n)}_{m, j_r} \in \mathcal G^{(n)}_m$, for $r \in \set{1, \cdots, N_m}$. 
\begin{lemma}\label{L:initial_convergence} For $\rho_\circ \in L^1_{\loc}(\RR)$, we have that
\begin{equation*}
    \rho_\circ^{(n)} \to \rho_\circ, \quad \text{pointwise, and $L^1_{\loc}(\RR)$}.
\end{equation*}
\end{lemma}
\begin{proof}
First, we note that since $\tfrac 12 \in \mathcal G^{(n)}_m$, for any $m \in \set {0, \cdots, M}$ and hence
\begin{equation}\label{E:order}
    \rho^{(n)}_{m, j},  \rho^{(n)}_{m, j+1} \le \tfrac 12 , \quad \text{or $\rho^{(n)}_{m, j},  \rho^{(n)}_{m, j+1} \ge \tfrac 12$}.
\end{equation}
In addition, by \eqref{E:iv_approx} we have that for any $x \in \RR$ there exists $m, j$ such that $x \in E_{m,j}^{(n)}$ and therefore,
\begin{equation*}
    \rho_\circ(x) - \rho_\circ^{(n)}(x) \le \rho^{(n)}_{m, j+1} - \rho^{(n)}_{m, j}.
\end{equation*}
Equation \eqref{E:order} implies that either $j, j+1 <0$ or $j, j+1 \ge 0$. Suppose first, they are positive ($j \ge 0)$. Then, we have that
\begin{equation*}\begin{split}
    \rho_\circ(x) - \rho_\circ^{(n)}(x) &\le  \tfrac 12 \left(1 + \sqrt{\tfrac{\abs{j+1}}{2^{n + 2}\gamma_{r_m}}}\right) - \tfrac 12 \left(1 + \sqrt{\tfrac{\abs{j}}{2^{n +2} \gamma_{r_m}}}\right)\\
    & \le \left(\tfrac{1}{\sqrt{2^{n + 2}\gamma_{r_m}}} \right) \left(\tfrac{1}{\sqrt{j + 1} + \sqrt j} \right)\\
    &\le  \left(2^{n + 2}\gamma_{r_m}\right)^{-\tfrac 12} \to 0 , \quad \text{as $n \to \infty$}.
\end{split}\end{equation*}
The same argument holds true for the case $j \le -1$. 
This proves the pointwise convergence of the claim. 
Furthermore, $\rho^{(n)}_\circ \le \rho_\circ$ by construction. Therefore, $\abs{\rho^{(n)}_\circ - \rho_\circ} \le 2 \rho_\circ \in L^1_{\loc}(\RR)$ which implies the $L^1_{\loc}(\RR)$ convergence.
\end{proof}
Let $(\gamma(x), z_\circ(x)) \in \mathcal W$ with
\begin{equation*}
    z_\circ(x) \Def \psi(\gamma(x), \rho_\circ(x)) ,\quad x \in \RR.
\end{equation*}
Then, the convergence $z_\circ^{(n)} \to z_\circ$ as $n \to \infty$ pointwise and in $L^1_{\loc}$ follows from the definition in a straightforward way. 
Furthermore, by construction of the approximate functions $\rho^{(n)}_\circ$ of $\rho_\circ$ and $z^{(n)}_\circ$ of $z_\circ$, we have that 
\begin{equation}
    \label{E:bounded_var_initial}
    \totvar{\rho^{(n)}_\circ}{\RR } \le \totvar{\rho_\circ}{\RR }, \quad \totvar{z^{(n)}_\circ}{\RR} \le \totvar{z_\circ}{\RR }.
\end{equation}
\subsection{Wave Interactions and Bounded Variation}
By the above construction of the grid points, for any fixed $n \ge N_\circ$ the solution $(\rho^{(n)}, y_n)\in L^1_{\loc}(\RR_+ \times \RR; [0,1]) \times W^{1,1}_{\loc}(\RR_+;\RR)$ of the Riemann problem 
\begin{equation*}
    \begin{split}
       & \rho_t + \fpartial{}{x}[f^{(n)}(\gamma, \rho)] = 0, \quad x \in \RR, \, t \in \RR_+\\
       & \rho_\circ(x) = \begin{cases} \rho^{(n)}_{m-1,l} &, x<0 \\
        \rho^{(n)}_{m, r} &, x>0 \end{cases} , \quad \gamma(x) = \begin{cases} \gamma_{r_m-1} &, x<0 \\
        \gamma_{r_m} &, x>0 \end{cases}\\
        & f^{(n)}(\gamma_{r_m}, \rho(t, y_n(t)))-\dot y_n(t) \rho(t, y_n(t)) \le F_\alpha^{(m)}(\dot y_n(t)),
    \end{split}
\end{equation*}
where, $t \mapsto y_n(t)$ is the solution of the ODE
\begin{equation*}\begin{split}
    \dot y(t) &= w(y(t), \rho^{(n)}(t, y(t)+)) = \min \set{V_b^{(m)},v(\gamma_{r_m}, \rho^{(n)}(t, y(t)+))}\\
    y_\circ &= 0 \in I_{r_{m-1}} ,
\end{split}\end{equation*}
for some $m \in \set{0, \cdots, M}$ and where  $\rho^{(n)}_{m-1,j}, \rho^{(n)}_{m, r} \in  \mathcal G^{(n)}_m$, is well-defined and is obtained from Definition \ref{def:Riemann_sol_new}.
In the presence of a discontinuity in the flux, the uniform bounded variation can be violated in the $\mathcal V$ space, \cite{temple1982global}. To resolve the issue, we employ a functional (known as Temple functional, \cite{temple1982global}) which is employed to prove the bounded variation in the $\mathcal W$ space. 

In this paper, however, we need to customize a suitable Temple functional due to the existence of the moving bottleneck and in particular a non-classical shock. Let's start settling on some notations which will be used throughout the section. 
The function $(\gamma, \rho) \mapsto f(\gamma, \rho)$ will be approximated by a piecewise continuous function $(\gamma, \rho) \mapsto f^{(n)}(\gamma, \rho)$ and the initial condition $x \mapsto \rho_\circ(x)$ by a piecewise constant function $x \mapsto \rho^{(n)}_\circ(x)$ as in \eqref{E:iv_approx}.
The solution of the approximate Riemann problem is denoted by $(t, x) \mapsto (\rho^{(n)}(t, x), y_n(t))$. In addition, the corresponding solution in $\mathcal W$ space is denoted by $z^{(n)}$ and is defined by 
\begin{equation}\label{E:z_piecewise_sol}
    z^{(n)}(t, x) \Def \psi(\gamma(x), \rho^{(n)}(t, x)). 
\end{equation}

Using \eqref{E:piecewise_approx}, in the rest of this section we will show that the solution to the following constrained Riemann problem is well-defined and satisfies some compactness properties. 
\begin{equation}
    \label{E:discrete_Riemann_multiple}
    \begin{split}
       & \rho_t + \fpartial{}{x}[f^{(n)}(\gamma, \rho)] = 0\\
       & \rho^{(n)}_\circ(x) = \begin{cases} \rho^{(n)}_{\circ,1} &, x\in [x_\circ, x_1) \\
        \rho^{(n)}_{\circ, 2} &, x \in [x_1, x_2) \\
        \vdots \\
        \rho^{(n)}_{\circ, n_\circ} &, x \ge x_{n_\circ -1} \end{cases} , \quad \gamma(x) = \begin{cases} \gamma_{r_\circ} &, x \in I_\circ \\
        \gamma_{r_1}&, x \in I_1\\
        \vdots \\
        \gamma_{r_M} &, x \in I_M \end{cases}\\
        & f^{(n)}(\gamma(y_n(t), \rho(t, y_n(t)))-\dot y_n(t) \rho(t, y_n(t)) \le F_\alpha(y_n(t), \dot y_n(t))
    \end{split}
\end{equation}
where $y_n$ is the solution of the bottleneck dynamic
\begin{equation}\label{E:ODE_BN}
    \dot y(t) = \min \set{V_b(y(t)), v(\gamma(y(t)), \rho^{(n)}(t, y(t)+))} , \quad y(0) = y_\circ.
\end{equation}
Here, $\rho^{(n)}_{\circ, r} \in \mathcal G^{(n)}, \, r \in \set{1, \cdots, n_\circ}$ (the grid points). From here on, the model defined by \eqref{E:discrete_Riemann_multiple} and \eqref{E:ODE_BN} is called the \textbf{n-Approximate Problem}. 

In order to understand the solution to such a Riemann problem, we need to investigate the interaction between different types of waves.
Under the definitions of our problem, the $\gamma$-fronts have zero speed of propagation and hence they do not collide. In addition, since the bottleneck belongs to one region at a time, two or more non-classical $z$-fronts (trajectory of the bottlenecks) also do not interact with each other. Furthermore, the result of the collision of any two or more classical $z$-fronts will be merely one $z$-front. Therefore, we only need to study the interactions of non-classical $z$-fronts with classical $z$-fronts, and $z$-fronts (classical and non-classical) with $\gamma$-fronts. 

\noindent Let $\front$ denote an individual front (discontinuity). We define a (Temple) functional
\begin{equation}\label{E:Temple}
    \temple(\front) \Def
    \begin{cases}
    \abs{\Delta z}   &, \text{if $\front$ is a $z$-front}\\
    \abs{\Delta \gamma} &, \text{if $\front$ is a $\gamma$-front with $z_L < z_R$}\\ 
    \tfrac 12 \abs{\Delta \gamma}   &, \text{if $\front$ is a $\gamma$-front with $z_L > z_R$}
    \end{cases}
\end{equation}
where, $\Delta z \Def z_R - z_L$ and $z_L$ and $z_R$ are the left and right states of discontinuity $\front$, $\hat z_m$, $\Delta \gamma = \gamma_R - \gamma_L$ and $\check z_m$ corresponds to $\hat \rho_m$ and $\check \rho_m$, respectively. In addition, we recall that by \eqref{E:z_gamma_relation} if $\front$ is a $\gamma$-front, 
\begin{equation*}
    \abs{\Delta z} = \tfrac 14 \abs{\Delta \gamma}.
\end{equation*}

Now for a family of fronts, we need to generalize the definition of the Temple functional. Let's fix $\bar \ft \in \RR_+$. In addition, suppose the solution $z^{(n)}(\bar \ft, x)$ in $(z, \gamma)$-space corresponding to the solution $(\rho^{(n)}(\bar \ft, x), y_n(\bar \ft))$ to the approximate problem \eqref{E:discrete_Riemann_multiple} exists (we will discuss this in details later in this paper). We define 
\begin{equation}
    \label{E:Temple_general}
    \begin{split}
    \temple(z^{(n)}(\bar \ft, \cdot)) &\Def \sum_{\front \in \mathcal F(\bar \ft)} \temple(\front) +  \varpi(\bar \ft)
\end{split}\end{equation}
where, $\mathcal F(\bar \ft)$ is the collection of all fronts at time $\bar \ft$. To define the function $\varpi$, let's fix $\varsigma>0$, a sufficiently small scalar. Then
\begin{enumerate}
    \item If $y_n(\bar \ft) \in I_m=[\sa_m, \sa_{m+1})$ for some $m \in \set{0, \cdots, M}$, and $\abs{y_n(\bar \ft) - \sa_{m +1}} > \varsigma$, then
\begin{equation}\label{E:Gamma_1}
\varpi(\bar \ft) \Def \begin{cases} 0 &, \text{if $\rho^{(n)}(\bar \ft, y_n(\bar\ft)-)= \hat \rho_m$ and} \\
& \quad \text{$\rho^{(n)}(\bar \ft, y_n(\bar \ft) +) = \check \rho_m$} \\
2(\hat z_m - \check z_m)  &, \text{otherwise} \end{cases}
\end{equation}
\item If $y_n (\bar \ft) \in I_m$ for some $m\in \set{0, \cdots, M}$, and $\abs{y_n(\bar \ft) - \sa_{m + 1}}< \varsigma$, then
\begin{equation}\label{E:Gamma_2}
\varpi(\bar \ft) \Def 2 (\hat z_{m + 1} -\check z_{m + 1}).
\end{equation}
\begin{remark}
    In the definition of $\varpi(\bar \ft)$ we are mainly interested in $\varsigma \to 0+$. In particular, if collision time is denoted by $\ft_\circ$, then \eqref{E:Gamma_1} covers the case in which the interacting fronts will remain in the same region $I_m$ at $\ft_\circ^-$ and $\ft_\circ^+$ and the case for which the states are in $I_{m+1}$ at $\ft_\circ^+$ and in $I_m$ at $\ft_\circ^-$ will be explained by \eqref{E:Gamma_2}. 
\end{remark}

\end{enumerate}

Then, from \eqref{E:Temple}, \eqref{E:Temple_general}, \eqref{E:Gamma_1} and \eqref{E:Gamma_2} we have that for any $t \ge0$
\begin{equation}
    \label{E:TV_Temple}
    \totvar{z^{(n)}(t, \cdot)}{\RR } \le\temple(z^{(n)}(t, \cdot)) \le  \totvar{z^{(n)}(t, \cdot)}{\RR} + \totvar{\gamma(\cdot)}{\RR} + C_{\eqref{E:TV_Temple}},
\end{equation}
where
\begin{equation*}
    C_{\eqref{E:TV_Temple}} \Def 2 \max_{m \in \set{0, \cdots, M}} (\hat z_m - \check z_m),
\end{equation*}
is independent of $n$ and $t$.
The left inequality in \eqref{E:TV_Temple} is clear by the definition of the Temple functional. The right inequality can also be deduced by the definition of Temple functional and the fact that $\max_{m} \left(\hat z_m - \check z_m \right)$ is finite. 

\begin{figure}
    \centering
    \includegraphics[width=2.9in]{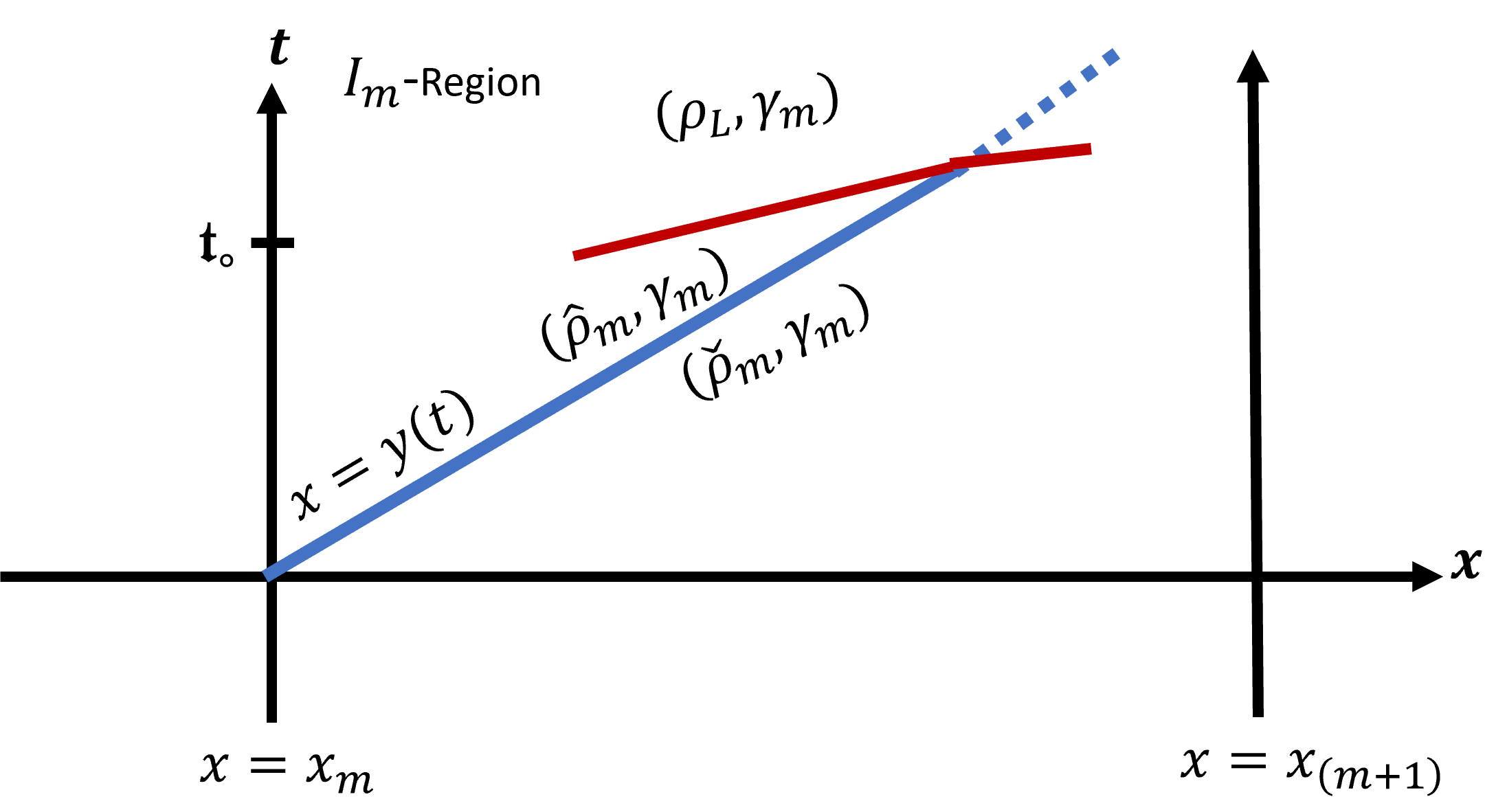}
    \includegraphics[width=2.9in]{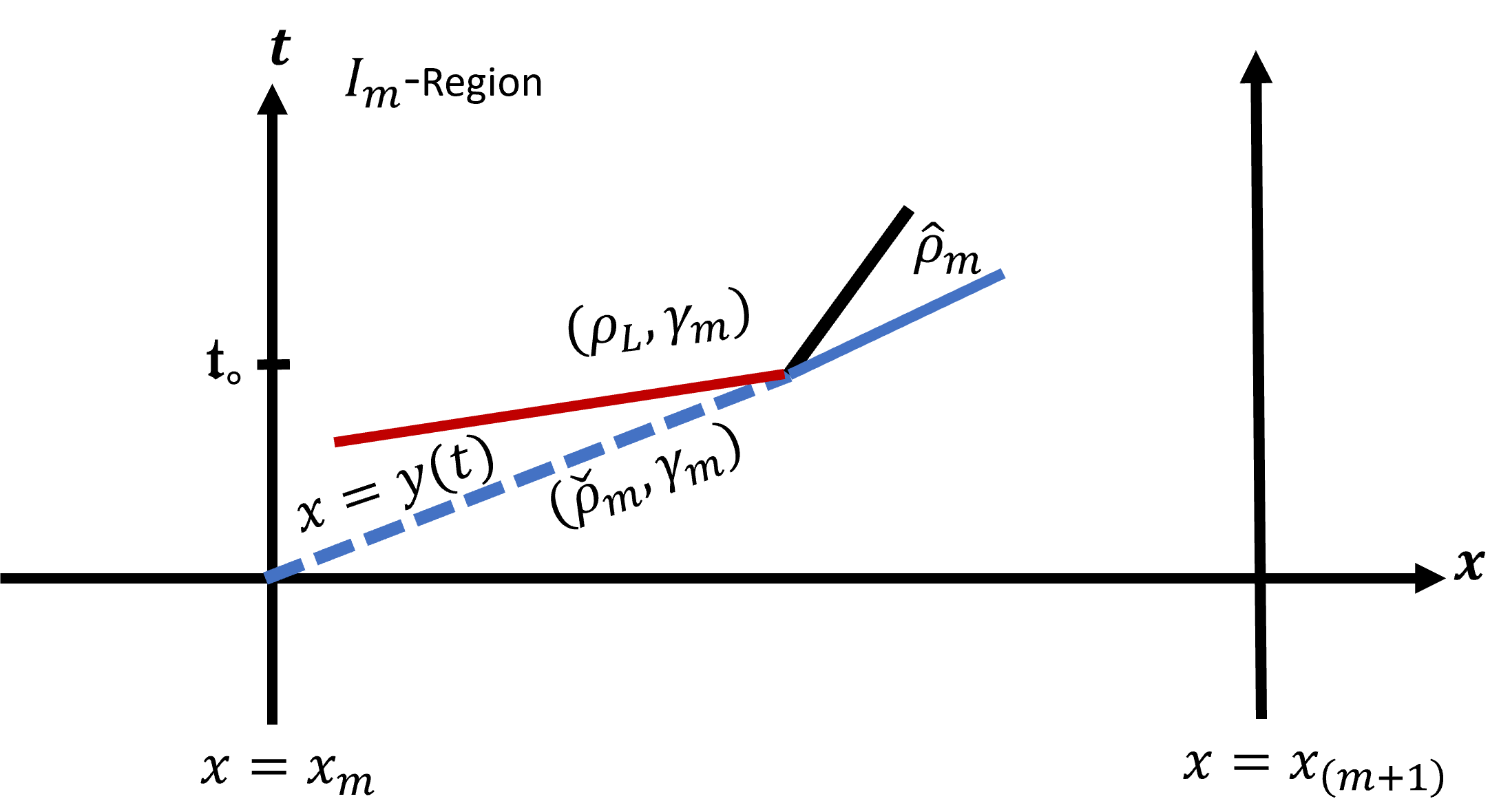}
    \caption{Instances of interaction between the bottleneck and $z$-fronts. In the left case, $\rho_L \in [0, \check \rho_m)$ and in the right case $\rho_L > \check \rho_m$. }
    \label{fig:interaction_88_89}
\end{figure}

\begin{figure}
    \centering
    \includegraphics[width=3in]{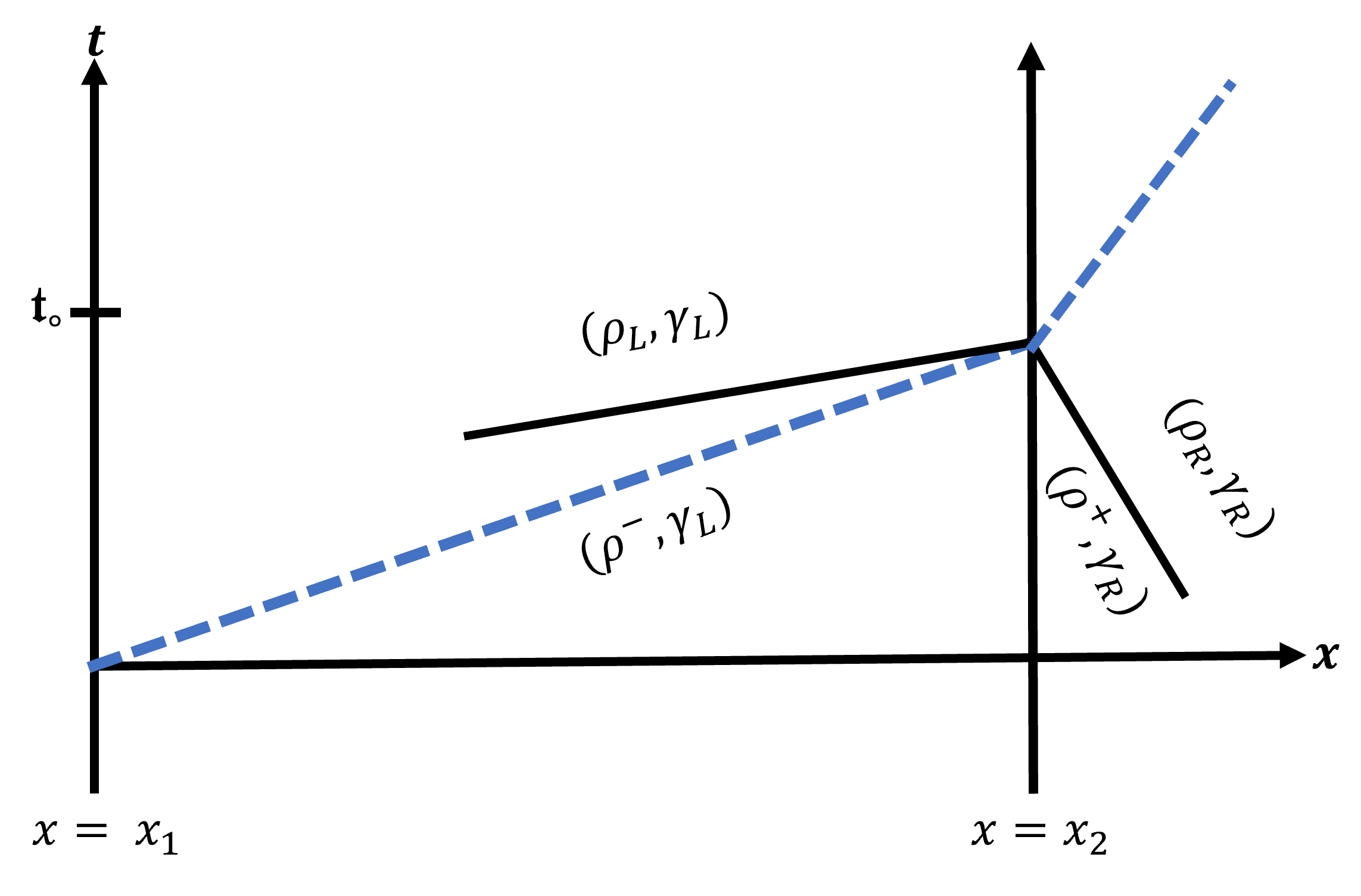}
    \caption{Interaction of fronts with $\gamma$ wave at time $\ft_\circ$. The dashed lines represent the bottleneck trajectories in each region.}
    \label{fig:interaction_87}
\end{figure}

In the following subsections, we will investigate the interaction of a $z$-front and a $\gamma$-front, non-classical $z$-front and $\gamma$-front, and classical and non-classical $z$-fronts to determine the state of the solution after collisions with other fronts. In addition, in each case, we show that the Temple function is decreasing. 
\paragraph{\textbf{A. Interaction of Waves and the Bottleneck Trajectory Inside a Region}}
Let's fix $m \in \set{0, \cdots, M}$ and consequently $\gamma_{r_m}$ and region $I_m$. In the interior of this region, $I_m^\circ$, all the interactions between the bottleneck trajectory $z$-fronts and also between other $z$-fronts follow from the classical case (see \cite{delle2014scalar, liard2021entropic} for the details of the wave interactions in this case). Figure \ref{fig:interaction_88_89} shows two possible interactions between the bottleneck trajectory and the $z$-fronts.
Using the definition of \eqref{E:Temple_general}, we can calculate that $\temple(z^{(n)}(\ft_\circ^+, \cdot)) - \temple(z^{(n)}(\ft_\circ^-, \cdot))$ for the left instance will be zero and for the right one will be less than $-2 \underline \delta^{(n)}$. The Temple functional shows the same decreasing behavior for all other cases. 
\paragraph{\textbf{B. Interaction of Waves with Bottleneck Trajectory Between Different Regions}}
In this work (in comparison with the conservation law without the flux constraint) there are relatively distinct potential states that can happen as the result of the wave collisions on the boundary of regions $I_m$ for all $m$ (also known as collision with $\gamma$-waves); see Figure \ref{fig:interaction_87}. Particular interest is in the collision of non-classical shock with the $\gamma$ fronts. 

To fully understand the interaction of the waves and in particular to ensure the decreasing behavior of the Temple functional in all cases, one needs to investigate all possibilities depending on the values of $\rho_L, \rho^-, \rho^+$ and $\rho_R$. More importantly, in the case of non-classical shock in any of the regions, the locations of $\check \rho_L, \hat \rho_L, \check \rho_R$ and $\hat \rho_R$ contribute to the creation of various cases that need to be studied (For instance, when $\rho_L = \hat \rho_L$ and $\rho^- = \check \rho_L$ in Figure \ref{fig:interaction_87} and the bottleneck creates a non-classical shock $\front[\hat \rho_R, \check \rho_R]$ after the collision with the boundary). 

\begin{figure}
    \centering
    \includegraphics[width=5in]{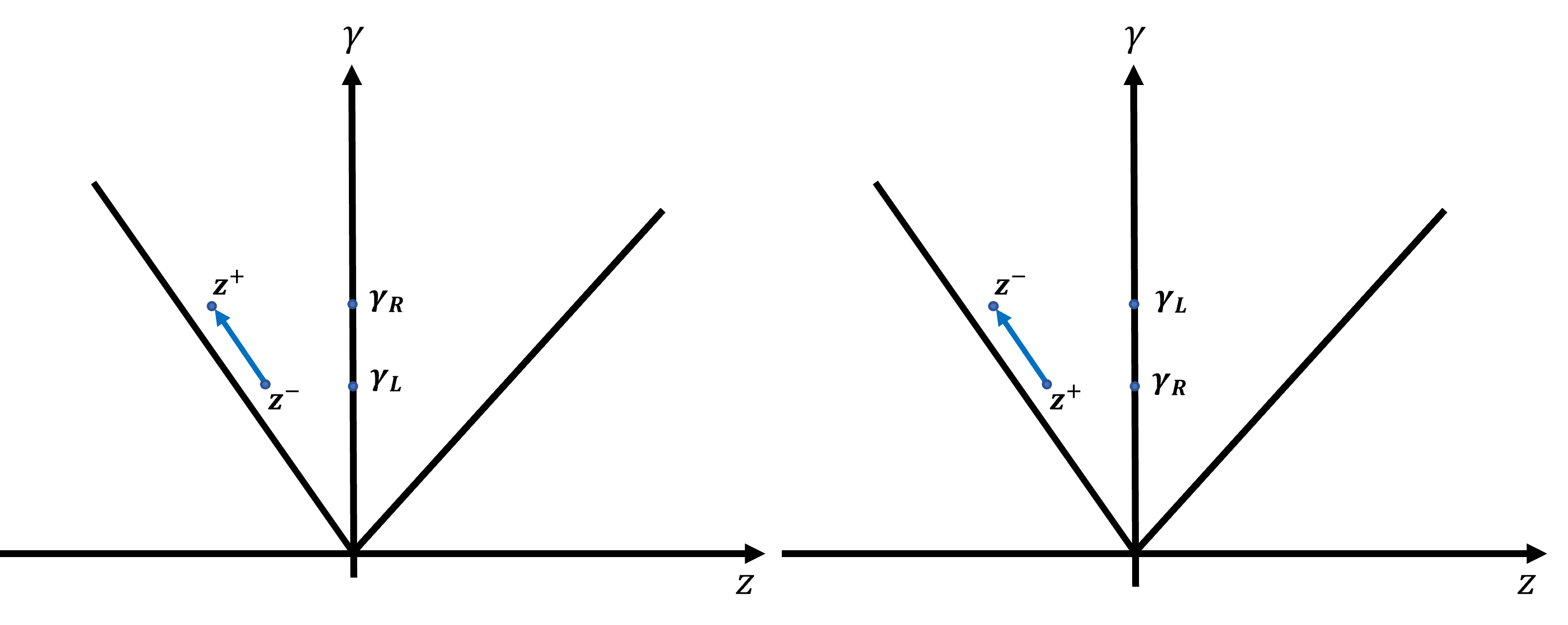}
    \caption{Front $\front[z^-, z^+]$ is negative. The categorization of the Riemann solution is based on the initial location of $z^-$ and $z^+$. Similarly, the front $\front[z^-, z^+]$ can be positive.}
    \label{fig:interaction_107}
\end{figure}

More precisely, all possible interactions can be enlisted in one of the four categories which are defined based on $\gamma_L \lessgtr \gamma_R$ and that the $\gamma$ front $\front[z^-, z^+]$ is located in the positive quadrant (both states are positive) or the negative one (both states are negative); see Figure \ref{fig:interaction_107}. See also Lemma \ref{L:gamma_jump_state_closedness} for justifying that these are the only possible cases. When the locations of $z^-$ and $z^+$ are determined, one may look for the admissible range of $z_L$ and $z_R$ (consult Figure \ref{fig:interaction_87} and Figure \ref{fig:interaction_90_91} and the detail will be elaborated in the following cases). In addition, the acceptable range of $\check \rho_R$ and $\hat \rho_R$ and consequently the Riemann solution $\Riemann^\alpha(\rho_L, \rho_R; \gamma_L, \gamma_R)(\cdot)$ will be identified. Finally, the decreasing behavior of the Temple functional will be concluded. 

In the following, we will study some of the most important cases which do not happen in the classical case, in depth. We will explain these cases in detail to provide insight into the complexities that the presence of the bottleneck can lead to.

\begin{case} Let's consider $\gamma_R > \gamma_L$, and front $\front[z^-, z^+]$ is negative with $\rho_L = \hat \rho_L$ and $\rho^- = \check \rho_L$ (see the left illustration of Figure \ref{fig:interaction_90_91}). In particular, we consider the non-classical shock hitting a $\gamma$ front from the left and a $z$ front hitting the $\gamma$ front from the right. This case, in particular, creates more complex conditions since as opposed to a classical rarefaction, the distance between $\check z_L$ and $\hat z_L$ are not necessarily limited to $\hat \delta^{(n)}$ and hence various new cases may happen in general. 

Next, we need to find the possible locations of $\hat z_L$ and $z_R$ to find the the Riemann solution $\Riemann^\alpha(\hat \rho_L, \rho_R; \gamma_L, \gamma_R)(\cdot)$. Considering that the slope of the non-classical shock should be positive and that of front $\front[\rho^+, \rho_R]$ negative, the admissible range of $\hat z_L$ and $z_R$ will be 
\begin{equation}\label{E:range}
    \hat z_L \in [\check z_L + \delta^{(n)}_+(\check z_L), - \check z_L - \delta^{(n)}_-(\check z_L)], \quad z_R \in [- \zplus + \delta^{(n)}_+(\zplus), \tfrac14 \gamma_R] \cup \set{\zplus}
\end{equation}
where, $\delta^{(n)}_{\pm}(\cdot)$ is defined in Notation \ref{N:PP_grid_points}.
\begin{figure}
    \centering
    \includegraphics[width=2.9in]{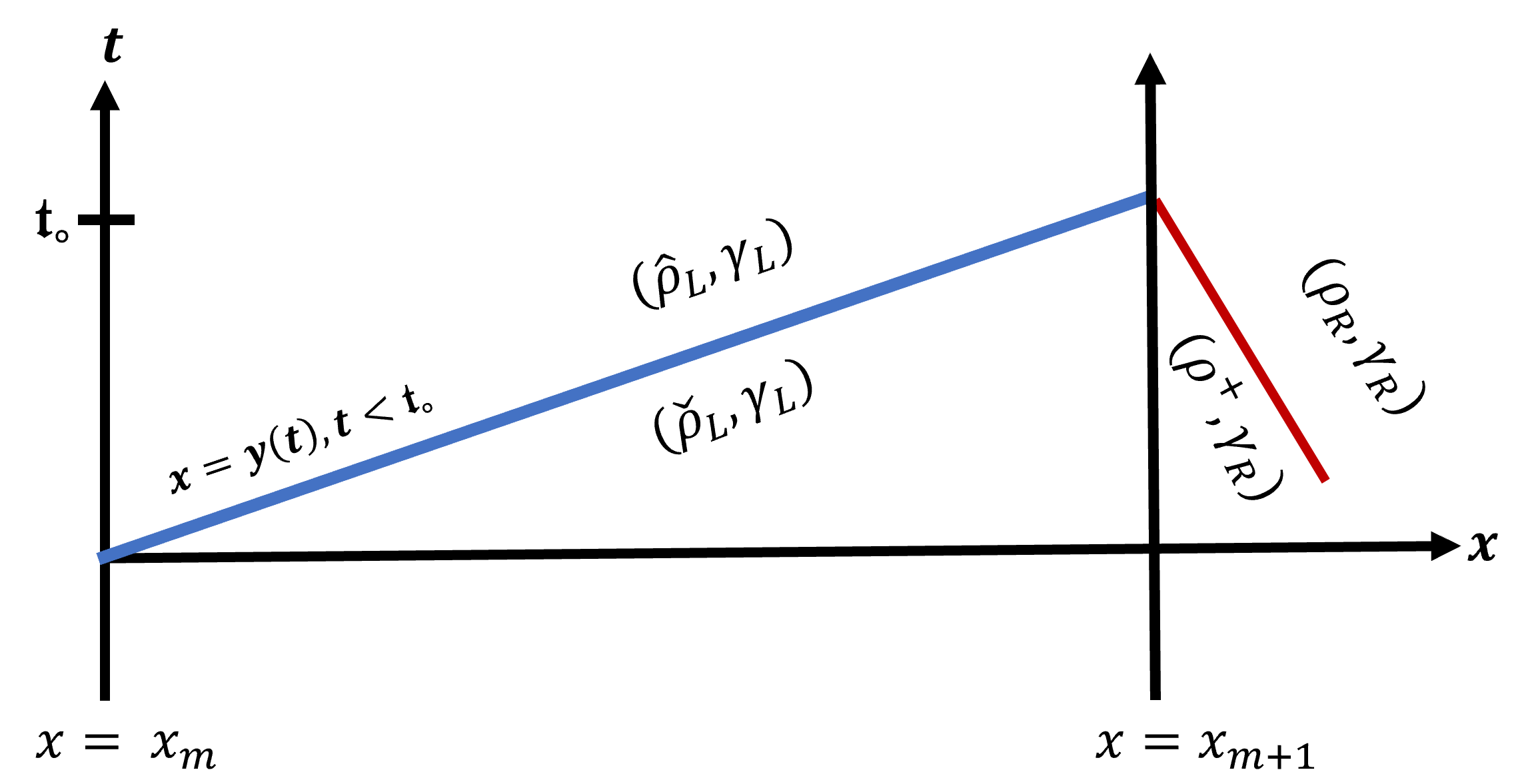}
    \includegraphics[width=2.9in]{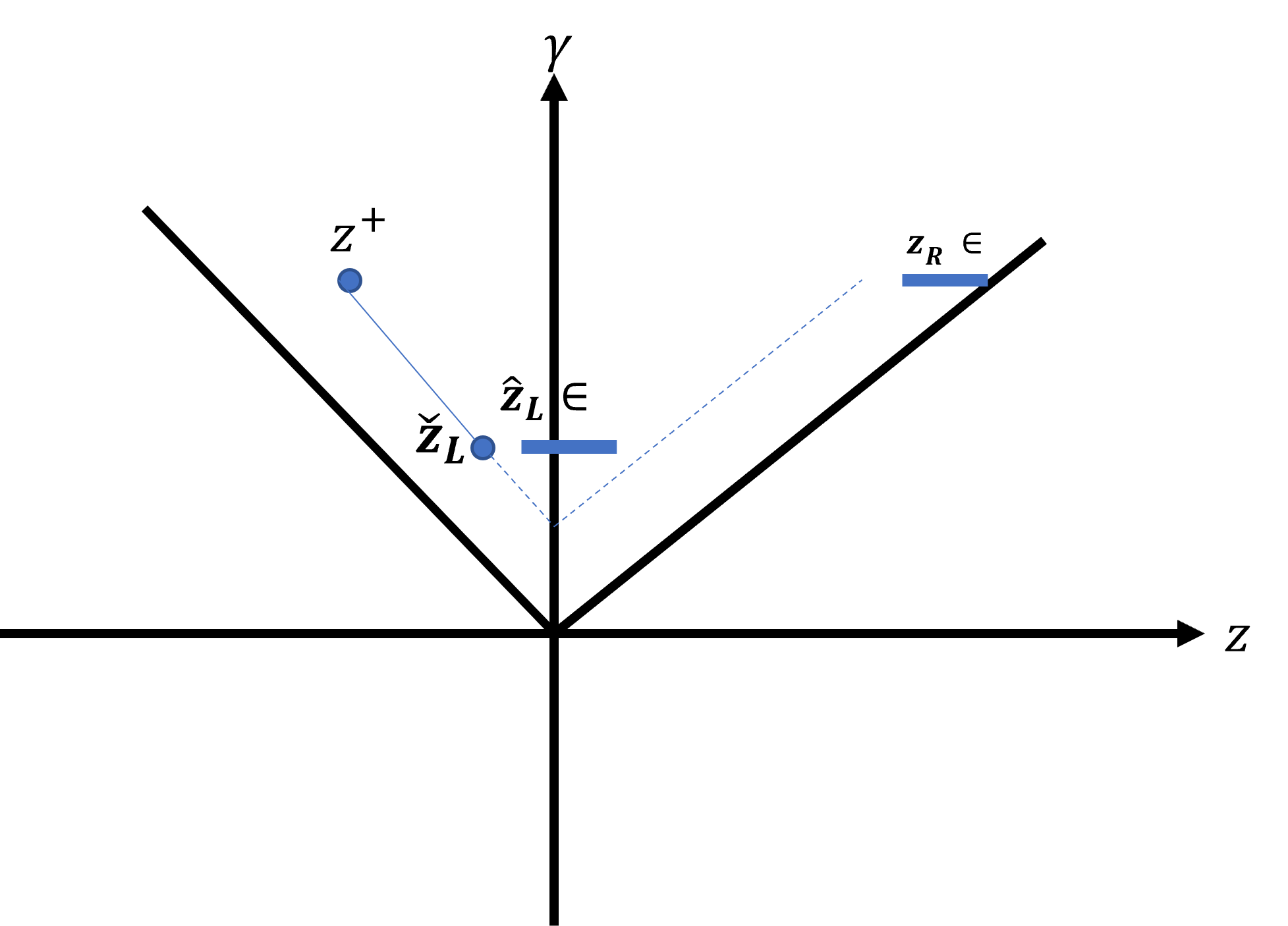}
    \caption{The left illustration shows the collision of a non-classical shock with a $\gamma$ front. The right illustration shows the admissible range of $\hat z_L$ and $z_R$.}
    \label{fig:interaction_90_91}
\end{figure}
The right illustration of Figure \ref{fig:interaction_90_91} shows the intervals \eqref{E:range}. It is worth noting that without the presence of the bottleneck, the admissible range would be $z_L \in [- \tfrac 14 \gamma_L, \zminus + \delta^{(n)}_+(\zminus)]$, with a null intersection with the case of non-classical shock. This implies that in the presence of the bottleneck, the composition of the solution and the variety of cases that can arise could be fundamentally different from the classical case. 

In this case, for all admissible $\hat z_L$, $f(\gamma_L, \check \rho_L) > f(\gamma_R, \rho_R)$ which means that the solution can be categorized as $z\gamma$-type (first a $z$ front and then a $\gamma$ front present the solution). Figure \ref{fig:interaction_92_93} illustrates the Riemann solution  $\Riemann^\alpha(\hat \rho_L, \rho_R; \gamma_L, \gamma_R) = \Riemann(\hat \rho_L, \rho_R; \gamma_L, \gamma_R)$ (when the bottleneck constraint is satisfied by the solution $\Riemann$) before and after the collision. Using these illustrations, we may show that $\temple(z^{(n)}(\ft_\circ^+, \cdot)) - \temple(z^{(n)}(\ft_\circ^-, \cdot)) = - 2 (\hat z_L - \check z_L) \le -2 \underline \delta^{(n)}$. 
\begin{figure}
    \centering
    \includegraphics[width=2.9in]{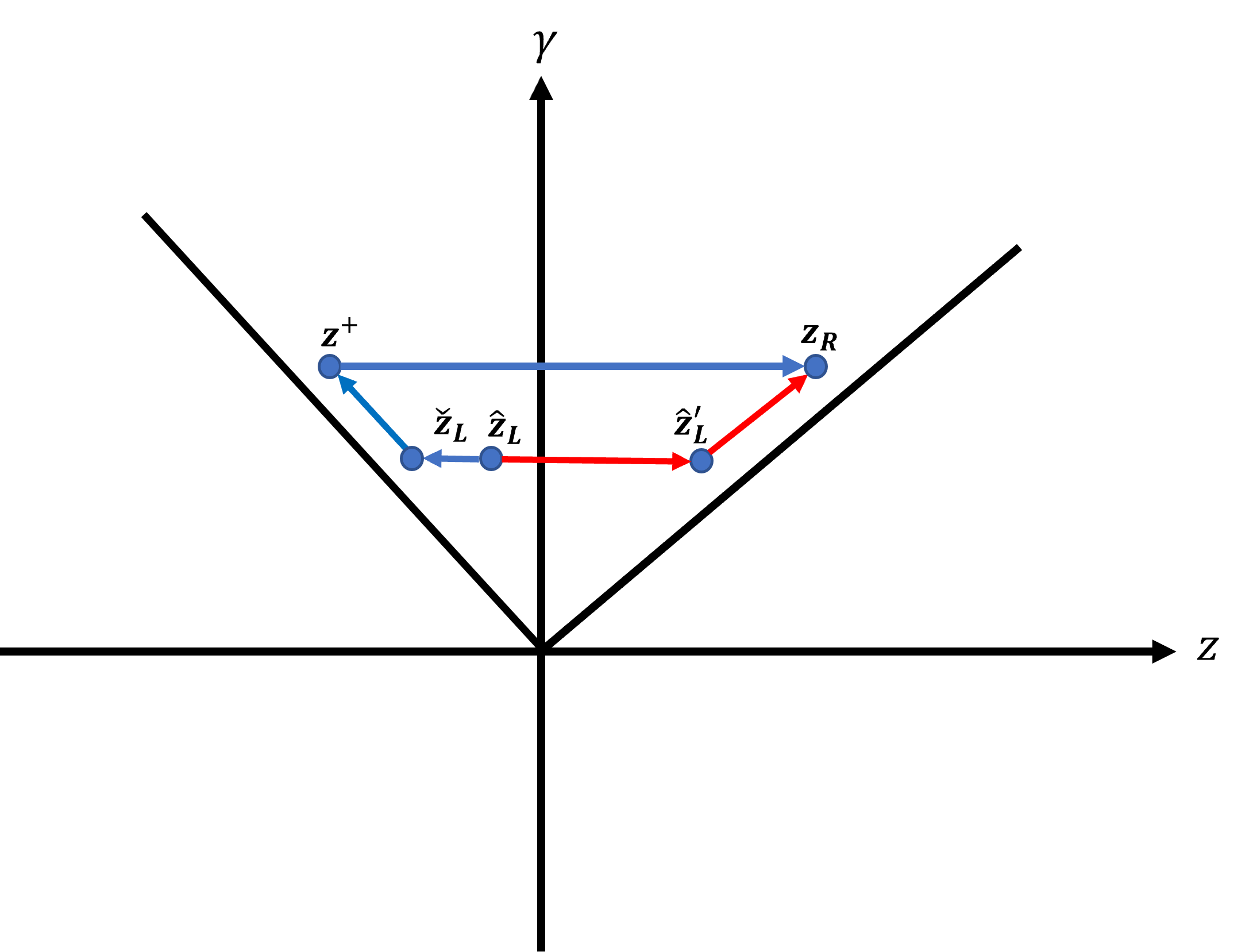}
    \includegraphics[width=2.9in]{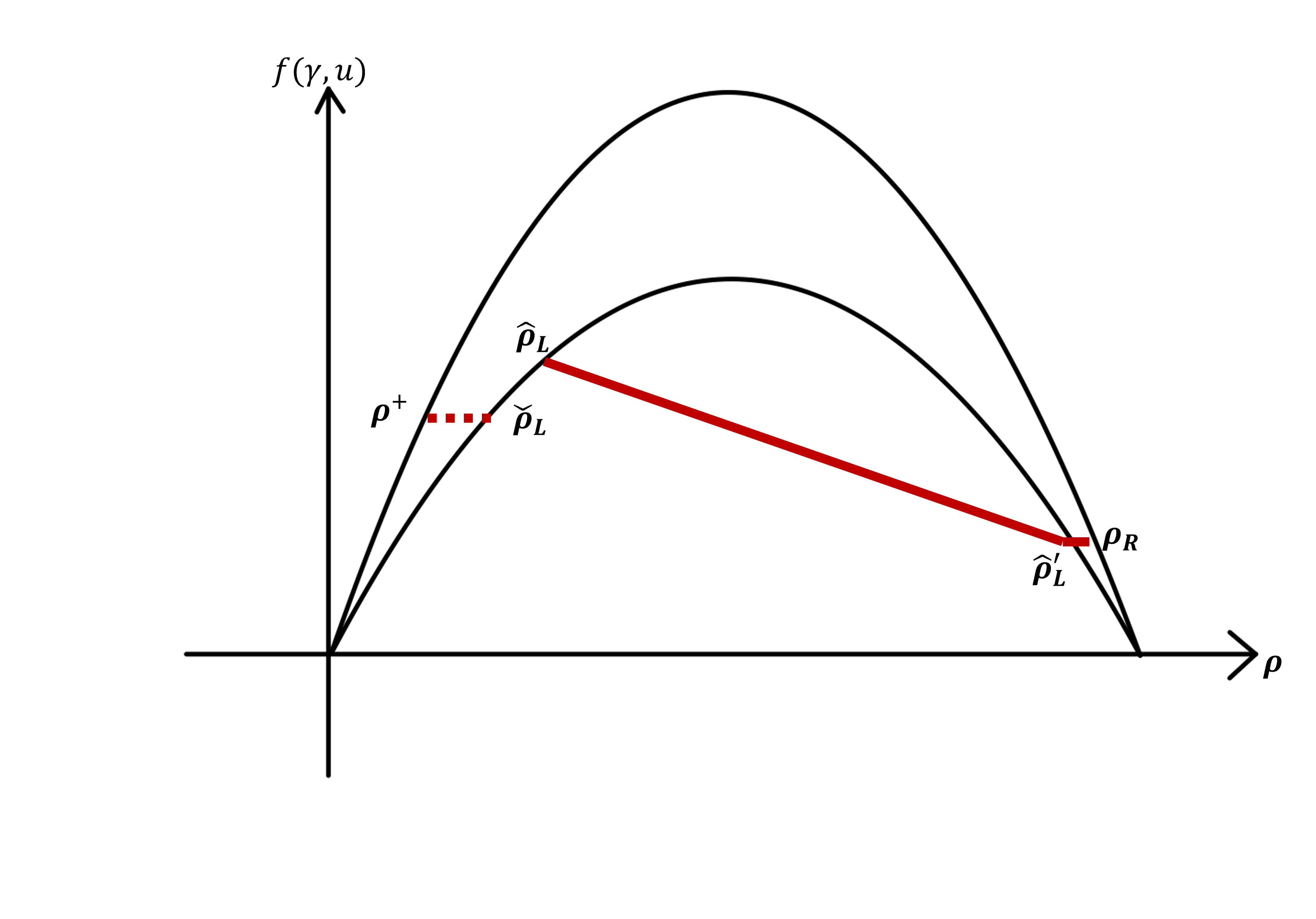}
    \caption{The Riemann solution $\Riemann(\hat \rho_L, \rho_R; \gamma_L, \gamma_R)(\cdot)$ in $\mathcal W$ space and on the fundamental diagram. In the left figure, the red lines illustrate the Riemann solution after collision, while the blue lines show the solution before the collision.}
    \label{fig:interaction_92_93}
\end{figure}
The Riemann solution $\Riemann^\alpha(\hat \rho_L, \rho_R; \gamma_L, \gamma_R) \ne \Riemann(\hat \rho_L, \rho_R; \gamma_L, \gamma_R)$ (i.e. the Riemann solution $\Riemann$ does not satisfy the bottleneck constraint) is illustrated in Figure \ref{fig:interaction_94_95} and \ref{fig:interaction_96}. 
\begin{remark} To find the Riemann solution $\Riemann^\alpha(\cdot)$, the values of $\check \rho_R$ and $\hat \rho_R$ should be determined (recall the Definition \ref{def:Riemann_sol_new}). Therefore, in general, the solution strongly depends on the values of $\check \rho_R$ and $\hat \rho_R$. This will be explained in more detail in Case \ref{case:2}. 
\end{remark}
In the case of this problem, however, by structure (the solution is of $z \gamma$-type), $\hat \rho_R > \rho_R$ and only one type of Riemann solution $\Riemann^\alpha(\hat \rho_L, \rho_R; \gamma_L, \gamma_R)$ follows. Moreover, for this Riemann solution, a similar decrease in Temple functional can be calculated.
\begin{figure}
    \centering
    \includegraphics[width=2.9in]{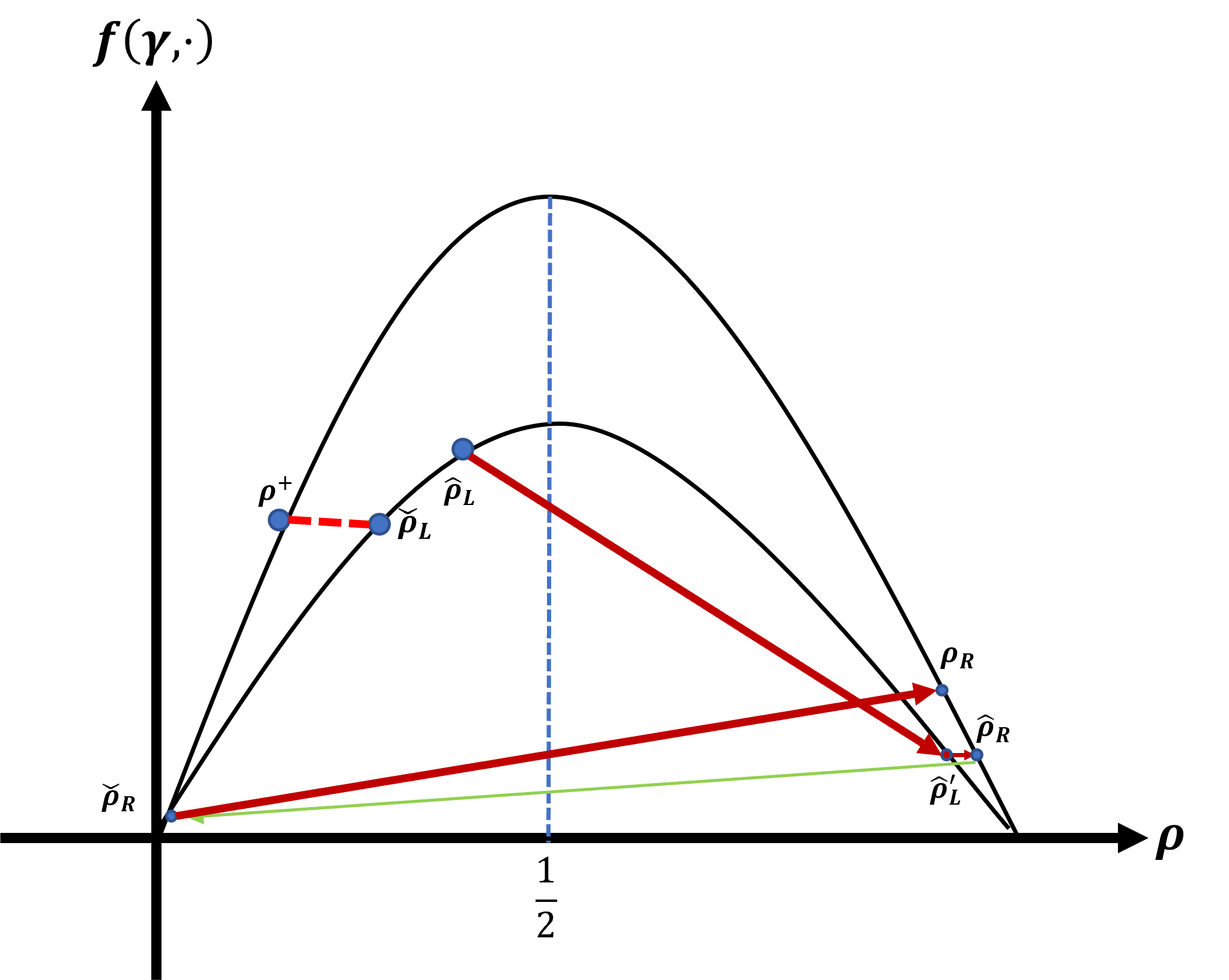}
    \includegraphics[width=2.9in]{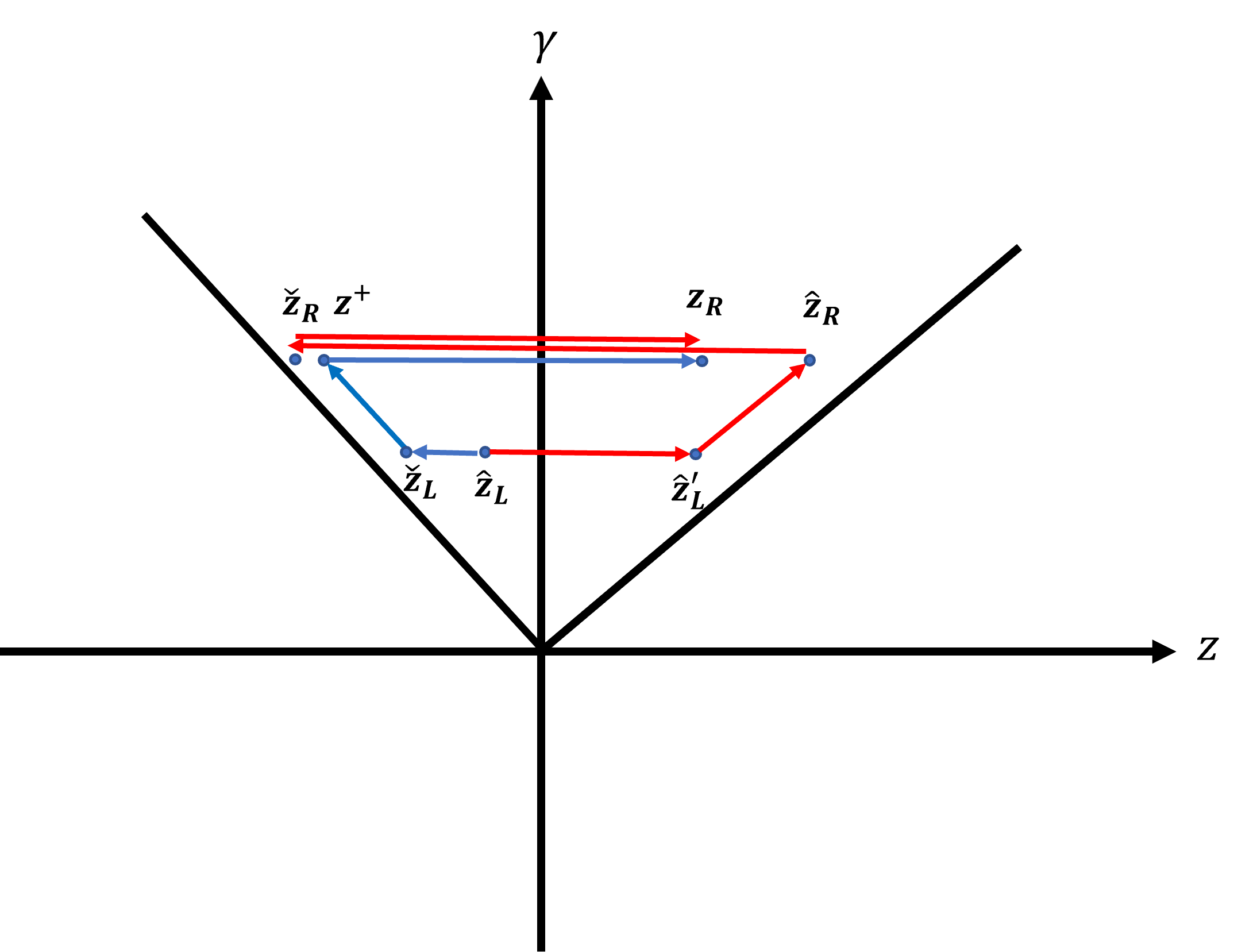}
    \caption{The Riemann solution $\Riemann^\alpha(\hat \rho_L, \rho_R; \gamma_L, \gamma_R)(\cdot)$ in $\mathcal W$ space and on the fundamental diagram. In the right figure, the red lines illustrate the Riemann solution after collision, while the blue lines show the solution before the collision.}
    \label{fig:interaction_94_95}
\end{figure}

\begin{figure}
    \centering
    \includegraphics[width=3in]{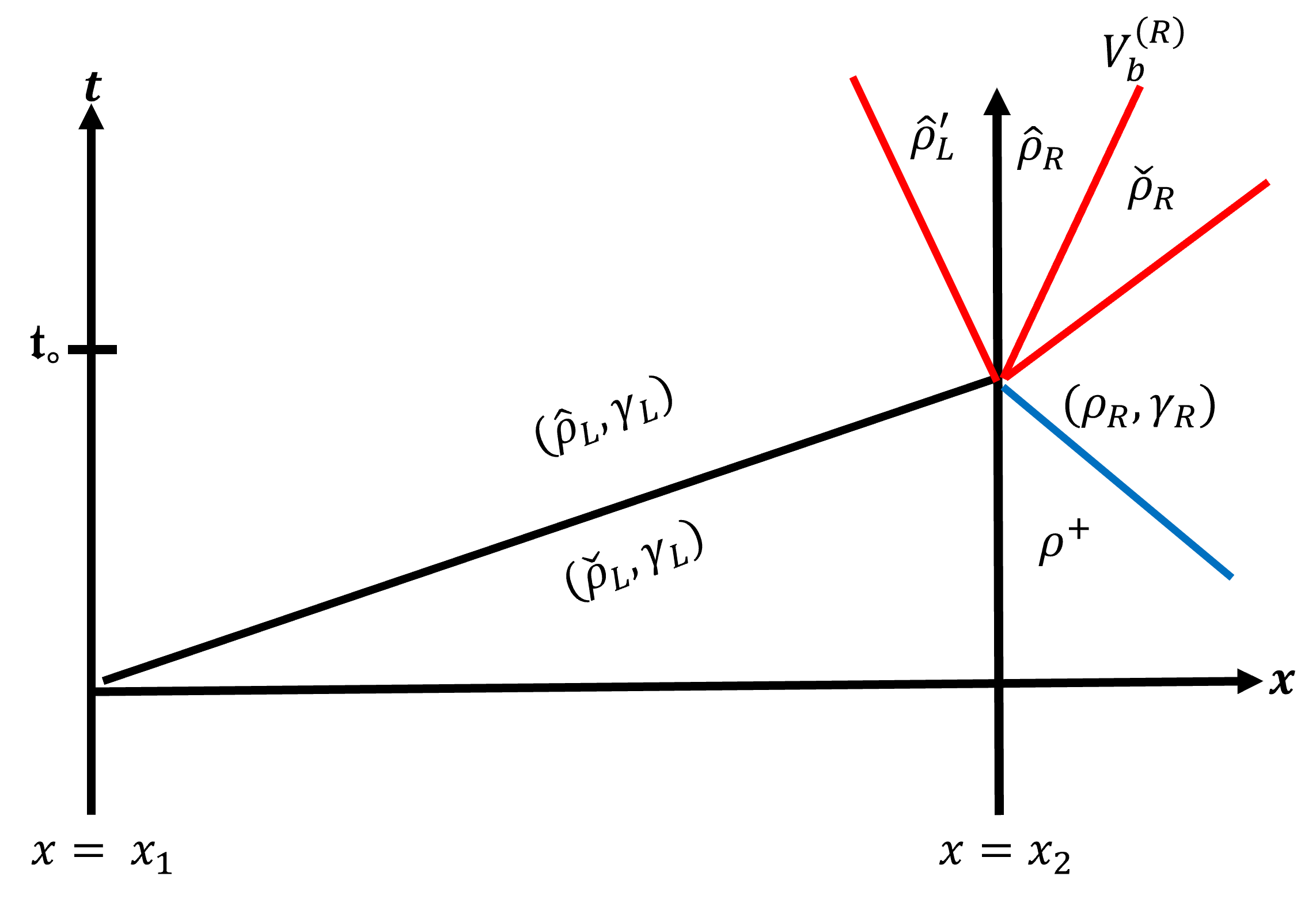}
    \caption{The Riemann solution $\Riemann^\alpha(\hat \rho_L, \rho_R; \gamma_L, \gamma_R)(\cdot)$. }
    \label{fig:interaction_96}
\end{figure}
\end{case}
\begin{case}\label{case:2}
For the next case, we consider $\gamma_R > \gamma_L$, the $\gamma$-front $\front[\rho^-, \rho^+]$ is negative, $\rho_L = \hat \rho_L$, $\rho^- = \check \rho_L$ and $\rho^+ = \rho_R$ (see the left illustrations of Figure \ref{fig:interaction_97_98}). The admissible range of $\hat z_L$, i.e. $\hat z_L \in [\check z_L + \delta^{(n)}_+(\check z_L) , -\check z_L - \delta^{(n)}_-(-\check z_L)]$, is illustrated in the right illustration of Figure \ref{fig:interaction_97_98}.
\begin{figure}
    \centering
    \includegraphics[width=3in]{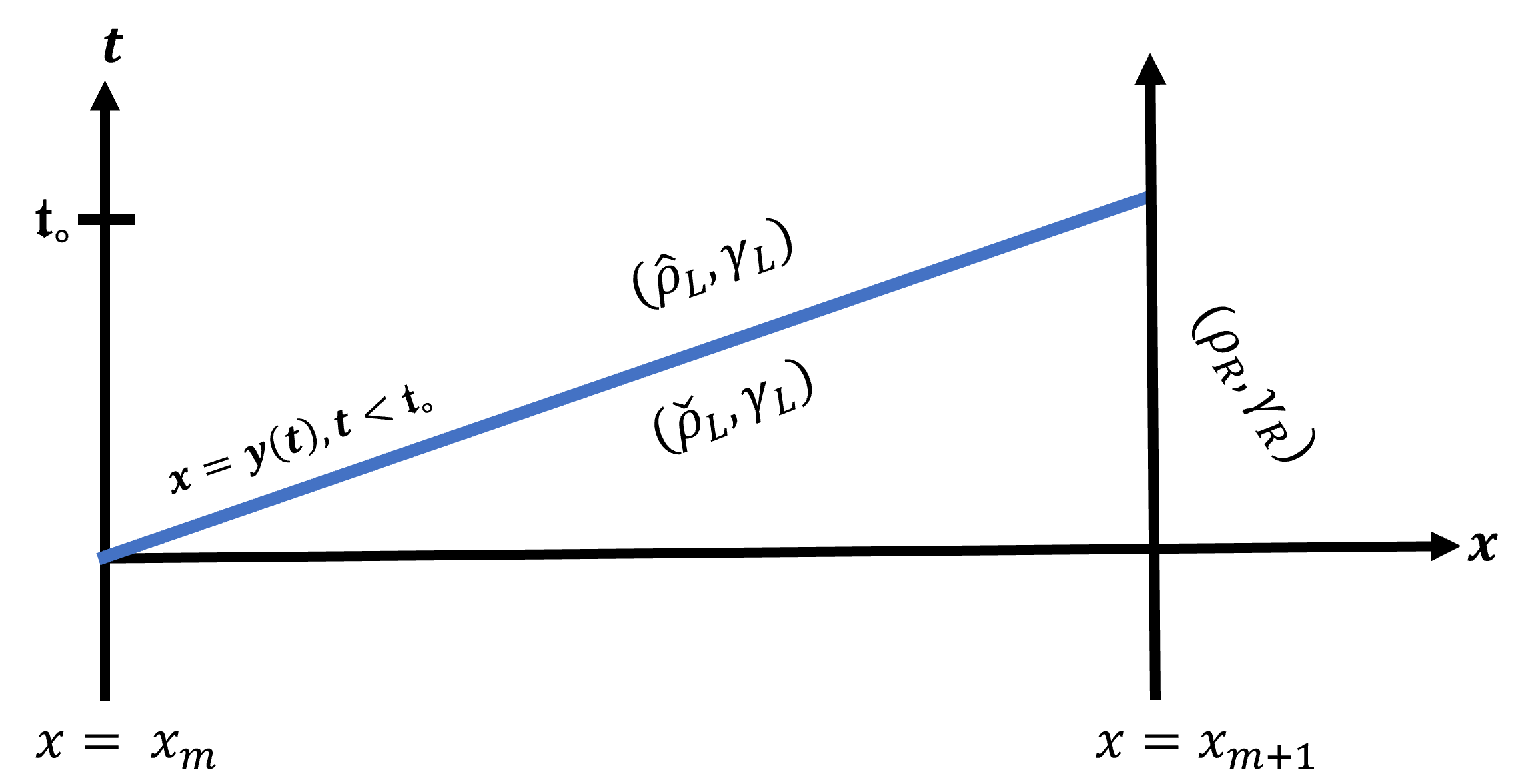}
    \includegraphics[width=3in]{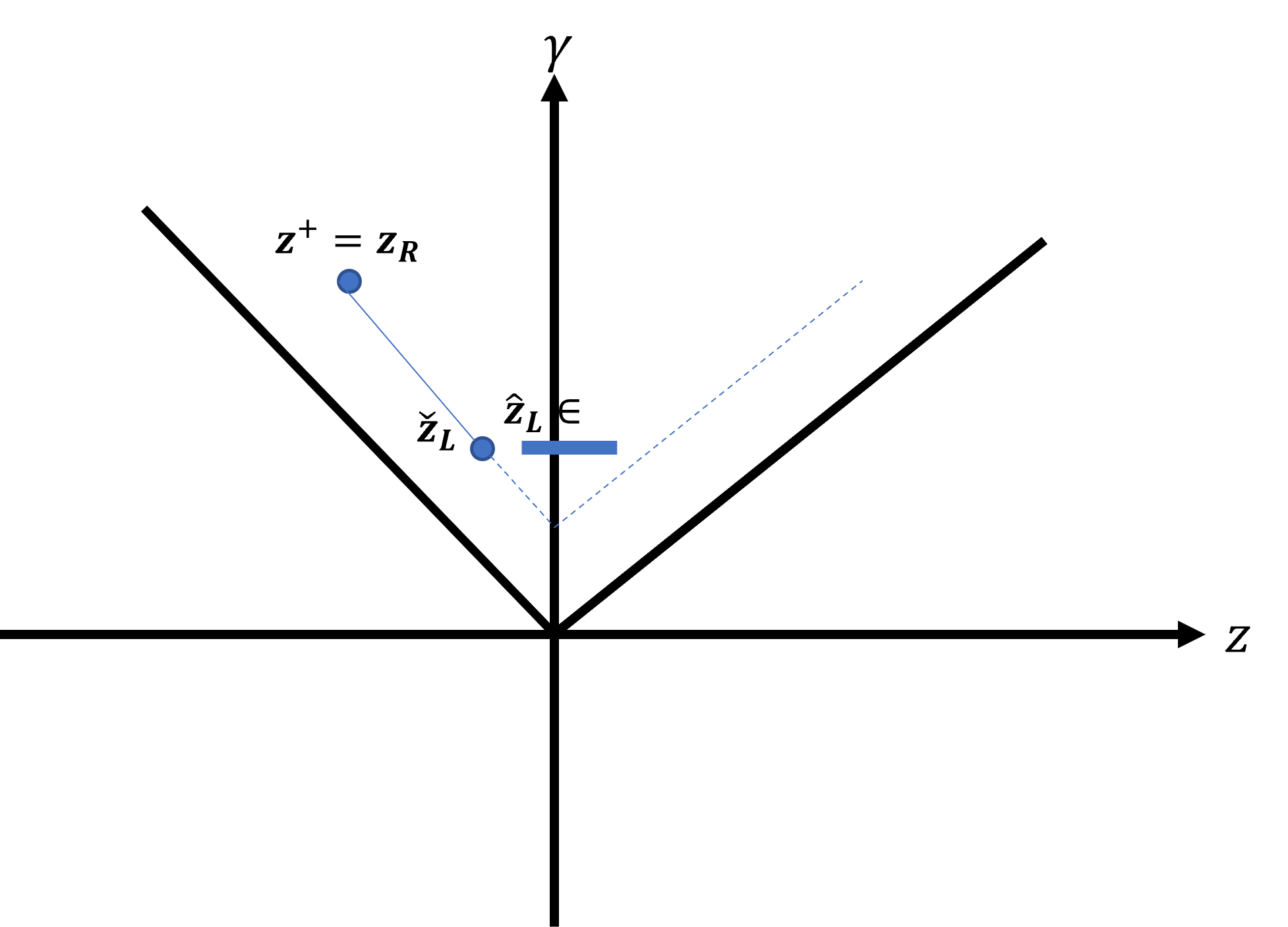}
    \caption{The interaction of non-classical shock and a $\gamma$-front in the left and the acceptable range of $\hat z_L$ in the right.}
    \label{fig:interaction_97_98}
\end{figure}

Let's choose $\hat z_L >0$ which creates some new possibilities which will be of interest in this paper. The Riemann solution $\Riemann(\hat \rho_L, \rho_R; \gamma_L, \gamma_R)$ is shown in Figure \ref{fig:interaction_99_100}. For the Temple functional remains unchanged in this case, i.e. $\temple(z^{(n)}(\ft_\circ^+, \cdot)) - \temple(z^{(n)}(\ft_\circ^-, \cdot))= 0$,
\begin{figure}
    \centering
    \includegraphics[width=3in]{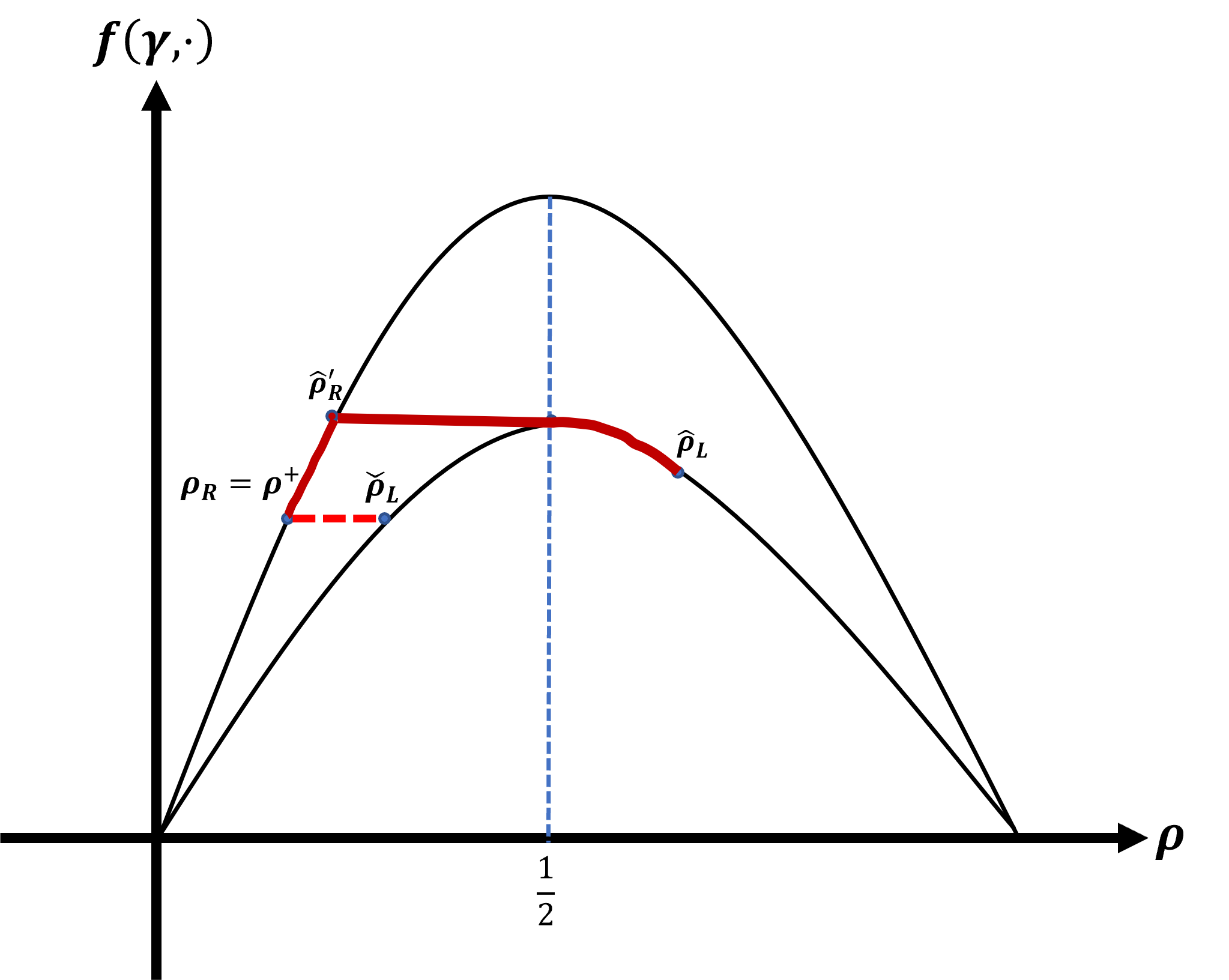}
    \includegraphics[width=3in]{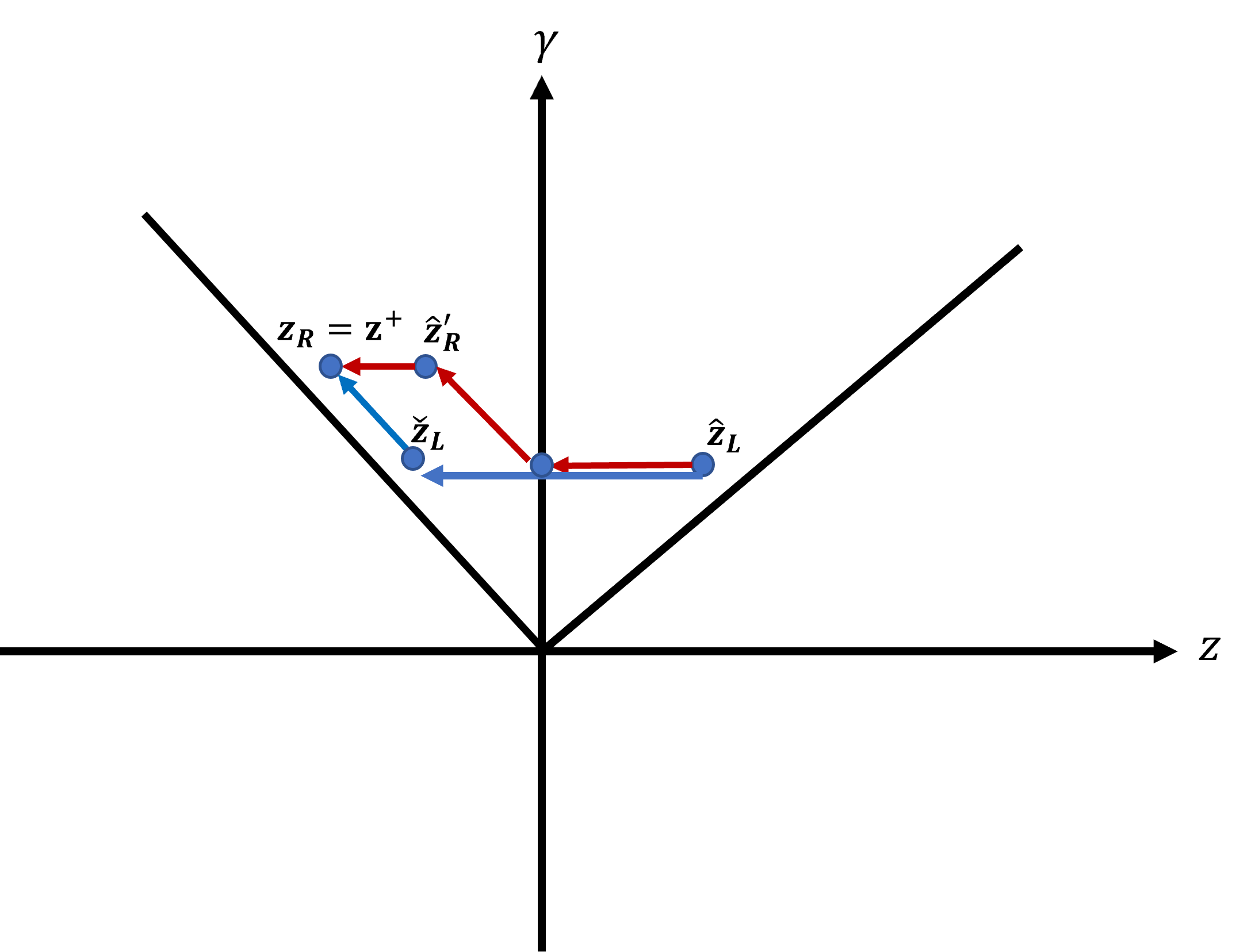}
    \caption{The Riemann solution $\Riemann(\hat \rho, \rho_R;\gamma_L, \gamma_R)$ for $\hat z_L >0$. In the right illustration, }
    \label{fig:interaction_99_100}
\end{figure}

\begin{remark}
It should be noted that the Riemann solution may create rarefactions both from $\hat \rho_L$ to $\rho = \tfrac 12$ and from $\hat \rho'_R$ to $\rho_R$. It is important to note that such a rarefaction solution can only be created when a non-classical shock hits a $\gamma$ front. This will be notable when we discuss the extension of the solution for all time $t \ge 0$.  
\end{remark}
One set of admissible ranges of $\check \rho_R$ and $\hat \rho_R$ is illustrated in Figure \ref{fig:interaction_101}. As illustrated in this figure, for some range of $\hat z_R$ the Riemann solution will be of $z\gamma z$-type and in other range of $z\gamma$-type. 
\begin{figure}
    \centering
    \includegraphics[width=3.2in]{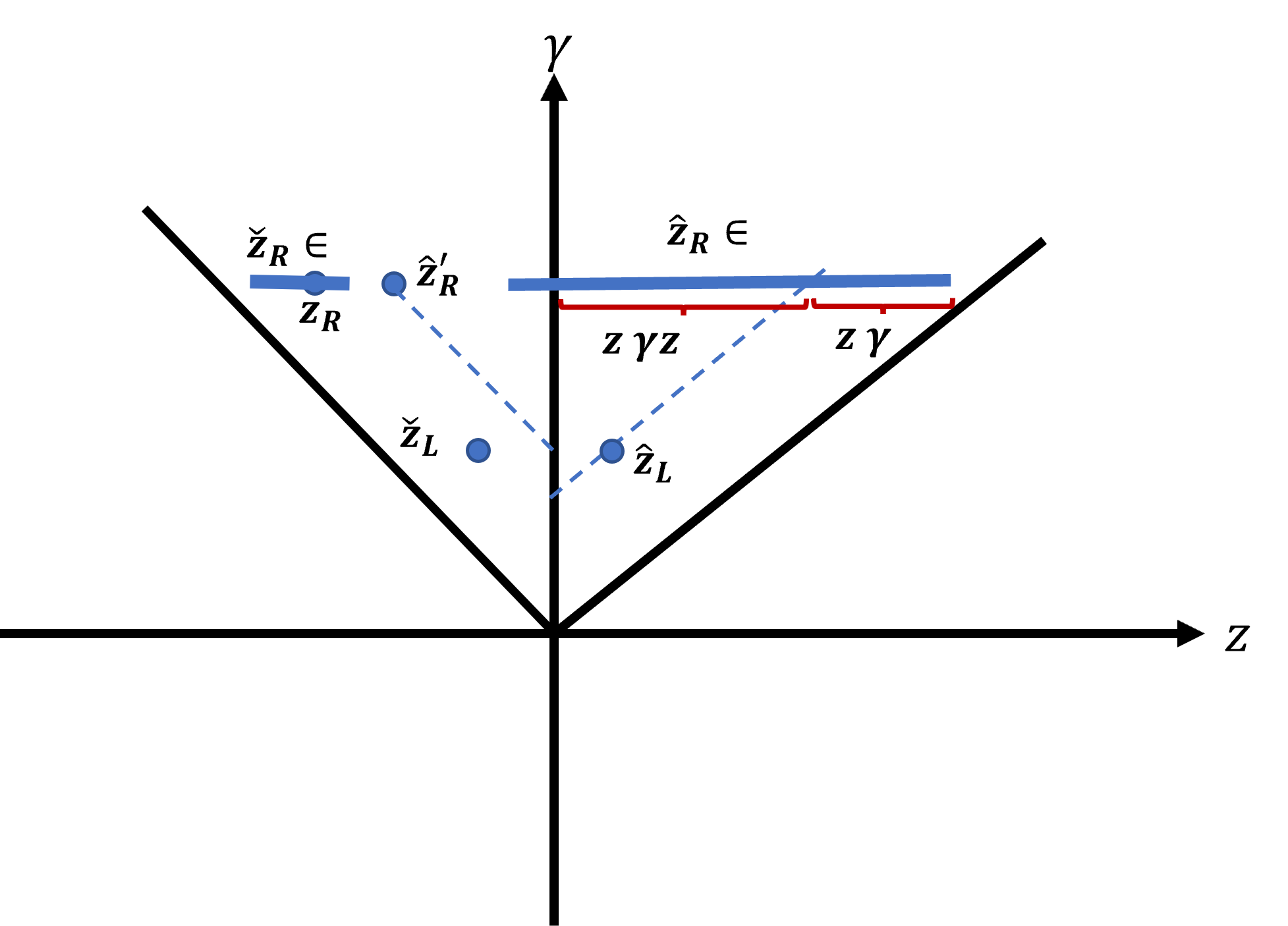}
    
    \caption{ Admissible range for $\check z_L$ and $\hat z_L$.}
    \label{fig:interaction_101}
\end{figure}
Figure \ref{fig:interaction_102_103} illustrates the Riemann solution of $z \gamma z$-type. In addition, Figure \ref{fig:interaction_104} shows the same Riemann solution in the $xt$-coordinates. 
\begin{figure}
    \centering
    \includegraphics[width=3in]{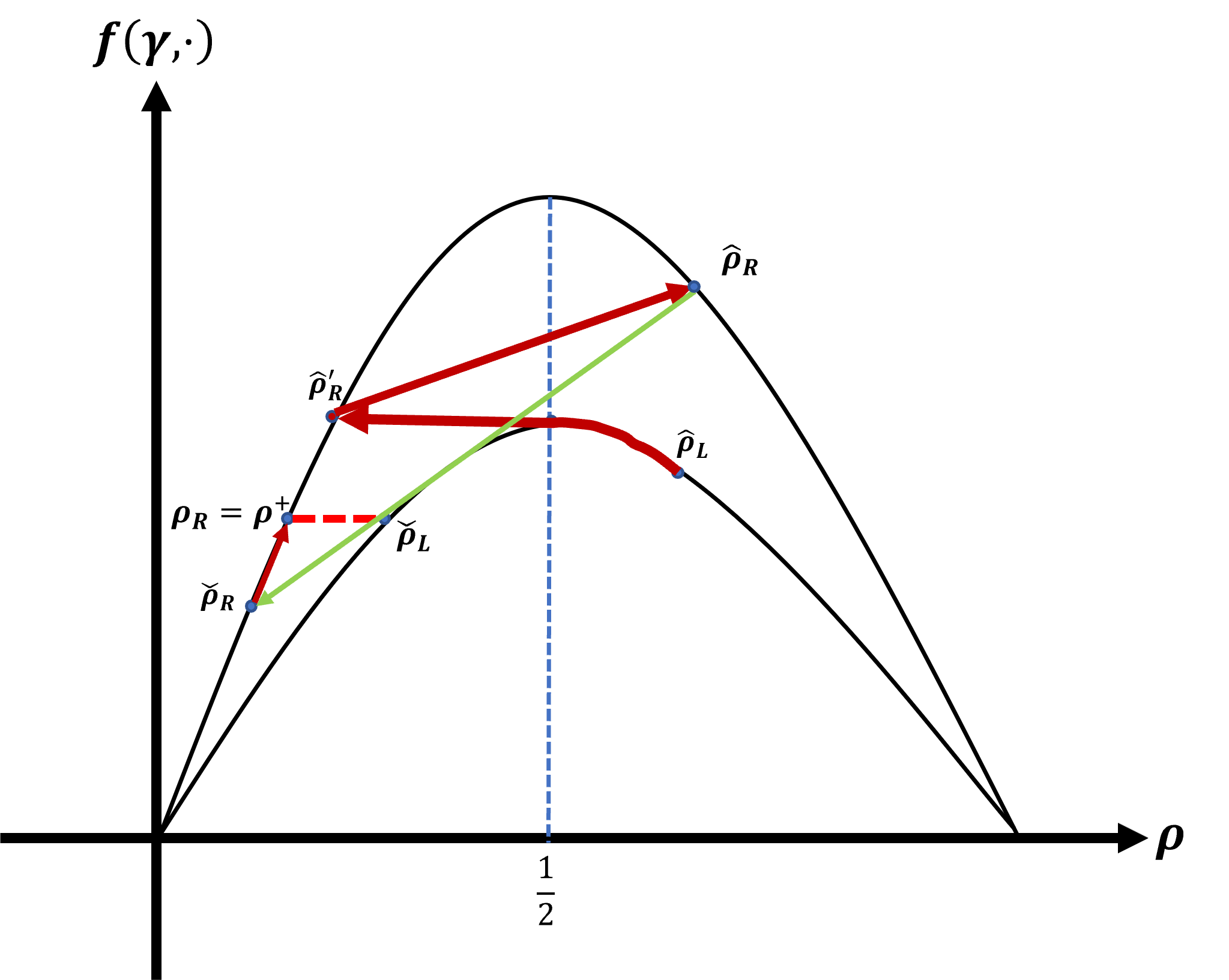}
     \includegraphics[width=3in]{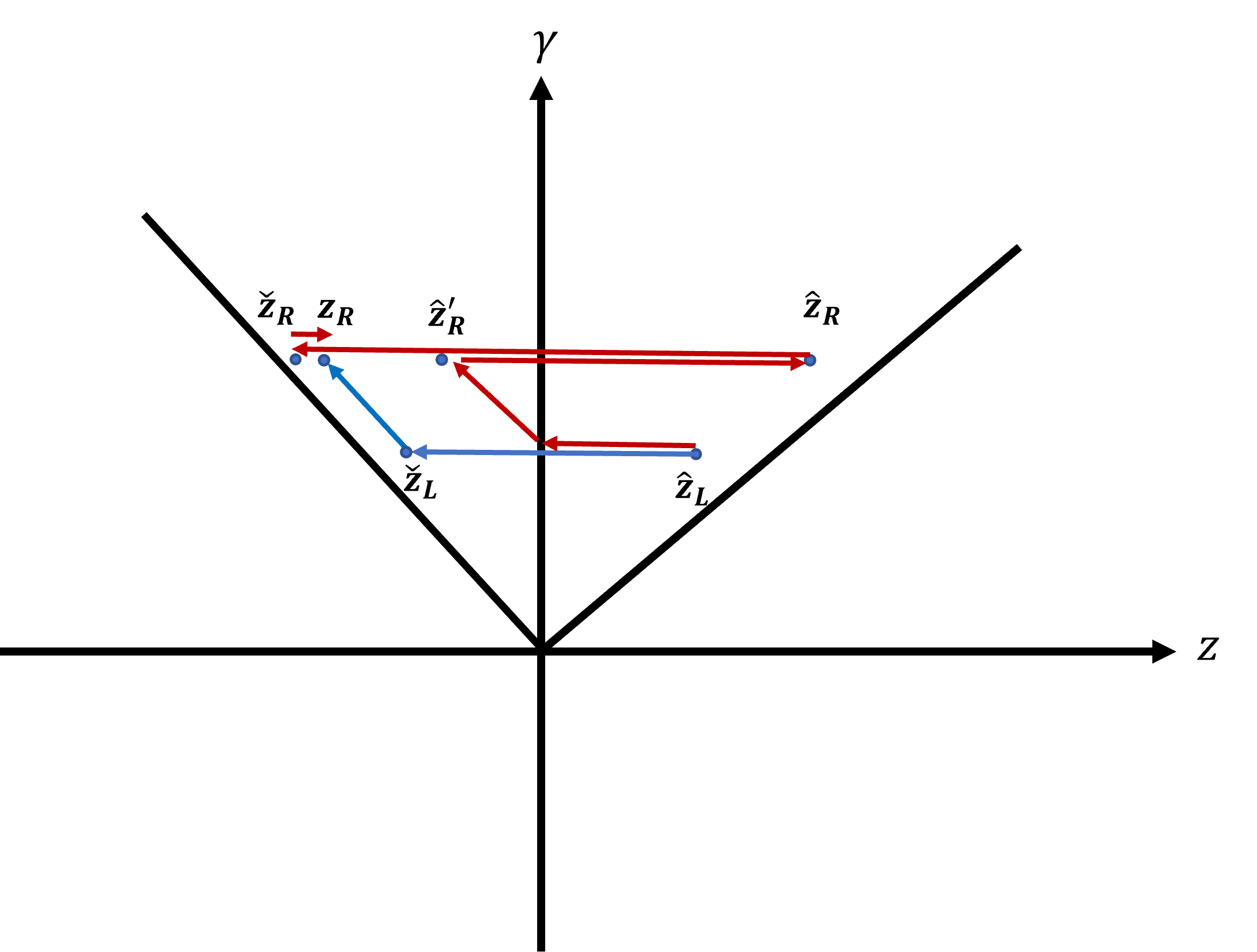}
    \caption{The Riemann solution of $z \gamma z$-type.}
    \label{fig:interaction_102_103}
\end{figure}
\begin{figure}
    \centering
    \includegraphics[width=3in]{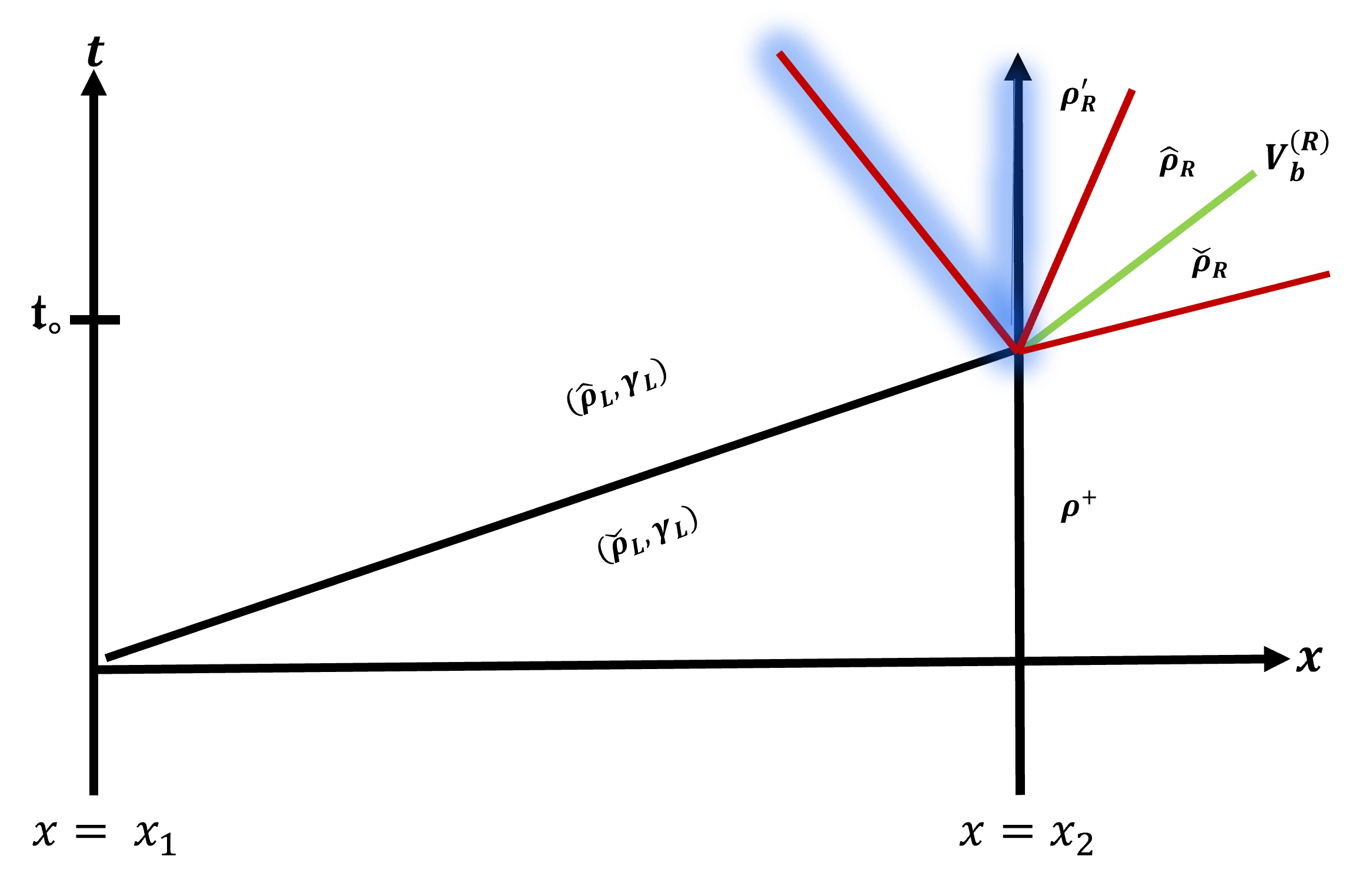}
    \caption{The Riemann solution in $xt$-coordinates. The shadow in the left region indicates the existence of a rarefaction.}
    \label{fig:interaction_104}
\end{figure}
There are a couple of notes: 
\begin{itemize}
    \item The Temple functional is decreasing and $\temple(z^{(n)}(\ft_\circ^+, \cdot)) - \temple(z^{(n)}(\ft_\circ^-, \cdot)) \le - 2 (\hat z_L - \check z_L) \le -  2 \underline \delta^{(n)}$. 
    \item The solution consists of rarefaction.

    \item If $\check \rho_R > \rho_R$ then, the front $\front[\check \rho_R, \rho_R]$ also consists of rarefaction. 
    
\end{itemize}
This in particular shows the dependence of the solution on the values of $\check \rho_R$. In other words, by changing the values of $\hat z_R$, the Riemann solution will be different. Figure \ref{fig:interaction_105_106} shows the Riemann solution of $z\gamma$-type. It should also be noted that similar to the previous case by changing the location, i.e. $\check z_R > \rho_R$, the Riemann solution may consist of a rarefaction. The Temple functional decreases in a similar way as in the previous case. 
\begin{figure}
    \centering
    \includegraphics[width=3in]{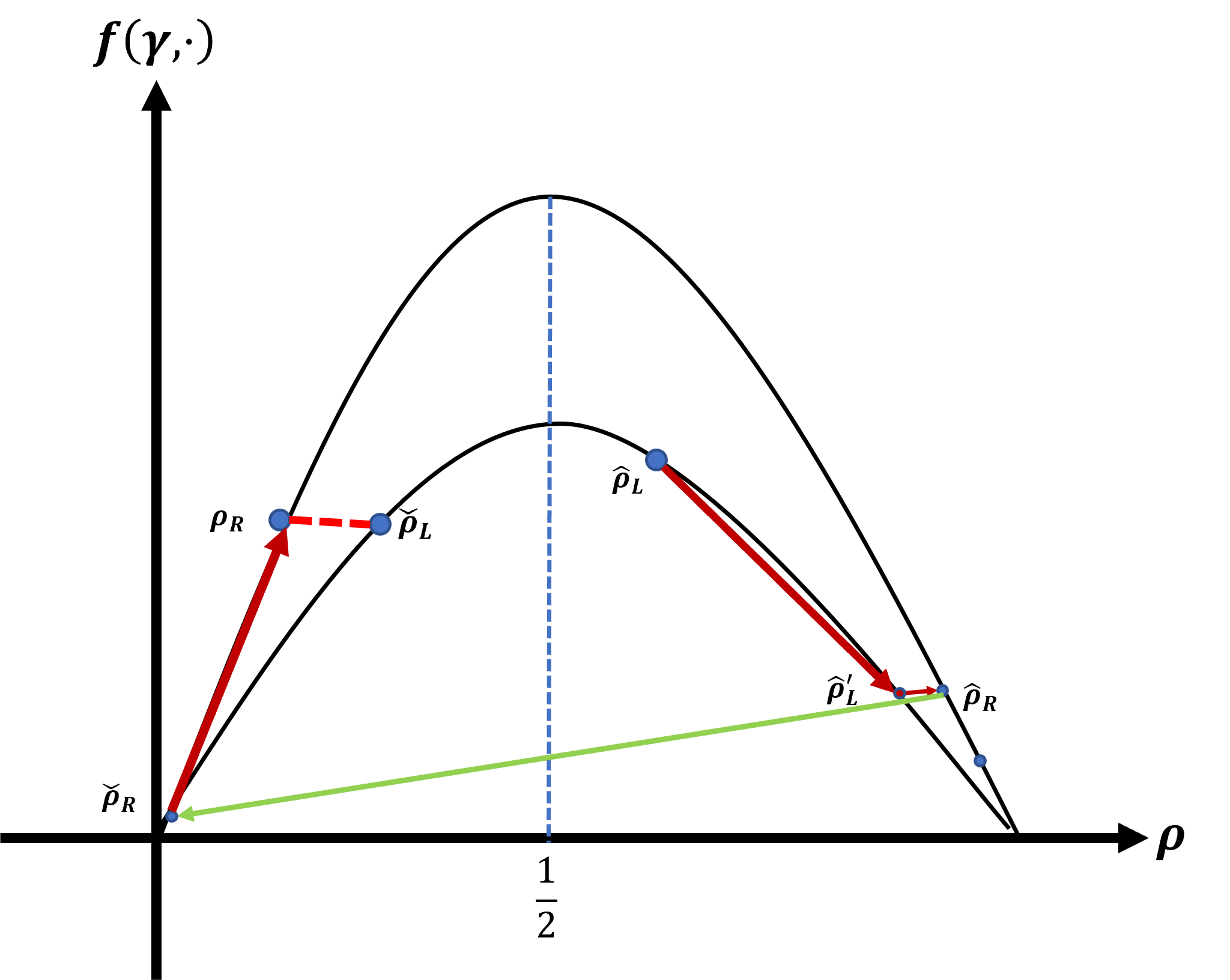}
     \includegraphics[width=3in]{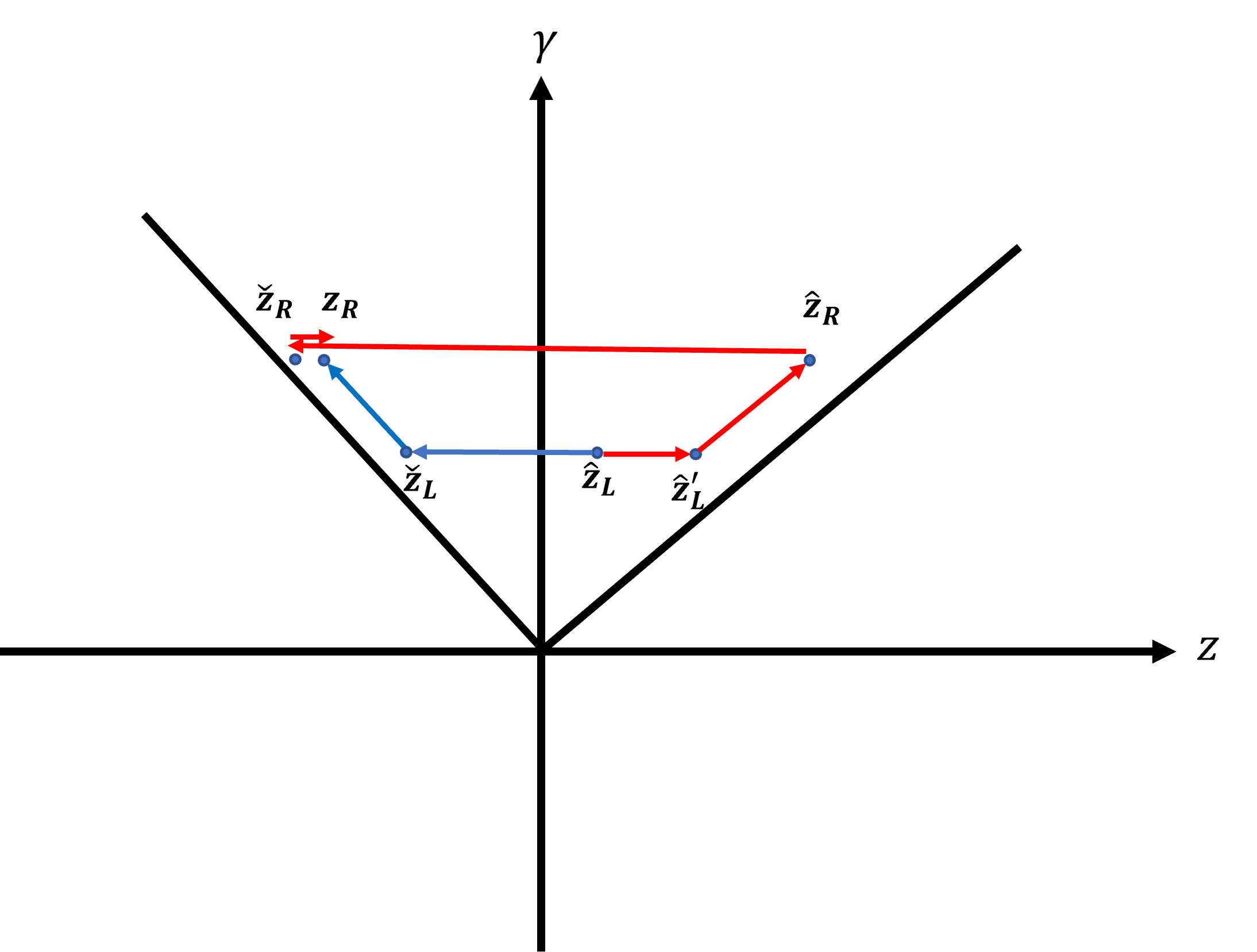}
    \caption{The Riemann solution of $z \gamma$-type.}
    \label{fig:interaction_105_106}
\end{figure}

\end{case}

A thorough investigation of all possible cases in each category concludes that the Temple functional is decreasing and hence the following result is an immediate consequence.
\begin{proposition}\label{p:Temple_decreasing}
Let $(\rho^{(n)}, y_n)$ be the solution of the $n$-approximate Cauchy problem \eqref{E:system} and $(\gamma, z^{(n)})$ the corresponding solution in $\mathcal W$. For any $0 <\check t < \hat t$ in the domain of definition of the solution, we have that 
\begin{equation*}
    \temple(z^{(n)}(\hat t, \cdot)) \le \temple(z^{(n)}(\check t, \cdot))
\end{equation*}
i.e. Temple functional is decreasing. In particular, for any $t \ge 0$ we have that
\begin{equation}
    \totvar{z^{(n)}(t, \cdot)}{\RR} \le \temple(z_\circ^{(n)}(\cdot)).
\end{equation}
\end{proposition}
\begin{lemma}\label{lem:TV_bounded} Let $(\gamma,z^{(n)}) \in \mathcal W$ corresponds to the n-approximate solution $(\rho^{(n)}, y_n)$. For any $t \ge 0$ in the domain of the solution,
\begin{equation*}
    \totvar{z^{(n)}(t, \cdot)}{\RR} \le \totvar{z_\circ(\cdot)}{\RR} + \totvar{\gamma(\cdot)}{\RR }+ C_{\eqref{E:TV_Temple}}.
\end{equation*}
\end{lemma}
 
\begin{proof}
At the time $t = 0$, the claim follows from \eqref{E:bounded_var_initial}. For $t>0$, using \eqref{E:TV_Temple}, we can write
\begin{align*}
    \totvar{z^{(n)}(t, \cdot)}{\RR} &\le  \temple(z^{(n)}(t, \cdot)) \\
    & \le \temple(z_\circ^{(n)}(\cdot)) \\
    &\le \totvar{z_\circ^{(n)}(\cdot)}{\RR} + \totvar{\gamma(\cdot)}{\RR} + C_{\eqref{E:TV_Temple}}\\
    & \le \totvar{z_{\circ}(\cdot)}{\RR}+ \totvar{\gamma(\cdot)}{\RR} + C_{\eqref{E:TV_Temple}},
\end{align*}
where the second inequality is by Proposition \ref{p:Temple_decreasing}, the third inequality is by \eqref{E:TV_Temple} and finally the last inequality is by \eqref{E:bounded_var_initial}. 
\end{proof}
This result implies that the $\totvar{z^{(n)}(t, \cdot))}{\RR}$ is bounded (uniformly) independent of $n$ and $t$. This is one of the necessary conditions for invoking Helly's compactness theorem which provides a tool to show the convergence of approximate solutions to the solution of the Cauchy problem (the next theorem states other necessary conditions).

\begin{theorem} \label{T:L1_Lipschitz_time}
The sequence of approximate solutions $\set{z^{(n)}(t, \cdot): n \in \NN, n \ge N_\circ}$ satisfies the following bounds:
\begin{align}
  \sup_{t>0} \norm{z^{(n)}(t, \cdot)}_{L^\infty(\RR)} \label{E:uniform_bound} & \le \frac 14 \max_{m} \gamma_{r_m} \le C_{\eqref{E:uniform_bound}}\\
    \norm{z^{(n)}(t, \cdot) - z^{(n)}(s, \cdot)}_{L^1(\RR)} & \le \mathbf C_\ell(t -s), \quad \text{for any $0< s< t$} \label{E:modulus_continuity_z}
\end{align}
where the constant $\mathbf C_\ell \Def\max_{m} \gamma_{r_m} \lb \totvar{z_{\circ}(\cdot)}{\RR}+ 2\totvar{\gamma(\cdot)}{\RR} + C_{\eqref{E:TV_Temple}} \rb$ \tup{cf. Lemma \eqref{lem:TV_bounded}} which is independent of $t$ and $n$. 
Moreover, the solution $z^{(n)}(t, x)$ exists for all time $t \ge 0$. 
\end{theorem}
\begin{proof}
Using the fact that $(t, x) \in \RR_+ \times \RR \mapsto \rho^{(n)}(t, x) \in [0, 1]$ of the \eqref{E:discrete_Riemann_multiple} and by the definition of $\psi$ as in \eqref{E:homeo_map}, we have that
\begin{equation*}\begin{split}
   \sup_{x \in \RR} \abs{z^{(n)}(t, x)} & = \tfrac 14\sup_{x \in \RR } \abs{\gamma(x) (2\rho^{(n)}(t, x) - 1)^2}\\
   & \le \tfrac 14 \sup_{x} \abs{\gamma(x)} (2 \rho^{(n)}(t, x) - 1)^2 \\
   & \le \tfrac 14 \max_{m} \gamma_{r_m}. 
\end{split}\end{equation*}
This proves the claim \eqref{E:uniform_bound} (the existence of solution for all $t >0$ will be discussed below in detail).

To prove \eqref{E:modulus_continuity_z}, the broad idea is to show first that the claimed bound is valid before any collision. Then keeping the wave interactions in mind, we show that the claimed bound \eqref{E:modulus_continuity_z} can be extended to the collision point. We start by setting $\tau$ to be the first time that any collision between the fronts happens. Using the initial data \eqref{E:piecewise_approx} and letting $y_\circ = 0$, on $[0, \tau)$ the solution $(\rho^{(n)}, y_n)$ of the $n$-approximate problem is defined by piecing the solutions of following Riemann problems $(\mathbf P_1)$ and $(\mathbf P_2)$ together.
\begin{equation*}
    (\mathbf P_1): \begin{cases} \partial_t \rho + \partial_x[f^{(n)}(\gamma(x), \rho(t, x))]  = 0 , \quad x \in \RR, \, t \in (0, \tau)\\
    \rho_\circ(x) = \begin{cases} \rho_{m_\circ - 1, l} & , x < y_\circ\\
    \rho_{m_\circ, r} & , x> y_\circ \end{cases}, \quad \gamma(x) = \begin{cases} \gamma_{r_{m_\circ} -1} &, x < y_\circ\\  \gamma_{r_{m_\circ}} &, x > y_\circ \end{cases} \\
    f(\gamma_{r_{m_\circ}}, \rho(t, y_n(t))) - \dot y_n(t) \rho(t, y_n(t)) \le F_\alpha^{(m_\circ)}(\dot y(t)), \quad \exists m_\circ \in \set{0, \cdots, M}\\
    \end{cases} \\
\end{equation*}
$\rho_{m_\circ - 1, l}, \rho_{m_\circ , r} \in \mathcal G^{(n)}$, and $y_n$ is the solution of 
\begin{equation}\label{E:BN_traj_region}\begin{split}
    \dot y(t) &= w(y(t), \rho^{(n)}(t, y(t)+)) \\
    y_\circ &= 0 \in I_{m_\circ}\end{split}
\end{equation}
It should be noted that since $y_\circ \in I_{m_\circ}$, the $y_n(t) \in I_{m_\circ}$ for $t \in [0, \tau)$. 
The other Riemann problems are presented by 
\begin{equation*}
    (\mathbf P_2): \begin{cases}  \partial_t \rho + \partial_x[f^{(n)}(\gamma(x), \rho(t, x))]  = 0 , \quad x \in \RR, \, t \in (0, \tau)\\
     \rho_\circ(x) = \begin{cases} \rho^{(n)}_{m, j} & , x \in [x_{m,j -1}, x_{m, j}) \\
    \rho^{(n)}_{m', j'} & , x \in [x_{m', j'-1}, x_{m', j'}) \end{cases}, \quad 
    \gamma(x) = \begin{cases} \gamma_{r_\circ} &, x \in I_\circ \\ \vdots  \\ \gamma_{r_M}  &, x \in I_M \end{cases} 
    \end{cases}
\end{equation*}
for any $m \in \set {0, \cdots, M-1}$, $m' \in \set{1, \cdots, M}$, $j \in \set{1,\cdots, N_m} \setminus \set {l}$ and $j' \in \set{1,\cdots, N_{m'}} \setminus \set{r}$ 
\begin{equation*}
    (m', j') = \begin{cases}
     (m, j +1) &, x_{m, j} < \sa_{m + 1} \\
     (m + 1, 1) &, x_{m, j} = \sa_{m + 1}
    \end{cases}
\end{equation*}
In particular, for $x_{m, j} < \sa_{m+1}$, there would be no $\gamma$-jump and the problem $(\mathbf P_2)$ is essentially the Riemann problem in one region while for $x_{m, j} = \sa_{m + 1}$ there is a $\gamma$-jump from the region $I_m$ to $I_{m +1}$ and the Riemann problems in $(\mathbf P_2)$ need to be solved accordingly.

By construction of the Riemann problem in the Definition \ref{def:Riemann_sol_new}, the solution $(\rho^{(n)}, y_n) \in L^1_{\loc}([0, \tau) \times \RR) \times W^{1,1}_{\loc}([0, \tau); \RR)$, calculated by piecing the solution of the problems $(\mathbf P_1)$ and $(\mathbf P_2)$ together, is a well-defined weak solution of the $n$-approximate problem of the Cauchy problem \eqref{E:system} over the time interval $[0, \tau)$; i.e. before any interaction between the fronts happens. 

Next, we show that \eqref{E:modulus_continuity_z} can be extended to the collision point at time $\tau$ (and consequently beyond $\tau$). To prove this, we need some primary results. 
\begin{lemma}\label{T:phase_1}
    Let's fix $0 < s< t < \tau$. Then, for a.e. $x \in \RR$
 \begin{equation}\label{E:weak_eta_eps_simp}
    \left(\rho^{(n)}(t, x) - \rho^{(n)}(s, x) \right) + \int_{r \in [s, t]}  \frac{\partial}{\partial x}[f^{(n)}(\gamma(x) , \rho^{(n)}(r, x))]  dr = 0.
\end{equation}
In the distributional sense. 
\end{lemma}
To keep the coherency of the discussion, the proof of this lemma is postponed to \ref{A:phase_1}.

\noindent To proceed with the rest of the proof, we need to recall some preliminary definitions in function spaces. 
\begin{remark}\label{R:BV_alternative}
    Let $\mathcal U$ be an open subset of $\RR^n$. We define a (continuous) linear functional $u \in L^1(\mathcal U)\mapsto I_u \in C_c^1(\mathcal U)^*$ by $I_u(\phi) \Def  \int_\mathcal U u\diverg \phi(x) dx$. Then, the $L^1(\mathcal U)$-norm can be defined alternatively considering $u$ as a linear operator on the space of $C_c^1(\mathcal U)$ by
\begin{equation}\label{E:alternative_L1}
    \norm{u}_{L^1(\RR)} = \sup \lb\int \phi(x) u(x) dx : \, \phi \in C_c^1(\RR), \, \norm{\phi} \le 1 \rb
    \end{equation}
    
Next, we recall a definition of a bounded variation function. We define a seminorm 
    \begin{equation}\label{E:seminorm}
        \begin{split}
             \nnorm{I_u} &\Def \sup \set{I_u(\phi): \phi \in C_c^1(\mathcal U), \, \norm{\phi}_\infty \le 1}\\
    & = \sup \set{\int_{\mathcal U} u(x) \diverg \phi dx: \phi \in C_c^1(\mathcal U), \, \norm{\phi }_\infty \le 1 }.
        \end{split}
    \end{equation}
    on the dual topological space $C_c^1(\mathcal U)^*$ associated with the strong dual topology on this space \tup{locally convex space generated by such seminorm on bounded sets}. Now defining $\totvar{u}{\mathcal U}\Def  \nnorm{I_u}$, a function $u \in L^1(\mathcal U)$ is of bounded variation, denoted by $u \in \BV(\mathcal U)$, if $\totvar{u}{\mathcal U} < \infty$. 

    In addition, one can define a function $u \in L^1(\mathcal U)$ is called of bounded variation if there exists a finite vector-valued Radon measure $\mu \in \mathcal M(\mathcal U, \RR^n)$ such that 
\begin{equation*}
    \int_{\mathcal U} u(x) \diverg \phi(x) dx = - \int_{\mathcal U} \action{\phi}{d \mu} , \quad \forall \, \phi \in C_c^1(\mathcal U; \RR^n).
\end{equation*}
This means the weak derivative of $u$ is a Radon measure and in fact $\totvar{u(\cdot)}{\mathcal U}$ can be stated as the norm of the weak derivative. 
\end{remark}
Employing \eqref{E:alternative_L1}, \eqref{E:weak_eta_eps_simp} and Remark \ref{R:BV_alternative} on the operator norm of the linear functional $I_u$, for any $s, t \in [0,\tau)$ with $s<t$ we have that
\begin{equation}\label{E:Lipsch_time_rho}\begin{split}
    \norm{\rho^{(n)}(t, \cdot) - \rho^{(n)}(s, \cdot)}_{L^1(\RR)} 
    & \le \int_s^t \nnorm{f^{(n)}(\gamma(\cdot), \rho^{(n)}(r, \cdot)} dr \\
    & = \int_s^t \totvar{f^{(n)}(\gamma(\cdot), \rho^{(n)}(r,\cdot))}{\RR} dr\\
    & \le \lb\totvar{z_{\circ}(\cdot)}{\RR}+ 2\totvar{\gamma(\cdot)}{\RR} + C_{\eqref{E:TV_Temple}} \rb (t - s)
\end{split}\end{equation}
The last inequality is by the next lemma and the Lemma \ref{lem:TV_bounded}.
\begin{lemma}\label{lem:TV_fn}
We have 
\begin{equation*}
     \totvar{f^{(n)}(\gamma(\cdot), \rho^{(n)}(t, \cdot)}{\RR} \le \totvar{{z^{(n)}(t, \cdot})}{\RR} + \tfrac 14\totvar{\gamma(\cdot)}{\RR}
\end{equation*}
\end{lemma}
See \ref{A:TV_fn} for the proof.

\noindent Noting that for any $\gamma$, the function $u \mapsto \psi(\gamma, u)$  
is smooth, for any $0 \le s < t < \tau$ using the mean value theorem
\begin{equation*}
    \begin{split}
        z^{(n)}(t, x) - z^{(n)}(s, x) &= \psi(\gamma(x), \rho^{(n)}(t,x)) - \psi(\gamma(x), \rho^{(n)}(s,x)) \\
        & = \fpartial{\psi}{\rho}(\gamma,\xi) \left(\rho^{(n)}(t, x) - \rho^{(n)}(s, x) \right)
    \end{split}
\end{equation*}
where, $\xi = \rho^{(n)}(\theta t + (1- \theta) s,x)$ for some $\theta \in (0, 1)$. Since $\rho \in [0, 1]$, $\abs{\frac{\partial \psi}{\partial \rho}} \le \max_{m} \gamma_{r_m}$ is bounded, we have that
\begin{equation}
\label{E:time_Lipschitz_finite}
\norm{z^{(n)}(t,\cdot) - z^{(n)}(s, \cdot)}_{L^1(\RR)} \le \mathbf C_\ell \abs{t - s}, \quad \text{for $0<s < t < \tau$}.
\end{equation}
So far, we have the desired result for $t \in [0, \tau)$. To extend this result to and beyond $t= \tau$, we need a couple of more steps. Let us consider a sequence $(t_m)_{m \in \NN} \in [0, \tau)$ such that $t_m \nearrow \tau$ as $m \to \infty$. In addition, let
\begin{equation} \label{E:sequence_before_collision}
\bar z^{(n)}_m (\cdot) \Def z^{(n)}(t_m, \cdot)\end{equation}
Equation \eqref{E:time_Lipschitz_finite} then implies that $\set{\bar z^{(n)}_m: m \in \NN}$ is a Cauchy sequence in the Fr\'echet space (locally convex, metrizable and complete) $L^1_{\loc}(\RR)$ endowed with the topology generated by the countable family of seminorms $\nnorm{u}_{\Omega_k}$ for open and bounded subsets $\Omega_k \subset\subset \Omega_{k + 1}$ and $\cup_{k \ge 1} \Omega_k = \RR$. This implies the convergence to some limit function in $L^1_{\loc}(\RR)$, denoted by $z^{(n)}(\tau, \cdot)$. In particular, 
\begin{equation}\label{E:uniqeness_limit}
    z^{(n)}(\tau, \cdot) \Def \lim_{t \nearrow \tau} z^{(n)}(t, \cdot)  , \quad \text{in $L^1_{\loc}(\RR)$}.
\end{equation}
Furthermore, we have that 
\begin{equation}\label{E:time_Lipschits_tau}
    \norm{z^{(n)}(\tau,\cdot) - z^{(n)}(s, \cdot)}_{L^1(\RR)} \le C_{\eqref{E:time_Lipschitz_finite}} \abs{\tau - s}, \quad \text{for any $0<s < \tau$}. 
\end{equation}
We consider the sequence $\set{\bar z^{(n)}_m:m \in \NN} \subset L^1_{\loc}(\RR)$ again. By using Lemma \ref{lem:TV_bounded} and \eqref{E:uniform_bound} we will have that
\begin{equation*}
    \totvar{\bar z^{(n)}_m(\cdot)}{\RR}\le \totvar{z_\circ(\cdot)}{\RR } +  \totvar{\gamma(\cdot)}{\RR} + C_{\eqref{E:TV_Temple}}, \qquad \sup_{m}\norm{z^{(n)}_m}_{L^\infty(\RR)} < \infty.
\end{equation*}
Therefore, by Helly's selection theorem, this sequence has a $L^1_{\loc}$-convergent subsequence with the limit function $\bar z^{(n)} \in L^1_{\loc}(\RR)$ with bounded variation. In particular, by uniqueness of limit \eqref{E:uniqeness_limit} in $L^1_{\loc}(\RR)$, we should have that 
\begin{equation}
    \label{E:limit_totvar}
    \totvar{z^{(n)}(\tau, \cdot)}{\RR} \le \totvar{z_\circ(\cdot)}{\RR} + \totvar{\gamma(\cdot)}{\RR} + C_{\eqref{E:TV_Temple}}. 
\end{equation}
By construction of the solution, at the time $t = \tau$, the first collision, the number of fronts is finite. Therefore, \eqref{E:time_Lipschits_tau}, \eqref{E:limit_totvar} and the fact that $z^{(n)}(\tau, \cdot)$ is piecewise constant and bounded, suggest that we can choose this function as the initial value of the Riemann problems after the first collision along with $y_n(\tau)$ as the initial value of the dynamics of the bottleneck. 

Before generalizing the claim to any time $t \ge 0$, we need to study the growth rate of the number of fronts after the collisions to ensure that the solution does not grow unboundedly in finite time. We consider $I_{r_\circ}$ to be the starting region for the bottleneck. The time $t = \tau_m$ denotes the first time that the bottleneck hits the boundary $\partial I_{r_m}$. We need to consider the growth rate of both classical and non-classical fronts. In particular, the interaction of $z$ fronts and the bottleneck trajectory can create a non-classical shock in some cases; see Figure \ref{fig:interaction_108}.
\begin{figure}
    \centering
    \includegraphics[width=3.5in]{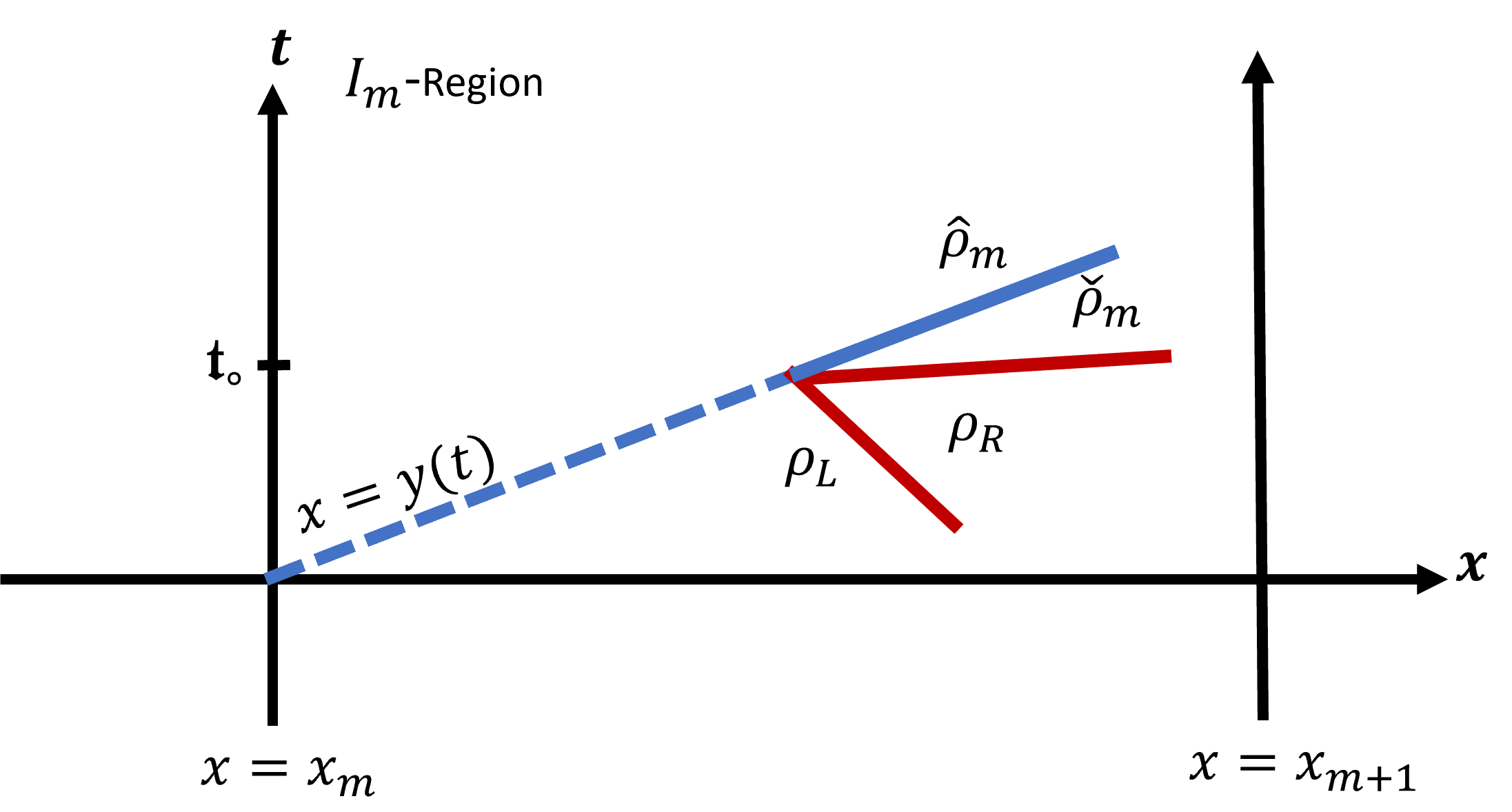}
    \caption{Here $\rho_L= \hat \rho_m$ and $\rho_R < \rho_L$. A new non-classical shock is created.}
    \label{fig:interaction_108}
\end{figure}
Therefore, the number of non-classical shocks can increase as a result of frequent collisions with $z$ fronts (see Figure \ref{fig:interaction_109}). On the other hand, it can be noted that
\begin{enumerate}
    \item Non-classical shocks do not contribute to increasing the number of classical $z$ fronts. In other words, the interaction of a $z$ front and the bottleneck trajectory does not create any new $z$ front in the interior of any region, 

    \item A non-classical shock can create new rarefaction at the collision with the $\gamma$ fronts; recall Case \ref{case:2}.

    \item The non-classical fronts are merely created as the result of collision of the bottleneck trajectory with $z$ fronts (in particular rarefaction) or collision of bottleneck trajectory with the $\gamma$ fronts. This means the number of generated non-classical fronts is bounded by the number of $z$ fronts and $\gamma$ fronts at any time or equivalently, the number of non-classical fronts cannot grow unbounded without $z$ fronts growing unbounded.
\end{enumerate}
\begin{figure}
    \centering
    \includegraphics[width=3.5in]{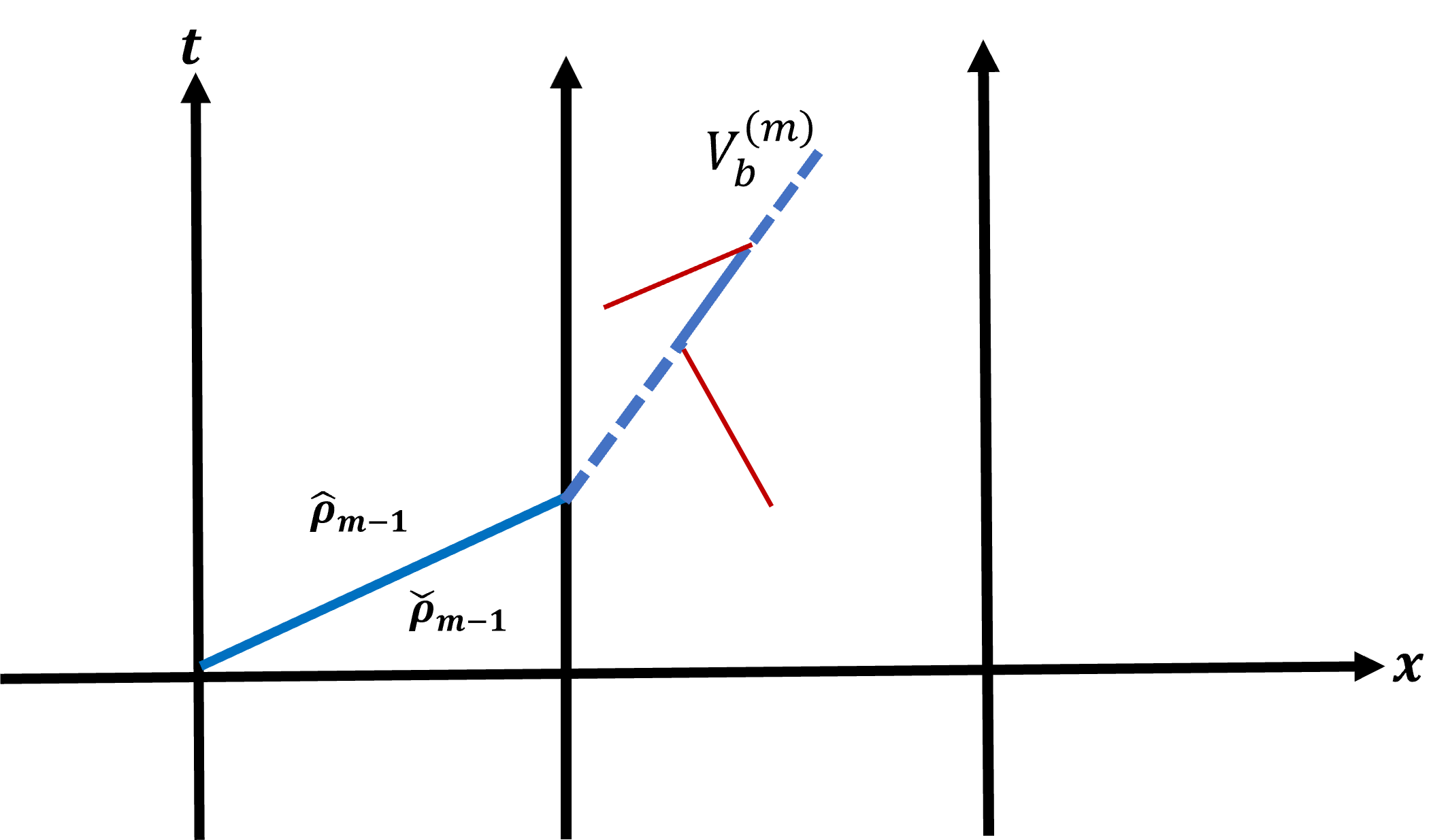}
    \caption{The non-classical shock is created and removed and does not contribute to increasing the number of fronts after collision.}
    \label{fig:interaction_109}
\end{figure}

From these points, one concludes that it is sufficient to understand the growth rate of the classical fronts. Let's suppose that the number of fronts, $N_\front(t) \to \infty$ as $t \to \hat \tau$ for some time $\hat \tau < \infty$. This implies that there should exist an interval $(\check \tau, \hat \tau)$ over which $N_\front(t)$ is strictly increasing. Let $t < \tau_1$. As a result of the collision of waves in each region, the number of fronts decreases. If a $z$ fronts hits a $\gamma$ front, the number of fronts at most remains the same (it may decrease if two or more waves hit the $\gamma$ front). Therefore, $N_\front(t)$ is decreasing on $(0, \tau_1)$. 

At $t = \tau_1$ when the bottleneck collides with a boundary, $N_\front(t)$ can increase at most in the order of $\mathcal O(\bar \delta^{(n)})$. However, $N_\front(t)$ will be decreasing on $(\tau_1, \tau_2)$ and this will be repeated until the bottleneck collides with the last boundary. Therefore, there is no region over which the number of fronts can be strictly increasing and hence $N_\front(t)$ cannot grow to infinity in finite time.

Putting all together, we can extend the Riemann solution $(t, x) \mapsto z^{(n)}(t, x)$ to any arbitrary time in this way. In particular, we have proven \eqref{E:modulus_continuity_z} for any $s, t >0$ and this completes the proof of Theorem 
\ref{T:L1_Lipschitz_time}. 
\end{proof}
This leads us to the main convergence result of the solution. We will prove the existence of the solution in several results.
\begin{theorem}\label{T:convergence_rho_y}
Let $(\rho^{(n)}, y_n)$, $n \in \NN$ be the solution to the $n$-approximate Cauchy problem constructed by the wave-front tracking scheme. Assume further that $\totvar{\rho_\circ}{\RR} < \infty$ be bounded and $\rho_\circ \in [0, 1]$. Then up to a subsequence, we have that
\begin{align}\label{E:rho_convergence}
   & \rho^{(n)} \to \rho , \quad \text{in $L^1_{\loc}(\RR_+ \times \RR; [0, 1])$}; \\
   \label{E:y_n_convergence}
    & y_n \to y , \quad \text{in $C([0, T]; \RR)$ for any $T >0$}
\end{align}
where $\rho \in C(\RR_+; L^1_{\loc}(\RR))$.

Furthermore, $\rho$ is a weak solution in the sense that
  \begin{equation}\label{E:weak_sol_rho}
        \int_{\RR_+} \int_{\RR } (\rho \partial_t \varphi + f(\gamma, \rho) \partial_x \varphi) dx dt + \int_{\RR}\rho_\circ(x) \varphi(0, x) dx = 0, \quad \forall \varphi \in C_c^\infty(\RR_+ \times \RR);
    \end{equation}
    and the inequality \eqref{E:entropy} holds.
\end{theorem}
\begin{proof}
By Lemma \ref{lem:TV_bounded}, Theorem \ref{T:L1_Lipschitz_time} and employing Helly's compactness theorem, there is a subsequence $\set{z^{(n_k)}: k \in \NN}$ and a function $z \in L^1_{\loc}(\RR_+ \times \RR) \cap L^\infty(\RR_+; \BV(\RR))$ such that 
\begin{equation}\label{E:convergent_subseq}
    z^{(n_k)} \to z, \quad \text{in $L^1_{\loc}(\RR_+ \times \RR)$},
\end{equation}
and $z$ satisfies the $L^1$-Lipschitz continuity \eqref{E:modulus_continuity_z} on $\RR_+$. 
Using inverse function $\psi^{-1}$, the fact that $\gamma$ is bounded away from zero, uniform boundedness of $z^{(n)}$ (as of \eqref{E:uniform_bound}) and applying dominated convergence theorem, we conclude that there exists a limit function $\rho\Def \psi^{-1}(\gamma, z)$(up to subsequence) of the sequence $\set{\rho^{(n)}:n \in \NN}$ such that $\rho \in L^1_{\loc}(\RR_+ \times \RR)$. In addition, $L^1$- Lispchitz continuity of $\rho$ with respect to time follows from passing the limit in \eqref{E:Lipsch_time_rho}. This completes the proof of \eqref{E:rho_convergence}. 

The sequence $\set{y_n: n \in \NN}$ calculated from the dynamics of \eqref{E:BN_traj_region}. Let $\mathcal F \subset C([0, T]; \RR)$ be the collection of such solutions. 
Noting that $w(y(t), \rho(t, y(t))) \le \max_m V_b^{(m)}$, the collection $\mathcal F$ is uniformly bounded as $\sup_{n\in \NN} \sup_{t \in [0 ,T]} \abs{y_n(t)} < \infty$. In addition, since $\mathcal F$ is equicontinuous, by Arzela-Ascoli theorem, $\mathcal F$ is totally bounded and hence has a Cauchy subsequence which converges in Banach space $C([0, T]; \RR)$ with respect to the uniform topology. This proves the desired result in \eqref{E:y_n_convergence}.


Finally, to see \eqref{E:weak_sol_rho}, first we note that $f^{(n)}(\gamma, \rho^{(n)}) \to f(\gamma, \rho)$ along the convergent subsequence \eqref{E:convergent_subseq}. In addition, $\rho^{(n)}$ and $f^{(n)}$ satisfy \eqref{E:weak_sol_rho}. By $L^1$-convergence of $\rho^{(n)}$, we can pass the limit to show the claim.

Finally, inequality \eqref{E:entropy} holds true for the approximate solution $(\rho_n, y_n)$ \cite{towers2000convergence, holden2015front}. Therefore, passing to the limit using the dominated convergence theorem, we conclude the inequality. 
\end{proof}

So far, we have shown the existence of the limit functions $\rho$ and $y$. Next, we need to show that $(\rho, y)$ solves the dynamic of the bottleneck trajectory.
To do so, we start with the following essential result. 
\begin{theorem}\label{T:doty_n}
    We have that 
    \begin{equation}\label{E:doty_convergence}
         \dot y_n \to \dot y, \quad \text{in $L^1([0, T]; \RR)$ for any $T >0$}.
    \end{equation}
    and $y \in W^{1,1}_{\loc}(\RR_+)$. 
\end{theorem}
\begin{proof}
For the case of continuous flux, using the fact that $\totvar{\dot y_n(\cdot)}{[0, T]} < \infty$ and invoking the compactness theorem the convergence \eqref{E:doty_convergence} follows, see e.g. \cite{delle2014scalar, liard2021entropic}. However, the bounded variation may not be employed in the case of this paper as new waves can travel beyond the region they originated from (in fact, it can be shown that the total variation of the $\dot y_n$ grows with $n$). 
Nevertheless, we can use the monotonically increasing bijection $[0, 1] \ni u\mapsto \psi(\gamma, u) \in [- \nicefrac{1}{4}\gamma,  \nicefrac{1}{4}\gamma]$ to investigate the variation of the bottleneck speed in $\mathcal W$. More precisely, we define
\begin{equation}\label{E:pos_variation_function}
    \xi_n(t) \Def  \psi(\gamma(y_n(t)), \nicefrac{\dot y_n(t)}{\gamma(y_n(t))})
\end{equation}
where map $\psi$ is defined as in \eqref{E:homeo_map}.
\begin{lemma}\label{T:positive_var}
 Let $\posvar{\xi(\cdot)}{\Omega}$ be the positive variation of a function $\xi$ over a set $\Omega$. We have that
\begin{equation}\label{E:total_var_xi}
    \totvar{\xi_n(\cdot)}{[0, T]} 
 \le 2\posvar{\xi_n(\cdot)}{[0, T]} + \norm{\xi_n}_{L^\infty([0, T])} < \infty 
\end{equation}
uniformly independent of $n$ and $t$. 
\end{lemma}
\begin{proof}[Proof of Lemma \ref{T:positive_var}]The first inequality is by directly using the definition of the total variation of $\xi_n(\cdot)$. In addition, $\norm{\xi_n}_{L^\infty([0, T])} < \infty $ since $\dot y_n(t) \le \max_m V_b^{(m)}$. Therefore, it only remains to show that $\posvar{\xi_n(\cdot)}{[0, T]} < \infty$. 

\begin{remark}[Positive Variation in $\xi_n(\cdot)$] 
By the fact that for any fixed $\gamma$, the map $u \mapsto \psi(\gamma, u)$ preserves the order \textup{monotone increasing function}, $\xi_n(t)$ positively varies if and only if $\dot y_n$ does so, for $y_n(t) \in I_m^\circ$. \end{remark}



We might have at most a finite number of positive changes in the bottleneck speed at the collision with the boundary points. Therefore, we focus on the variation of the bottleneck speed in the interior of the regions. First, we recall the only possible case in which $\dot y_n(t)$ has a positive variation in a region $I_m^\circ$ in Figure \ref{fig:positive_variation}; see \cite{delle2014scalar, liard2021entropic} for more detail. 
\begin{figure} 
    \centering
\includegraphics[width=3in]{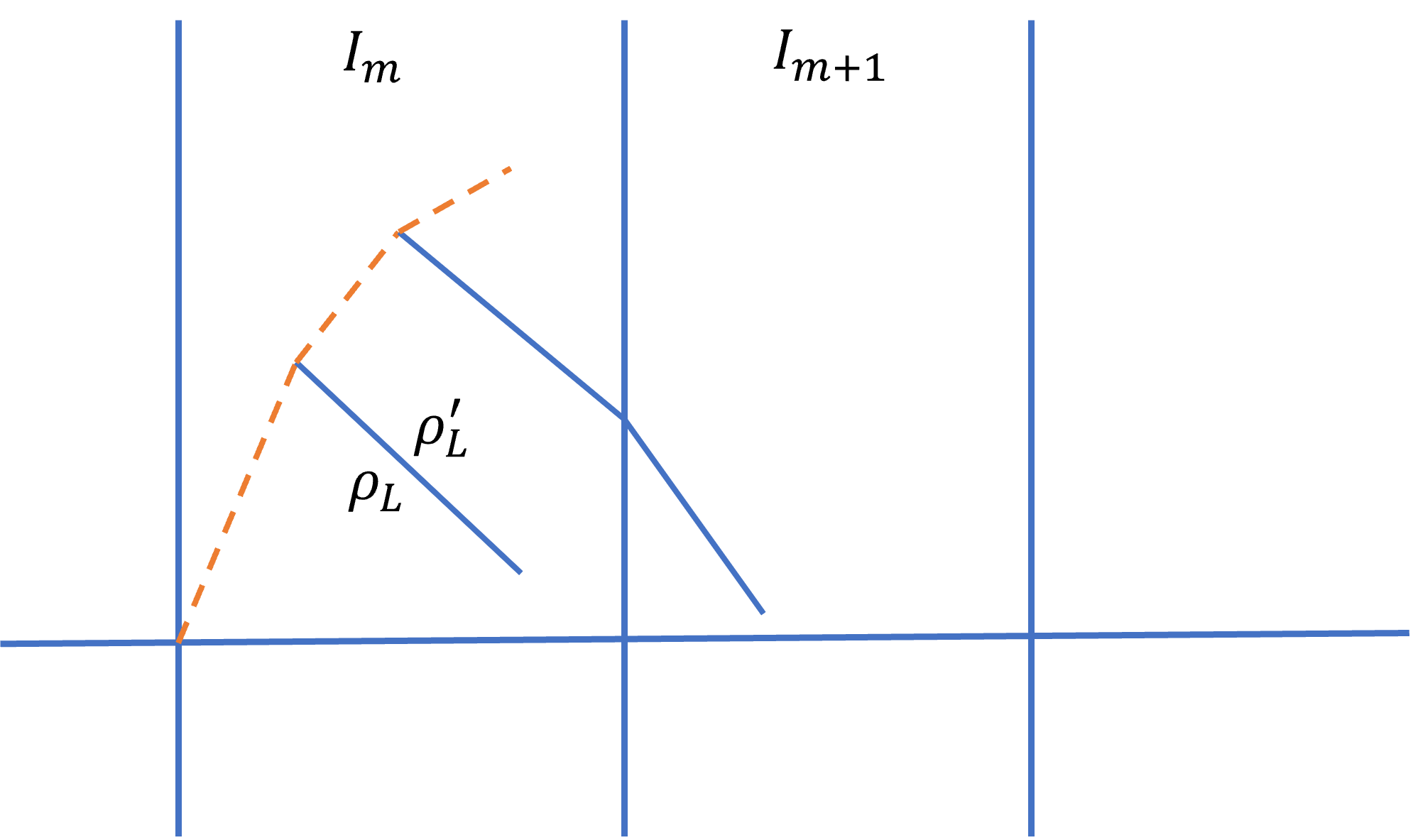}
    \caption{For $\rho_L > \rho_L' \ge \rho^*$ the speed of the bottleneck increases after collision. Such waves can be from inside region $I_m$ or they can be originated in other regions.}
    \label{fig:positive_variation}
\end{figure}
In this case, the collision between a rarefaction shock and the bottleneck trajectory takes place from the right and $\rho_L > \rho_L'\ge \rho^*$. By \eqref{E:velocity}, the positive variation of the bottleneck speed in this case is
\begin{equation}\label{positive_collision}
\dot y_n(\ft_p+) - \dot y_n(\ft_p-) = \gamma_{r_m} (\rho_L- \rho'_L),
\end{equation}
at a time $\ft_p$ of positive variation in $\dot y_n$. In this case, using \eqref{E:pos_variation_function}, we have that
\begin{equation}
    \label{E:positive_collision_trans}
    \xi_n(\ft_p+) - \xi_n(\ft_p-)  = z_L - z'_L.
\end{equation}
The rarefaction waves that can positively change the bottleneck speed, can be originated in several situations which need to be studied carefully.
\paragraph{\textbf{Reflection of Rarefaction Fronts}} Reflections of the waves on the boundary can only be shock waves and hence they cannot positively change the bottleneck speed (see Figure \ref{fig:reflection_shock}). 
\begin{figure}
    \centering
    \includegraphics[width=3in]{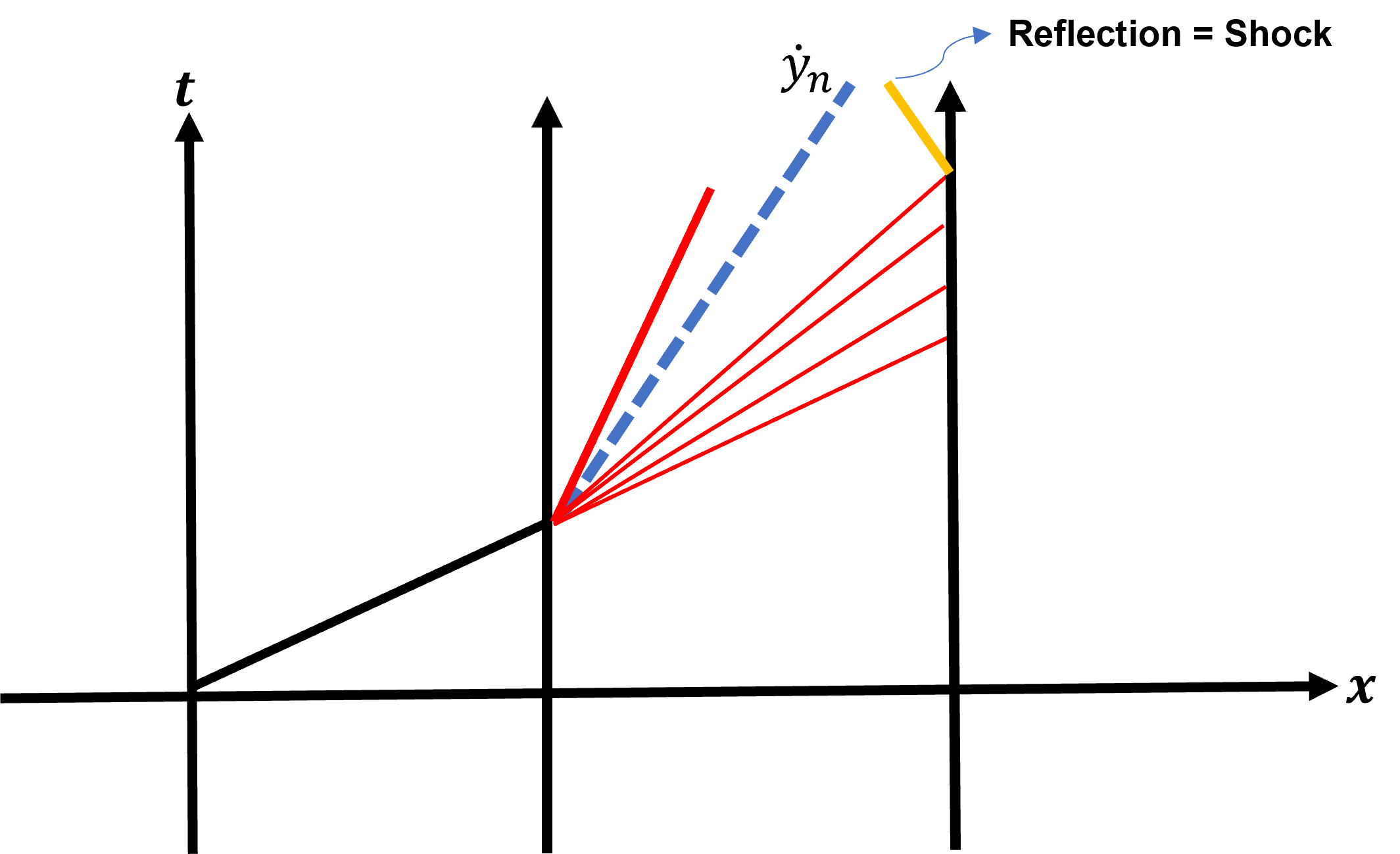}
    \caption{The reflection originated from the collision of a rarefaction with the $\gamma$ front can only be a shock.}
    \label{fig:reflection_shock}
\end{figure}

\paragraph{\textbf{Rarefaction Fronts Initiated at $t=0$}} Considering $y_n(t) \in I_m$, the positive variation in the bottleneck velocity caused by the rarefactions generated in the same region, is at most bounded by
\begin{equation}\label{E:pos_bound_same_region}
     \totvar{z_\circ(\cdot)}{I_m} \le \totvar{z_\circ(\cdot)}{\RR} < \infty, 
\end{equation}
where $\tau_{m+ 1}$ denotes the first time that the bottleneck trajectory approaches the boundary point $x = \sa_{m + 1}$, the boundary of $I_{m+ 1}$ region. Here, without loss of any generality, we assume that $T> \tau_{m + 1}$, where $T$ is as of \eqref{E:doty_convergence}. 
\begin{remark}
    It should be noted that as soon as the bottleneck trajectory is hit by any wave created outside region $I_m$, the rarefaction waves inside $I_m$ will not be able to modify the bottleneck speed positively since this would require interaction with other waves and result in loss of rarefaction waves; see Figure \ref{fig:positive_variation}.
\end{remark}
\paragraph{\textbf{Rarefaction Fronts Generated by Non-classical Shock}} In this case the rarefactions generated by the collision of the non-classical shock with the $\gamma$-wave can remain rarefactions after passing through the $\gamma$-wave. See Case \ref{case:2} for the creation of such rarefactions and see figure \ref{fig:rarefaction_positive_variation}. 
\begin{figure}
    \centering
    \includegraphics[width= 3in]{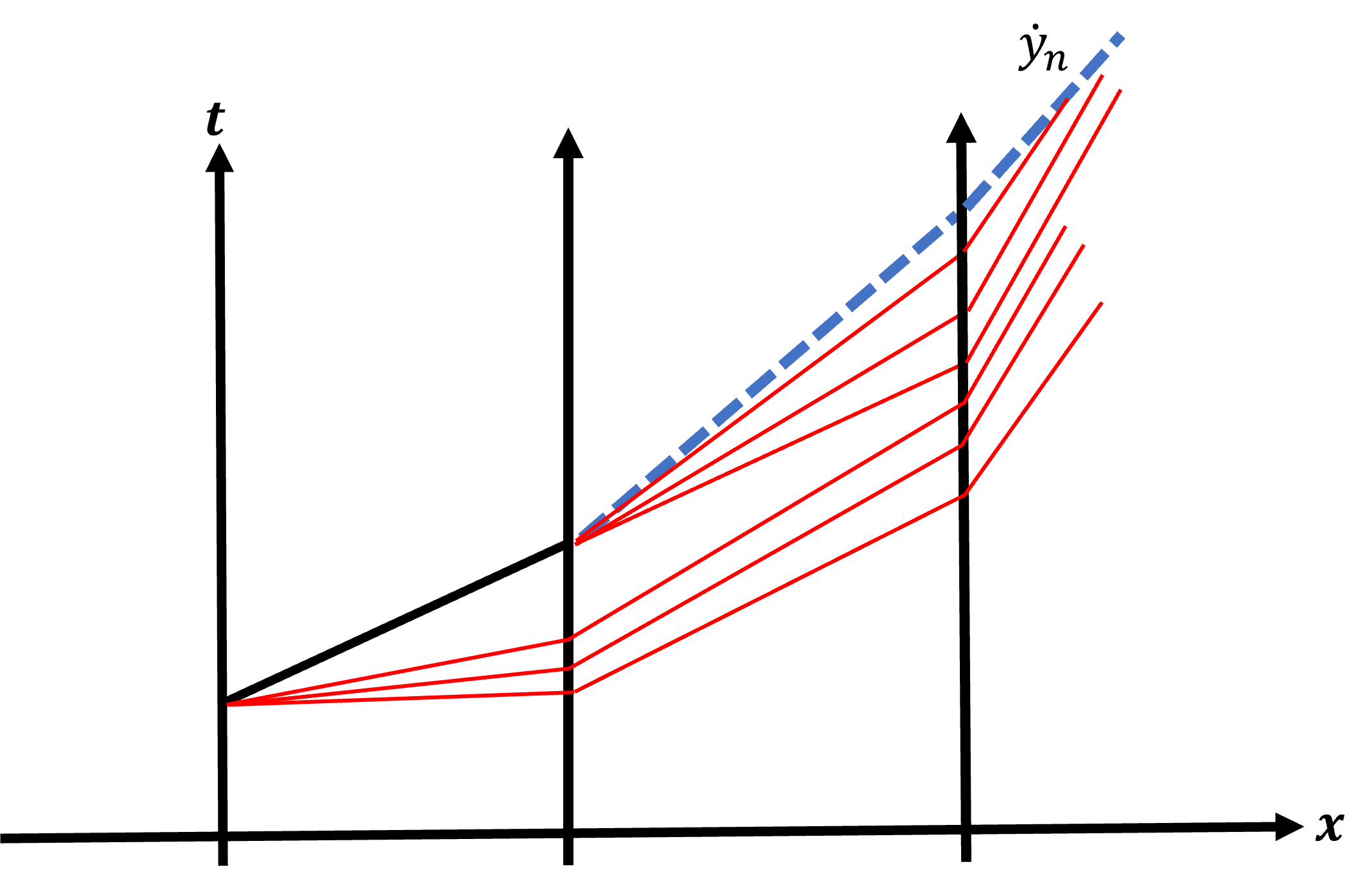}
    \caption{Rarefaction fronts remain rarefaction after collision with $\gamma$ front and can hit the bottleneck trajectory to create positive variation in the speed.}
    \label{fig:rarefaction_positive_variation}
\end{figure}
First, we note that a rarefaction-front can remain a rarefaction after collision with a $\gamma$-front from the left in two cases. One of the instances is shown in Figure \ref{fig:rarefaction_rarefaction} for $\gamma_L > \gamma_R$ and the other case ($\gamma_L < \gamma_R$) is similar. In this case, if the line $x = V_b^{(m)} t$ passes through $\rho^+$, then $\rho'_R > \rho_R \ge \rho^*$ and hence $\dot y_n > \lambda(\rho'_R, \rho_R)$ which implies that the collision is possible. 
\begin{figure}
    \centering
    \includegraphics[width=5in]{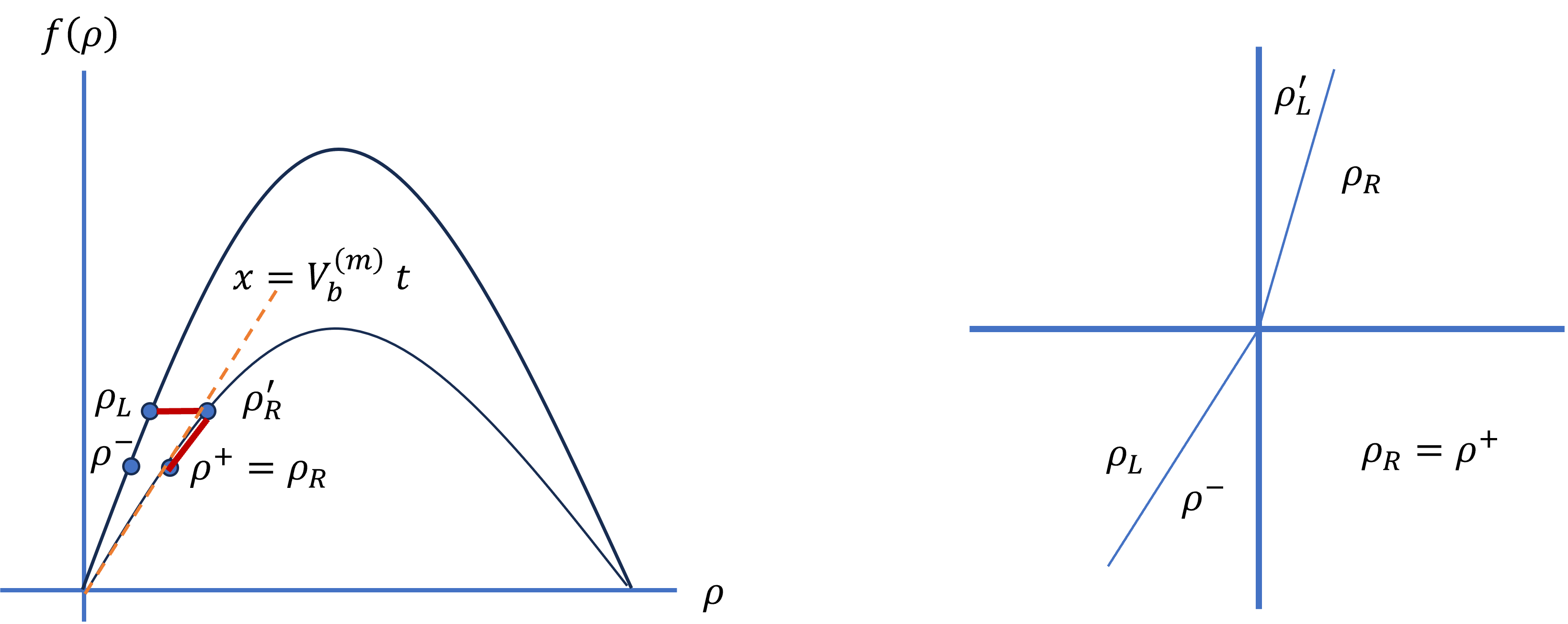}
    \caption{One of the cases in which a rarefaction front remains a rarefaction after collision with a $\gamma$-front.}
    \label{fig:rarefaction_rarefaction}
\end{figure}

Regarding the positive variation of $\xi_n$ in this case and in the presence of several groups of rarefactions (consulting Figure \ref{fig:rarefaction_positive_variation}), it should be noted that the largest value of positive variation corresponds to the case where all of the rarefaction fronts hit the trajectory in order (otherwise the fronts will collide and the number of existing rarefactions decreases). Let $\gamma_{\max} \Def \max\set{\gamma_{r_\circ}, \cdots, \gamma_{r_M}}$. Then, any group of rarefactions contains at most $\gamma_{\max} 2^{n + 1}$ fronts and hence considering \eqref{E:positive_collision_trans} and the Remark \ref{R:UL_delta}, the largest possible positive variation will be $M \gamma_{\max}$. 

\paragraph{\textbf{Rarefaction Fronts Originated in Other Regions}}
We start this part by recalling one instance of the waves generated in other regions that can positively affect the speed $\dot y_n$ (there are only two such cases and due to the similarity of the analysis we only consider one of them here). Figure \ref{fig:positive_variation_otherspaces} shows an arrangement between $\rho_L$, $\rho^+$ and $\rho_R$ (right illustration) such that the collision between the rarefaction front $\front[\rho^+, \rho_R]$ with the $\gamma$-front creates the proper rarefaction wave $\rho_L, 
\rho_L'$ that can positively affect the speed $\dot y_n$; see Figure \ref{fig:positive_negative_variation}.
\begin{figure}
    \centering
    \includegraphics[width=5in]{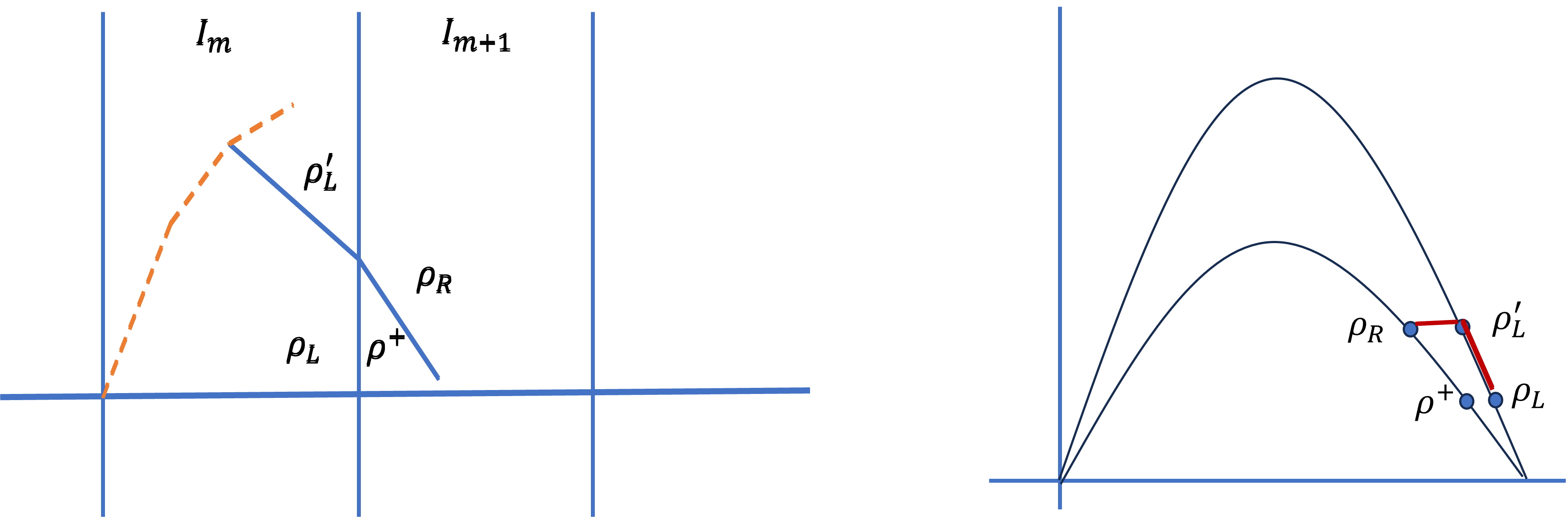}
    \caption{The left illustration shows one of the cases where the collision of a rarefaction $\front[\rho^+, \rho_R]$ with the $\gamma$-front can create a rarefaction with positive speed sign which can affect the speed $\dot y_n$. The illustration on the right shows the relation and position of $\rho_L$, $\rho^+$ and $\rho_R$ that can create the wave in the left figure.}
    \label{fig:positive_variation_otherspaces}
\end{figure}

\begin{figure}
    \centering
    \includegraphics[width=5in]{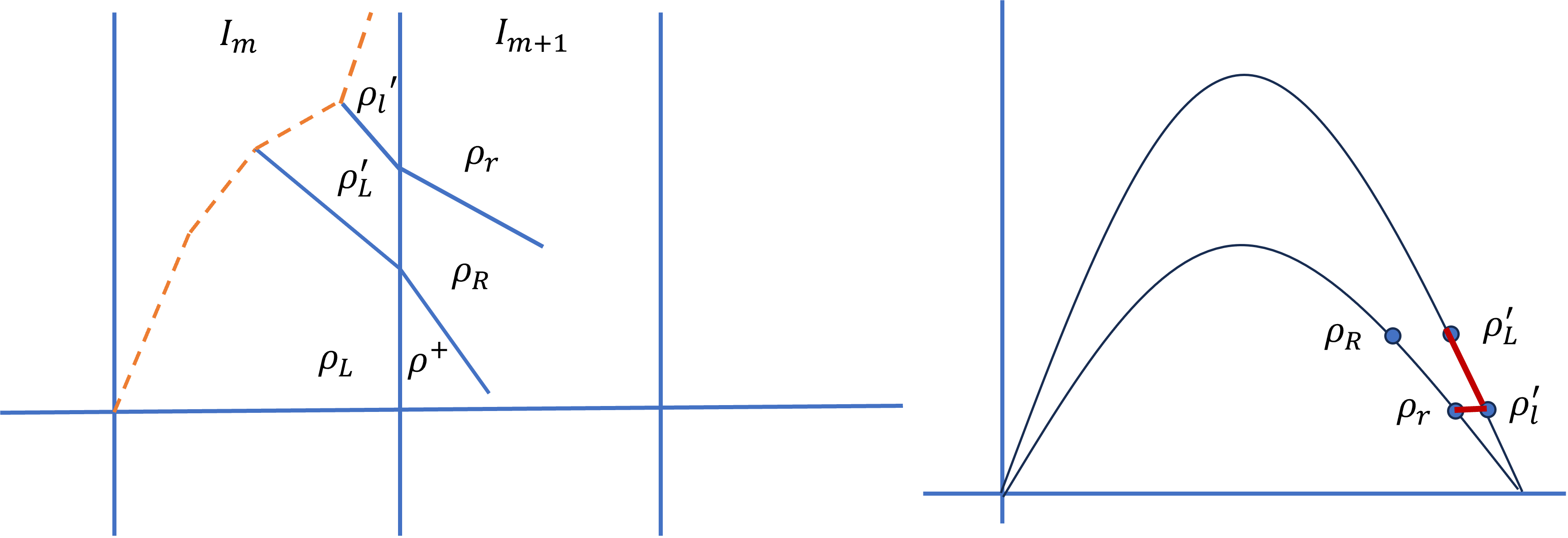}
    \caption{There can be waves hitting the bottleneck trajectory that decrease the speed. Therefore, the process of increasing and decreasing the variation can continue. In this case, $\rho_l' > \rho^*$ and $\rho_L' \in [\hat \rho_m , \rho_l']$. In this case, the speed $\dot y_n$ decreases after the collision; see \cite{liard2021entropic}.}
    \label{fig:positive_negative_variation}
\end{figure}
Therefore, to calculate the positive variation of $\xi_n$ in this case, we need to have an estimate on the number of $z$-fronts that originated in other regions and can potentially hit the bottleneck trajectory $y(t) \in I_m^\circ$ for some $m \in \set{0, \cdots, M}$. Let $\front_z$ denote an individual $z$ front (in particular, we are concerned with rarefactions) and $\Gamma_m(t)$ a collection of $z$ fronts that potentially approach the $\gamma$ front boundaries of region $I_m$ and can change $\dot y_n$ positively at a time $t$.

If a $\front_z$ enters the region $I_m$, either it interacts with the bottleneck, with other waves, or with the other boundary of this region. In any of these cases, after the interaction, the front will not be able to affect the $\dot y_n$ positively. This implies that $\Gamma_m(t)$ will be updated after any collision time $t = \ft_h$ at which a $\front_z$ enters $I_m$. Let's define a functional
\begin{equation*}
    \zeta_m(t) \Def \sum_{z \in \Gamma_m(t)} \abs{\Delta z}.
\end{equation*}
Therefore, 
\begin{equation}
    \label{E:zeta_bound}
    \zeta_m(t) \le \totvar{z^{(n)}(t, \cdot}{\RR} \le \totvar{z_\circ(\cdot)}{\RR} < \infty. 
\end{equation}
We may note that at $t = \ft_h+$, i.e. after $\front_z \in \Gamma_m(t)$ enters the region $I_m$, 
\begin{equation}\label{E:Gamma_decrease}
    \Gamma_m(\ft_h+) \subseteq \Gamma_m(\ft_h-) \setminus \set{\front_z}
\end{equation}
that is, the set is reduced by at least one $z$ front. In addition, since $\abs{\Delta z} \ge \underline \delta^{(n)}$, \eqref{E:Gamma_decrease} asserts that 
\begin{equation}\label{E:zeta_decrease}
    \zeta_m(\ft_h + ) - \zeta_m(\ft_h-) \le - \underline\delta^{(n)}.
\end{equation}
Letting 
\begin{equation}\label{E:passage_time}
    \tau_m \Def \inf \set{t > 0: y(t) \in I_m},
\end{equation}
and using \eqref{E:zeta_decrease}, we can write 
\begin{equation*}
    0 \le \zeta_m(\tau_{m + 1}) \le \zeta_m(\tau_m) - \mathfrak c \underline \delta^{(n)}
\end{equation*}
where $\mathfrak c = \abs{\Gamma_m(\tau_m)}$, the Cardinality of this set. Therefore, by \eqref{E:zeta_bound}
\begin{equation}\label{E:number_front_bound}
    \mathfrak c \le \frac{\zeta_m(\tau_m)}{\underline \delta^{(n)}} \le \frac{\totvar{z^{(n)}(\tau_m, \cdot)}{\RR}}{\underline \delta^{(n)}} \le \frac{\totvar{z_\circ(\cdot)}{\RR}}{\underline \delta^{(n)}}
\end{equation}

This result holds true for any $m \in \set{0, \cdots, M}$ and hence considering \eqref{E:positive_collision_trans}, Remark \ref{R:UL_delta}, and \eqref{E:number_front_bound}, the positive variation of $\xi_n$ in each region that can happen as a result of interaction of the bottleneck trajectory with the rarefactions originating in other regions, can be bounded by
\begin{equation}\label{E:pos_var_bound}
    \mathfrak c \bar \delta^{(n)} \le 4 \totvar{z_\circ(\cdot)}{\RR} < \infty.
\end{equation}
\begin{remark} \label{R:pos_variation_boundaries}
    It should be noted that the positive variation of the bottleneck at the boundary points, if happens, would be only finite by the structure of the problem. 
\end{remark}
Collecting all of these cases together, we conclude that $\posvar{\xi_n(\cdot}{\RR} < \infty$ and consequently from \eqref{E:total_var_xi}, $\totvar{\xi_n(\cdot)}{[0, T]} < \infty$ uniform in $n$ and $t$ and this completes the proof of Lemma \ref{T:positive_var}.
\end{proof}

By the result of Lemma \ref{T:positive_var} and consequently invoking (Helly's) compactness theorem, we conclude that 
\begin{equation}
    \xi_n \to \xi, \quad L^1([0, T], \RR).
\end{equation}
Using the continuity of bijection $u \mapsto \psi^{-1}(\gamma, u)$ for each $\gamma$, we can conclude that $\dot y_n \to \dot y$ in $L^1([0, T], \RR+)$. Here without loss of generality, we are assuming that $T > \tau_{M}$, where $\tau_m$ is as in \eqref{E:passage_time}. This completes the proof. 
\end{proof}

The next step is to show \eqref{E:BN_ODE} and \eqref{E:BN_cap} in the definition of the Cauchy solution will be satisfied by $(\rho, y)$. Fix $m \in \set{0, \cdots, M}$ and we prove the desired results in the interior $I_m^\circ$ of a fixed region $I_m$. We start by introducing some new notations and concepts.

By uniform convergence of $y_n$ to $y$, for a fixed $\eps_\circ >0$, for sufficiently large $n$ we have that $\sup_{t \in [0, \tau_{m+1}]} \abs{y_n (t) - y(t)} < \eps_\circ$ (in other words, $y_n \in \ball_{\norm{\cdot}_{L^\infty([0, \tau_{m+1}])}}(y , \eps_\circ)$) where $\tau_{m}$ is defined in \eqref{E:passage_time}. Therefore, we can consider the trajectories on a set $I_m^{\eps_\circ}$ defined by 
    \begin{equation}
        \label{E:approx_set}
        I_m^{\eps_\circ} \Def \set{x \in I_m^\circ: dist(x, \partial I_m) \ge \eps_\circ}.
    \end{equation}
    In particular, we will study the trajectory of $y$ over $I_m^{\eps_\circ}$ which ensures that both integral curves of $y_n$ and $y$ remain in $I_m^\circ$. In addition, we define the first and the last time that the trajectory of $y$ hits the boundary of $I_m^{\eps_\circ}$. More precisely, we initially investigate the results (see \eqref{E:BN_ineq} and \eqref{E:BN_limit_traj}) on $(\mathcal T_{\eps_\circ}^e, \mathcal T_{\eps_\circ}^o)$, where
    \begin{equation}
        \label{E:exit_time_interior}
        \mathcal T_{\eps_\circ}^e \Def \inf \set{t >0: y(t) \in \partial I_m^{\eps_\circ}}, \quad  \mathcal T_{\eps_\circ}^o \Def \sup \set{t >0: y(t) \in \partial I_m^{\eps_\circ}}
    \end{equation}
    and extend the results to the $I_m^\circ$ by letting $\eps_\circ \to 0$ as $n \to \infty$. 
Now, we are ready to prove the next main result.
\begin{theorem}\label{T:flux_constraint}
For a.e. $t \in \RR_+$ we have 
\begin{equation}\label{E:BN_ineq}
    \lim_{x \to y(t) \pm} f(\gamma(x), \rho(t,x)) - \dot y(t) \rho(t,x) \le F_\alpha(y(t), \dot y(t))
\end{equation}
\end{theorem}
\begin{proof}
Let's fix $m \in \set {0, \cdots, M}$. It is sufficient to prove the claimed result in the domain $I_m^\circ$. We define
\begin{equation*}\begin{split}
    &\Omega_n^- \Def \set{(t, x)\in (\mathcal T_{\eps_\circ}^e, \mathcal T_{\eps_\circ}^o) \times I_m^\circ: x < y_n(t)}, \quad
    \Omega^- \Def \set{(t, x)\in (\mathcal T_{\eps_\circ}^e, \mathcal T_{\eps_\circ}^o) \times I_m^\circ: x < y(t)} 
\end{split}\end{equation*}
where $\mathcal T_{\eps_\circ}^e$ and $\mathcal T_{\eps_\circ}^o$ are as in \eqref{E:exit_time_interior}. Let $\varphi \in C_c^\infty( (\mathcal T_{\eps_\circ}^e, \mathcal T_{\eps_\circ}^o)\times I_m^{\eps_\circ}; \RR)$ and $\varphi \ge 0$. Therefore, $\supp(\varphi(t, \cdot)) \subseteq [\underline \alpha, \bar \alpha] \subset [\eps_\circ, \sa_{m+1} - \eps_\circ]$. Now, we define 

For simplicity of the notations and without loss of generality, we assume the boundary location $\sa_m = y_\circ = 0$. 
We consider a function $\beta \in C_c^\infty([-1, 0]; [0, 1])$ such that $\beta \ge 0$, strictly increasing and $\int \beta(u)du = 1$. For $\eps >0$ we define  
\begin{equation*}
    \delta_\eps(x) \Def \tfrac 1\eps \beta(\tfrac{x}{\eps}), \quad \alpha_\eps(x) \Def \int_{-\infty}^{x} \delta_\eps(u) du
\end{equation*}
In particular, as $\eps \to 0$, $\delta_\eps \to \delta$ (Dirac distribution) and $\alpha_\eps(x) \to H(x)$ (Heaviside function) in the distribution sense. For sufficiently small $\eps < \eps_\circ$ such that $\underline \alpha - \eps \ge \eps_\circ$, we define
\begin{equation}\label{E:defn_test}
    \varphi_\eps(t, x) \Def \set{\alpha_\eps(x - \underline \alpha) - \alpha_\eps(x - y_n(t)) }^+
\end{equation}
where, $\set{u}^+ \Def \max \set{0, u}$. It can be noted that 
\begin{equation}
    \label{E:supp_varphi_eps}
    \supp(\varphi_\eps) \subset [\mathcal T_{\eps_\circ}^e, \infty) \times [\underline \alpha - \eps , y_n(t)]
\end{equation}
In addition, we define a compactly supported function
\begin{equation*}
     \phi_\eps(t, x) \Def \varphi_\eps(t,x) \varphi(t,x) \in W^{1,1}_\circ((\mathcal T_{\eps_\circ}^e, \mathcal T_{\eps_\circ}^o) \times I_m^\circ; \RR_+),
\end{equation*}
where for an open 
set $\Omega \subset \RR$, $W^{1,1}_\circ(\Omega) \Def \overline{C_c^\infty(\Omega)}^{\mathcal T_{\norm{\cdot}_{W^{1,1}(\Omega)}}}$, i.e. the closure of $C_c^\infty(\Omega)$ with respect to the topology of the Sobolev space $W^{1,1}(\Omega)$. 

\noindent In addition, considering that $\rho^{(n)}$ is a weak solution, we have 
\begin{multline}
    \label{E:weak_form}
    \int_{\RR_+} \int_\RR \left(\rho^{(n)}(t,x) \partial_t \phi_\eps(t,x) + f(\gamma(x), \rho^{(n)}(t, x)) \partial_x \phi_\eps \right) dx dt\\
    + \int_\RR \phi_\eps(0, x) \rho^{(n)}_\circ(x) dx= 0.
\end{multline}
The second integral vanishes due to the support of $\phi_\eps$. Let's define 
\begin{equation*}
    \mathfrak T \Def \set{t \in [\mathcal T_{\eps_\circ}^e,\mathcal T_{\eps_\circ}^o]: y_n(t) >\underline \alpha}, \quad \mathcal U_{\mathfrak T} \Def \mathfrak T \times \RR
\end{equation*}
With a slight abuse of notation, we consider $\mathfrak T$ and $\mathcal U_{\mathfrak T}$ for $y(t)$ instead of $y_n(t)$ when dealing with the solution $(\rho, y)$ and domain $\Omega^-$. 

Replacing $\phi_\eps$, on the set $\mathcal U_{\mathfrak T}$, we have that 
\begin{multline}\label{E:extended}
     \int_{\RR_+} \int_\RR \rho^{(n)}(t, x) \dot y_n(t) \delta_\eps(x - y_n(t))  \varphi(t,x) dx dt 
      +\int_{\RR_+} \int_\RR \rho^{(n)}(t, x) \varphi_\eps(t,x) \partial_t \varphi(t, x) dx dt \\
      + \int_{\RR_+} \int_\RR f(\gamma(x), \rho^{(n)}(t, x)) (\delta_\eps(x - \underline \alpha) - \delta_\eps(x - y_n(t))) \varphi(t, x) dx dt \\
      + \int_{\RR_+} \int_\RR f(\gamma(x), \rho^{(n)}(t, x)) \varphi_\eps(t, x) \partial_x\varphi(t, x) dx dt = 0.
\end{multline}
Remarking that $\supp(\delta_\eps(\cdot - y_n(t))) = cl(y_n(t)-\eps, y_n(t))$ and letting $\eps \to 0$, the first and second integral in \eqref{E:extended} read
\begin{equation}\label{E:lim_part1}
    \mathcal I = \int_{[\mathcal T_{\eps_\circ}^e, \mathcal T_{\eps_\circ}^o]\cap \mathfrak T} \rho^{(n)}(t, y_n(t)-) \dot y_n(t)\varphi(t, y_n(t))dxdt+ \int_{\Omega_n^- \cap \mathcal U_{\mathfrak T}} \rho^{(n)}(t, x) \partial_t \varphi(t, x)dxdt
\end{equation}
Similarly, for the last two integrals of \eqref{E:extended}, on the set $\mathcal U_{\mathfrak T}$, we can write
\begin{equation}\label{E:lim_part2}
    \begin{split}
        \mathcal J & = \int_{\RR_+}f(\gamma(x), \rho^{(n)}(t, x))  (\delta_{\underline\alpha} - \delta_{y_n(t)})(x) \varphi(t, x) dxdt\\
        & \quad + \int_{\RR_+}\int_\RR f(\gamma(x), \rho^{(n)}(t, x))  \bOne_{[\underline\alpha, y_n(t))} \partial_x \varphi(t, x)dxdt \\
        & = -\int_{[\mathcal T_{\eps_\circ}^e, \mathcal T_{\eps_\circ}^o] \cap \mathfrak T} f(\gamma_{r_m}, \rho^{(n)}(t, y_n(t)-)\varphi(t,y_n(t)-) dxdt+ \int_{\Omega_n^- \cap \mathcal U_{\mathfrak T}} f(\gamma(x), \rho^{(n)}(t, x))\partial_x\varphi(t,x)dx dt
    \end{split}
\end{equation}
where, we are using the fact that $\varphi(t, \underline \alpha) = 0$.
Therefore, replacing \eqref{E:lim_part1} and \eqref{E:lim_part2} in \eqref{E:extended}, we conclude that  
\begin{multline*}
    \int_{\Omega_n^- \cap \mathcal U_{\mathfrak T}} \left(\rho^{(n)}\partial_t \varphi + f(\gamma(x), \rho^{(n)}) \partial_x \varphi \right) dx dt \\= \int_{[\mathcal T_{\eps_\circ}^e, \mathcal T_{\eps_\circ}^o] \cap \mathfrak T} \left(f(\gamma_{r_m},\rho^{(n)}(t, y_n(t)-)) - \rho^{(n)}(t, y_n(t)-) \dot y_n(t) \right) \varphi(t, y_n(t)) dt
\end{multline*}
By Theorem \ref{T:convergence_rho_y} $(\rho, y)$ is a weak solution and hence following the same approach we will have that 
\begin{equation*}
    \int_{\Omega^- \cap \mathcal U_{\mathfrak T}} \left(\rho \partial_t \varphi + f(\gamma(x), \rho) \partial_x \varphi \right) dx dt= \int_{[\mathcal T_{\eps_\circ}^e, \mathcal T_{\eps_\circ}^o] \cap \mathfrak T} \left(f(\gamma_{r_m}, \rho(t, y(t)-)) - \rho(t, y(t)-) \dot y(t) \right) \varphi(t, y(t)) dt
\end{equation*}
In addition, since $\varphi \ge 0$ we have that  
\begin{multline*}
    \int_{[\mathcal T_{\eps_\circ}^e, \mathcal T_{\eps_\circ}^o] \cap \mathfrak T} \left(f(\gamma_{r_m}, \rho^{(n)}(t, y_n(t)-)) - \rho^{(n)}(t, y_n(t)-) \dot y_n(t) \right) \varphi(t, y_n(t)) dt \\
    \le \int_{[\mathcal T_{\eps_\circ}^e, \mathcal T_{\eps_\circ}^o] \cap \mathfrak T} F_\alpha(y_n(t), \dot y_n(t)) \varphi(t, y_n(t))dt
\end{multline*}
Collecting all together, on the set $\mathfrak T$ we can write
\begin{equation}\label{E:final}
    \begin{split}
        \int_{[\mathcal T_{\eps_\circ}^e, \mathcal T_{\eps_\circ}^o]} \left(f(\gamma_{r_m}, \rho(t, y(t)-)) - \rho(t, y(t)-) \dot y(t) \right) & \varphi(t, y(t)) dt  \\
        &= \int_{\Omega^-} \left(\rho \partial_t \varphi + f(\gamma(x), \rho \partial_x \varphi \right) dx dt\\
        & = \lim_{n \to \infty} \int_{\Omega_n^-} \left(\rho^{(n)} \partial_t \varphi + f(\gamma(x), \rho^{(n)}) \partial_x \varphi\right) dx dt \\
        & \le \lim_{n \to \infty} \int_{[\mathcal T_{\eps_\circ}^e, \mathcal T_{\eps_\circ}^o]} F_\alpha(y_n(t), \dot y_n(t)) \varphi(t, y_n(t)) dt \\
        & = \int_{[\mathcal T_{\eps_\circ}^e, \mathcal T_{\eps_\circ}^o]} F_\alpha(y(t), \dot y(t)) \varphi(t, y(t)) dt,
    \end{split}
\end{equation}
where the last equality is by the Theorem \ref{T:doty_n}. On the other hand, \eqref{E:final} is satisfied on $[\mathcal T_{\eps_\circ}^e, \mathcal T_{\eps_\circ}^o] \setminus \mathfrak T$ as $\varphi (t, x)$ vanishes on this set. This means the \eqref{E:final} holds true on $[\mathcal T_{\eps_\circ}^e, \mathcal T_{\eps_\circ}^o]$.

Furthermore, $\varphi$ and $\eps_\circ$ are chosen arbitrarily, the result can be concluded on the region $I_m^\circ$. 
In addition, by appropriately considering the positive values of function $\beta$ and $\bar \alpha$, the same approach can be applied to $\Omega_n^+$. This completes the proof.
\end{proof}
\begin{theorem}\label{T:BN_sol}
We have that for $a.e. \, t \in \RR_+$ 
\begin{equation}\label{E:BN_limit_traj}
    y(t) = y_\circ + \int_0^t w(y(s), \rho(s, y(s)+) ds
\end{equation}
where $\rho$ and $y$ are defined as in Theorem \ref{T:convergence_rho_y}. 
\end{theorem}
\begin{proof}
   It is sufficient to prove the result for $y(t) \in I_m^\circ$, for a fixed $m \in \set{0, \cdots, M}$.
 We mostly adopt the same approach as in \cite[Section 3.3]{liard2021entropic} for proving this theorem. However, due to \textit{lack of bounded variation and the possibility of waves crossing the boundaries}, we need to redefine and reprove some of the results. 

In particular, we reprove the claim of \cite[Lemma 4]{liard2021entropic} in which the proof depends on the assumption of the total variation of the solution $\rho(t, \cdot)$ (and is not directly applicable to the case of the present paper). 

Let $\eps >0$. From Theorems \ref{T:convergence_rho_y} and \ref{T:doty_n} there exists a measure zero set $\mathcal N$ such that for any $\bar t \in (\mathcal T_{\eps_\circ}^e, \mathcal T_{\eps_\circ}^o) \setminus \mathcal N$, we have that 
\begin{itemize}
    \item $\lim_{n \to \infty} \rho^{(n)}(\bar t, x) = \rho(\bar t, x)$ pointwise (up to some subsequence) for almost every $x \in \RR$, 
    \item $y(\cdot)$ is differentiable at $t = \bar t$, 
    \item $\lim_{n \to \infty}  \dot y_n(\bar t) = \dot y(t)$ pointwise (up to a subsequence) and in particular $\abs{y_n(\bar t) - y(\bar t)} \le \bar t \eps$ for sufficiently large $n$, 
    \item For all $n \in \NN$, $\dot y_n(\bar t) = \min\set{V_b^{(m)}, v(\gamma_{r_m}, \rho^{(n)}(\bar t, y_n(\bar t) + ))}$
\end{itemize}
For simplicity of the notation, let $\rho_- \Def \rho(\bar t, y(t) - )$ and $\rho_+ \Def \rho(\bar t, y(t) + )$. The following lemma provides a range for $\rho$ which is essential for the rest of the results. In what follows $\ball^-(x, r) \Def (x - r, x)$, $\ball^+(x, r) \Def (x , x + r)$ and $\ball(x, r) \Def (x - r, x + r)$ for some $x \in \RR$ and $r >0$. 
\begin{lemma}\label{T:ball_bound}
    Fix $\bar t \in (\mathcal T_{\eps_\circ}^e, \mathcal T_{\eps_\circ}^o) \setminus \mathcal N$, and $\eps >0$. For $\rho_-  , \rho_+ \in [0, 1]$, there exists a $\varsigma >0$ such that 
    \begin{equation*}
        \rho(\bar t, x) \in \begin{cases}
            \left(\max\set{\rho_- - \frac{\eps}{2}, 0}, \min\set{\rho_- + \frac{\eps}{2} , 1}\right) &, x \in \ball^-(y(\bar t), \varsigma) \\
            \left(\max\set{\rho_+ - \frac{\eps}{2}, 0}, \min\set{\rho_+ + \frac{\eps}{2} , 1}\right) &, x \in \ball^+(y(\bar t), \varsigma) 
        \end{cases}
    \end{equation*}
\end{lemma}
\begin{proof}
    Let's define $z_- \Def \lim_{x \nearrow y(\bar t)} z(\bar t, x)$ and $z_+ \Def \lim_{x \searrow y(\bar t)} z(\bar t, x)$. First, we note that the map $z \mapsto \psi^{-1}(\gamma_{r_m}, z)$ is uniformly continuous on $[- \frac 14 \gamma_{r_m}, \frac 14 \gamma_{r_m}]$. Therefore, 
    \begin{equation*}
        \rho_- = \lim_{x \nearrow y(\bar t)} \rho(\bar t, x) = \lim_{x \nearrow y(\bar t)}\psi^{-1}(\gamma_{r_m}, z(\bar t, x)) = \psi^{-1}(\gamma_{r_m} , z_-)
    \end{equation*}
    Similarly, we can show that $\rho_+ = \psi^{-1}(\gamma_{r_m}, z_+)$. 
    The uniform continuity of $z \mapsto \psi^{-1}(\gamma_{r_m} , z)$ then asserts that for the given $\eps$ there exists an $\bar \eps >0$ sufficiently small such that 
    \begin{equation} \label{E:last_lemma}
    z \in \ball(z', \bar \eps) \subset [- \tfrac{\gamma_{r_m}}{4}, \tfrac{\gamma_{r_m}}{4}], \quad \text{for any $z'$, implies $\rho = \psi^{-1}(\gamma_{r_m}, z) \in \ball(\rho', \eps)$}\end{equation} 
    where $\rho' = \psi^{-1}(\gamma_{r_m}, z')$. Now, for $\bar \eps$ and by employing the fact that $\totvar{z(\bar t, \cdot)}{\RR} < \infty$, we choose $\varsigma>0$ sufficiently small such that
    \begin{equation*}
        \totvar{z(\bar t, \cdot)}{(y(\bar t) - \varsigma, y(\bar t))} < \frac{\bar \eps}{2}
    \end{equation*}
    and consequently, 
    \begin{equation*}
        z(\bar t, x) \in \ball(z_-, \frac{\bar \eps}{2}) , \quad \text{for $x \in \ball^-(y(\bar t), \varsigma)$}.
    \end{equation*}
    Therefore, using the fact that $z \mapsto \psi^{-1}(\gamma_{r_m}, z)$ is monotonically increasing (and hence preserves the order) and \eqref{E:last_lemma}, we have that 
    \begin{equation*}
        \rho(\bar t, x) \in \ball(\rho_-, \frac{\eps}{2}), \quad \text{for $x \in \ball^-(y(\bar t), \varsigma)$}.
    \end{equation*}
    A similar argument proves the claimed range of $\rho(\bar t, \cdot)$ over $\ball^+(y(\bar t), \varsigma)$. This concludes the proof.
\end{proof}
Employing the bound derived in Lemma \ref{T:ball_bound}, the rest of the proof of Theorem \ref{T:BN_sol} follows directly from Appendix \ref{A:min_distance} and the proofs in \cite[Section 3.3]{liard2021entropic}. In fact, since all the waves entering the $I_m^{\eps_\circ}$ region are either shocks or rarefaction, all the analytical results of \cite[Section 3.3]{liard2021entropic} follow immediately. Finally, letting $\eps_\circ \to 0$, the results remain valid on $I_m^\circ$. This completes the proof of the Theorem \ref{T:BN_sol}.  
\end{proof}
\noindent Now, we have all the results to prove the main theorem.
\begin{proof}[Proof of Theorem \ref{T:existence_Cauchy}]
Putting the results of Theorems \ref{T:convergence_rho_y}, \ref{T:flux_constraint} and \ref{T:BN_sol} together, the claim follows.
\end{proof}
\newpage
\appendix
\section{Proof of Lemma \ref{T:phase_1}}\label{A:phase_1}

Let $\eta_\eps$ be the smooth approximation of the step function $\bOne_{[s,t]}$. In particular, we define $u \in \RR \mapsto \eta_\eps(u)$ by
\begin{equation*}
    \eta_\eps \Def \lambda_\eps \star \bOne_{[s, t]}
\end{equation*}
where $\star$ is the mathematical convolution, 
\begin{equation*}
    \lambda_\eps(u) = \tfrac{1}{\eps } \lambda(\tfrac u\eps), \quad \lambda(u) = \begin{cases} c \exp \lb\tfrac{1}{u^2 -1}\rb &, \abs{u} <1 \\
0 &, \abs{u} \ge 1 \end{cases}
\end{equation*}
and constant $c$ is defined such that the integral $\int_\RR \lambda(u) du = 1$. By definition, for sufficiently small $\eps>0$, we have that
\begin{equation}\label{E:eta_eps}
    \eta_\eps \in C_c^\infty(\RR), \quad \supp(\eta_\eps)\subset \supp(\lambda_\eps) + \supp(\bOne_{[s, t]}) \subset [s- \eps , t + \eps] \subsetneq [0, \tau).
\end{equation}
We are interested in properties of $\eta_\eps$. Let $\mathcal D'(\RR)$ denote the space of distributions (the continuous dual of topological vector space $C_c^\infty(\RR)$). We have that
\begin{equation}\label{E:eta_eps_convergence}
    \lim_{\eps \to 0} \eta_\eps = \lim_{\eps \to 0} \lambda_\eps \star \bOne_{[s,t]} = \bOne_{[s, t]} \in \mathcal D'(\RR)
\end{equation}
where the convergence is in the distributional sense (the convergence can be shown to be valid pointwise and in $L^p$ as well). In other words, $\lim_{\eps \to 0} \lambda_\eps(r) = \delta(r)$ where the limit is defined with respect to the strong topology on the space of distributions. 

Now, let's define
\begin{equation*}
    \varphi_\eps(r, x) \Def \eta_\eps(r) \phi(x)
\end{equation*}
for a function $\phi \in C_c^\infty(\RR)$. Therefore, by \eqref{E:eta_eps} we should have $\varphi_\eps \in C_c^\infty(\RR\times \RR)$ with $\supp(\varphi_\eps) \subset [0, \tau) \times \supp(\phi)$. Let's define a function $t \mapsto q^{(n)}(t, \cdot)$ which extends $t \mapsto \rho^{(n)}(t, \cdot)$ to zero for $t \notin [0, \tau)$ (function $\rho(t, \cdot)$ is not defined beyond $\tau$ at this point). Noting that the solution $(\rho^{(n)}, y_n)$ is a well-defined over $[0, \tau)$ and $\rho^{(n)}$ is piecewise constant, we should have  
\begin{equation}\label{E:weak_eta_eps}
    \int_\RR \int_\RR \left(q^{(n)}(r, x) \fpartial{\varphi_\eps}{r} + f^{(n)}(\gamma(x),q^{(n)}(r, x)) \fpartial{\varphi_\eps}{x} \right) dx dr+ \int_\RR \varphi_\eps(0, x) q^{(n)}(0, x) dx = 0,
\end{equation}
First, we note that the last integral vanishes due to the definition of its support on $[0, \tau)$. Next, we note that 
\begin{equation*}
    \fpartial{\varphi_\eps}{r}(r, x) = \fpartial{\eta_\eps}{r}(r) \phi(x)
\end{equation*}
and in addition, we have that 
\begin{equation}\label{E:charac_deriv}
    \fpartial{\bOne_{[s,t]}}{r} = \delta_s - \delta_t,
\end{equation}
in the sense of distributions and $\delta_t$ is the Dirac distribution concentrated at $t$.  
Let's proceed with the first term in the first integral of \eqref{E:weak_eta_eps}. By changing the order of integrals (employing the Fubini theorem), for a.e. $x \in \RR$, the inside integral reads 
\begin{equation*}
    \begin{split}
        \int_{r \in \RR}  q^{(n)}(r, x) \fpartial{\eta_\eps(r)}{r} dr & =  \int_{r \in \RR} q^{(n)}(r, x) \left(\fpartial{\lambda_\eps}{r} \star \bOne_{[s,t]}\right)(r) dr\\
        & = \int_{r \in \RR } q^{(n)}(r, x) \int_{u \in \RR} \fpartial{\lambda_\eps(r - u) }{r} \bOne_{[s, t]}(u) du dr\\
        & \overset{\langle 1 \rangle}{=}- \int_{r \in \RR} q^{(n)}(r, x) \int_{u \in \RR}\fpartial{\lambda_\eps(u - r) }{u} \bOne_{[s, t]}(u) du dr \\
        & \overset{\langle 2 \rangle}{=}  \int_{r \in \RR} q^{(n)}(r, x) \left(\lambda_\eps(t - r) - \lambda_\eps(s -r) \right) dr\\
        & = (\lambda_\eps \star q^{(n)}(\cdot, x))(t) - (\lambda_\eps \star q^{(n)}(\cdot, x))(s)\\
    \end{split}
\end{equation*}
where equality $\langle 1 \rangle$ is by the fact that 
\begin{equation*}
    \frac{d}{dr}\lambda_\eps(r - \tau) = \exp \lb\frac{-1}{1 - \left(\frac{(r - \tau)^2}{\eps^2}\right)^2} \rb  \frac{- 2 (r - \tau)}{\eps^2 \left(\frac{(r - \tau)^2}{\eps^2}\right)^2} = - \frac{d}{d\tau } \lambda_\eps(\tau - r)
\end{equation*}
and the equality $\langle 2 \rangle $ is by using \eqref{E:charac_deriv}.
Therefore, for a.e. $x \in \RR$
\begin{equation*}
   \lim_{\eps \to 0} \int_{r \in \RR} q^{(n)}(r, x) \fpartial{\eta_\eps(r)}{r} dr = q^{(n)}(t, x) - q^{(n)}(s,x)
\end{equation*}
in the sense of distributions. 
This implies that by passing $\eps \to 0$ (using dominated convergence theorem), the first term of the first integral of \eqref{E:weak_eta_eps} will be
\begin{equation*}\begin{split}
  (I) & \Def \int_{x \in \RR} \phi(x) \left(q^{(n)}(t, x) - q^{(n)}(s, x) \right)dx = \int_{x \in \RR} \phi(x) \left(\rho^{(n)}(t, x) - \rho^{(n)}(s, x) \right)dx
\end{split}\end{equation*}
where the second equality follows from the definition of $q^{(n)}$ and the fact that $s, t \in [0, \tau)$. 
For the second term of the first integral of \eqref{E:weak_eta_eps}, letting $\eps \to 0$ and using the convergence \eqref{E:eta_eps_convergence} in the pointwise sense, we have 
\begin{equation*}
    (II) \Def \int_{r \in [s, t]} \int_{x \in \RR } f^{(n)}(\gamma(x) , q^{(n)}(r, x)) \phi'(x) dx dr = \int_{r \in [s, t]} \int_{x \in \RR } f^{(n)}(\gamma(x) , \rho^{(n)}(r, x)) \phi'(x) dx dr
\end{equation*}
where the second equality is using $r \in [s, t] \subset [0, \tau)$. Collecting all together, from \eqref{E:weak_eta_eps}, $(I)$ and $(II)$, we have 
\begin{equation*}
   \int_{x \in \RR} \phi(x) \left(\rho^{(n)}(t, x) - \rho^{(n)}(s, x) \right)dx + \int_{r \in [s, t]} \int_{x \in \RR } f^{(n)}(\gamma(x) , \rho^{(n)}(r, x)) \phi'(x) dx dr = 0,
\end{equation*}
which proves the claim.

\section{Proof of Lemma \ref{lem:TV_fn}}\label{A:TV_fn}
Let $x_1< x_2 <\cdots< x_{N_{\mathcal P}}$ be any finite sequence $\mathcal P$ of points in $\RR$. Then, by definition
\begin{equation*}\begin{split}
        &\totvar{f^{(n)}(\gamma(\cdot), \rho^{(n))}(t, \cdot)}{\RR} \\
        &= \sup_{\mathcal P} \sum_{k = 0}^{N_{\mathcal P}-1} \abs{f^{(n)}(\gamma(x_{k +1}), \rho^{(n)}(t, x_{k +1})) -f^{(n)}(\gamma(x_{k}), \rho^{(n)}(t, x_{k})) }\\
         &\le \sup_{\mathcal P} \sum_{k = 0}^{N_{\mathcal P}-1} \Bigg\{\abs{f^{(n)}(\gamma(x_{k +1}), \rho^{(n)}(t, x_{k +1})) -f^{(n)}(\gamma(x_{k+1}), \tfrac 12) } \\
         &\quad + \abs{f^{(n)}(\gamma(x_{k}), \rho^{(n)}(t, x_{k})) -f^{(n)}(\gamma(x_{k}), \tfrac 12) } \\
         &\quad + \abs{f^{(n)}(\gamma(x_{k +1}, \tfrac 12) - f^{(n)}(\gamma(x_{k}), \tfrac 12)} \Bigg\}\\
        &\le \sup_{\mathcal P} \sum_{k = 0}^{N_{\mathcal P}-1} \Bigg\{ \abs{\sign \left(\tfrac 12 - \rho^{(n)}(t, x_{k+1})\right)\left(f^{(n)}(\gamma(x_{k +1}), \rho^{(n)}(t, x_{k +1})) -f^{(n)}(\gamma(x_{k+1}), \tfrac 12)\right) } \\
          & \qquad + \abs{\sign \left(\tfrac 12 - \rho^{(n)}(t, x_{k}) \right)\left(f^{(n)}(\gamma(x_{k}), \rho^{(n)}(t, x_{k})) -f^{(n)}(\gamma(x_{k}), \tfrac 12)\right) } \\
          & \qquad + \tfrac 14\abs{\gamma(x_{k +1}) - \gamma(x_{k})} \Bigg\}\\
          & \le \totvar{{z^{(n)}(t, \cdot})}{\RR} + \tfrac 14 \totvar{\gamma(\cdot)}{\RR}. 
          \end{split}
\end{equation*}
This completes the proof of the lemma. 


\section{Minimum Distance Between $\rho_L$ and $\rho_R$}\label{A:min_distance}
One of the main concepts that play a major role in the proof of most of the results of \cite{liard2021entropic} is the minimum distance to go from a density $\rho_L$ to $\rho_R$ in the solution $\rho^{(n)}$. In the context of the present paper, this concept requires validation as waves can travel between regions. 
 
Let's fix $m \in \set{0, \cdots, M}$ and $\rho_L, \rho_R \in \mathcal G_m^{(n)}$ (defined as in Notation \ref{N:proj_FD_m}) with $\rho_L > \rho_R$ and a fixed time $\bar t \in (\mathcal T_{\eps_\circ}^e, \mathcal T_{\eps_\circ}^o)$. Then, a set $\mathcal A_m(\rho_1, \rho_2, \bar t) \subset \mathcal G_m^{(n)} \times I_m^{\eps_\circ} \times I_m^{\eps_\circ}$ is defined by
   \begin{equation*}
       (\rho_\circ^{(n)}, x_L, x_R) \in \mathcal A_m(\rho_1, \rho_2, \bar t) , \quad \text{if}
   \end{equation*}
  
   \begin{enumerate}
       \item[(i)] $ x_L < x_R$ with $\rho^{(n)}(\bar t, x_i) = \rho_i$ for $i \in \set{L, R}$, 
       \item[(ii)] Either, for all $x \in [x_L, x_R]$,  $\rho(\bar t, x-) - \rho^{(n)}(t, x+) = \zeta_+^{(n)} (\rho_R)$ 
       \item[(iii)] Or $\rho^{(n)}(\bar t, x-) \le \rho^{(n)}(\bar t, x+)$, 
   \end{enumerate}
   where $\rho^{(n)}(\bar t, \cdot)$ is the wave-front tracking solution at time $\bar t$ with initial data $\rho^{(n)}_\circ$ and $\zeta_+^{(n)}(\cdot) \Def \psi^{-1}(\gamma_{r_m}, \delta_+^{(n)}(\cdot))$, defined as of Notation \ref{N:PP_grid_points}. In other words, $\zeta_+^{(n)}(\cdot)$ denotes the distance to the next grid point in $\mathcal V$ space.

   In particular, conditions (ii) and (iii) assert that over $\set{\bar t} \times [x_L, x_R]$, shocks and rarefactions are the only admissible waves. The difference between this definition and the one in \cite{liard2021entropic} is here waves can originate from other regions. However, in a similar fashion, we claim that for any $\rho_L, \rho_R \in \mathcal G_m^{(n)}$ and such that $\rho_L > \rho_R$ and $\bar t \in (\mathcal T_{\eps_\circ}^e, \mathcal T_{\eps_\circ}^o)$, we have that 
   \begin{equation}
       \label{E:min_distance_waves}
       \varsigma^{(n)}(\rho_L, \rho_R, \bar t) \Def \inf_{(\rho_\circ^{(n)}, x_L, x_R) \in \mathcal A_m(\rho_L, \rho_R, \bar t)} (x_R - x_L) \ge 2\bar t (\rho_L - \rho_R - \bar \zeta^{(n)})
   \end{equation}
where, $\bar \zeta^{(n)} \Def \psi^{-1}(\gamma_{r_m},\bar \delta^{(n)})$ (cf. the Notation \ref{N:PP_grid_points}). More precisely, $\varsigma^{(n)}(\rho_L, \rho_R, \bar t)$ is the minimum distance in space at time $\bar t$ to go from $\rho_L$ to $\rho_R$ with respect to the solution $\rho^{(n)}$ only using shocks and rarefactions. 

\begin{figure}
    \centering
    \includegraphics[width=4in]{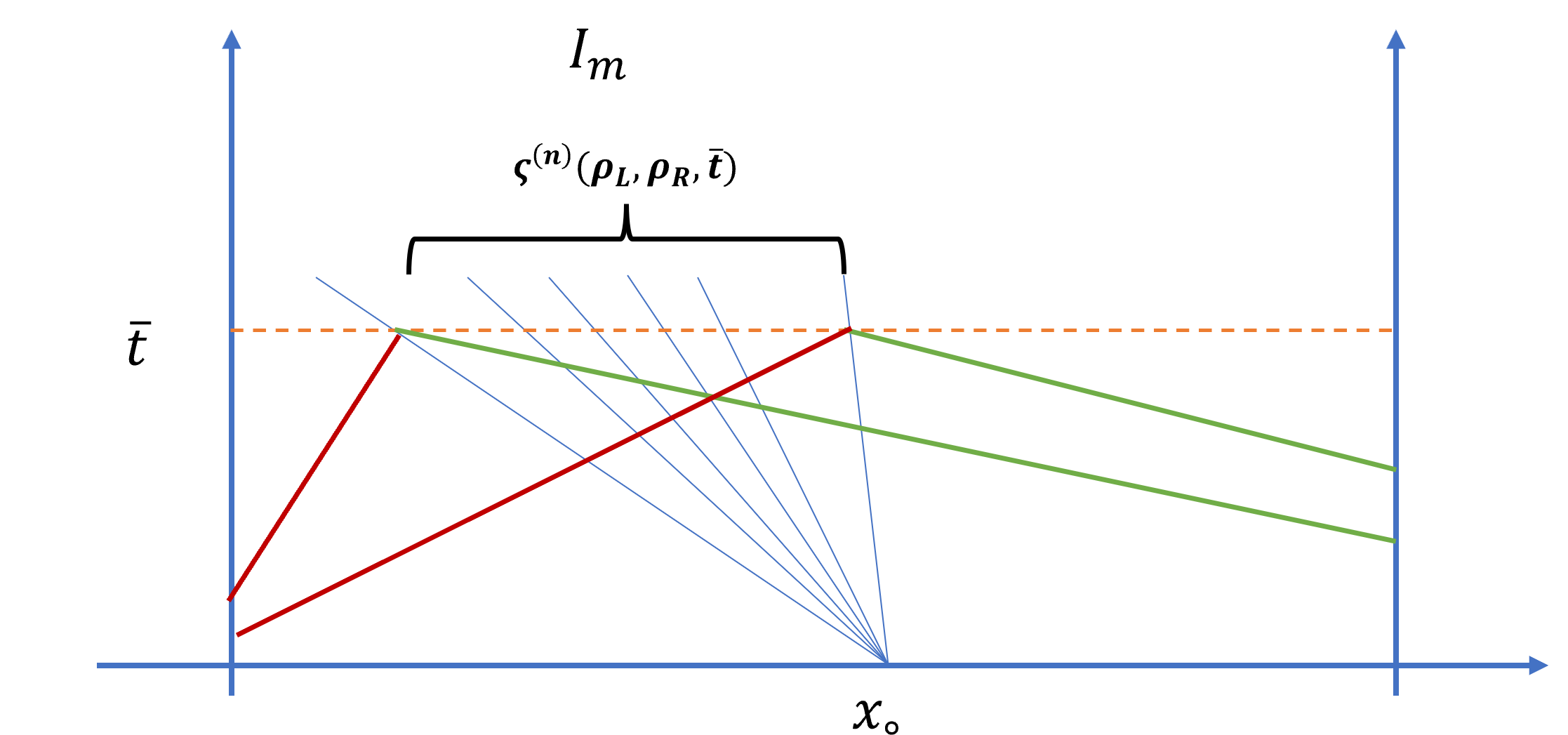}
    \caption{The minimum distance between the $x_L$ and $x_R$ such that $(\rho_\circ^{(n)}, x_L, x_R) \in \mathcal A_m(\rho_L, \rho_R, \bar t)$ can be presented by a rarefaction shock starting from a point $x_\circ \in I_m^\circ$. }
    \label{fig:min_distance_waves}
\end{figure}

As mentioned before, in the case of the problem in these notes, in addition to the waves in the $I_m$ region, we might have other waves crossing the boundaries of $I_m$. However, as can be noticed in Figure \ref{fig:min_distance_waves} even if $x_L$ and $x_R$ are created by waves from outside the region $I_m$, the minimum distance between these two points can be calculated by considering a rarefaction starting from any point $x_\circ \in I_m^{\eps_\circ}$. Therefore, by a similar calculation to the \cite[Lemma 1]{liard2021entropic}, we will have the claimed lower bound \eqref{E:min_distance_waves}.

\newpage
 \bibliographystyle{elsarticle-num} 
\bibliography{references}






\end{document}